\documentclass[12pt]{article}
\usepackage{amsmath,amssymb,amsthm,amscd,latexsym,}
\usepackage{graphicx}

\makeatletter
\renewcommand{\mod}[1]{\allowbreak \if@display \mkern 8mu \else
\mkern 5mu\fi {\operator@font mod}\,\,#1}
\makeatother

\newcommand{\bc}{\mathbb C}

\newcommand{\bz}{\mathbb Z}

\newcommand{\bp}{\mathbb P}

\newcommand{\aaa}{\mathbb A}
\newcommand{\ddd}{\mathbb D}

\DeclareMathOperator{\rk}{rk}

\newtheorem{theorem}{Theorem}

\newcommand\SSS{\mathfrak S}
\newcommand\AAA{\mathfrak A}

\begin{document}
\title{Classification of Picard lattices of K3 surfaces \footnote{With support by Russian
Sientific Fund N 14-50-00005.}}
\date{5 September 2017}
\author{Viacheslav V. Nikulin}
\maketitle

\begin{abstract}
Using results of our papers  \cite{Nik9}, \cite{Nik10} and \cite{Nik11} about
classification of
degenerations of K\"ahlerian K3 surfaces with finite
symplectic automorphism groups, we classify Picard lattices
of K\"ahlerian K3 surfaces. By classification we understand
classification depending on their possible finite symplectic automorphism groups
and their non-singular rational curves if a Picard lattice is negative
definite.
\end{abstract}

\centerline{Dedicated to the memory of Igor Rostislavovich Shafarevich}

\section{Introduction}
\label{sec1:introduction}

Let $X$ be a K\"ahlerian K3 surface with negative definite Picard lattice $S_X$.
Then the full symplectic automorphism group $G=Aut\ (X)_0$ of $X$ is finite and
and its coinvariant sublattice $S_G=((S_X)^G)^\perp_{S_X}$ is negative definite.
We denote by $P(X)\subset S_X$ the set of classes of all irreducible
non-singular rational curves
of $X$. They have the square $(-2)$. The primitive sublattice $S=MS_X=[S_G,P(X)]_{pr}\subset S_X$ of $S_X$
generated by $S_G$ and $P(X)$ gives the {\it main part} of the Picard lattice $S_X$.
The remaining part $(S)^\perp_{S_X}$ can be considered as {\it exotic part} of
the Picard lattice $S_X$.

If a K3 surface $X$ is algebraic, instead of $S_X$, one can take the
orthogonal complement $(h)^\perp_{S_X}$ to a nef element $h\in S_X$ with $h^2>0$.

We classify the main parts $S$ of Picard lattices of K3 if the group $G$ is large enough, it has order
$|G|\ge 8$, equivalently, it is different from $D_6$, $C_4$, $(C_2)^2$, $C_3$, $C_2$ and $C_1$
where $C_m$ is a cyclic group of order $m$ and $D_m$ is a dihedral group of order $m$.
It was known (for all $G$) if the set $P(X)$ is empty. See \cite{Hash}, \cite{Kon}, \cite{Muk},
\cite{Nik0},  \cite{Xiao}. The case, when $P(X)$ is non-empty, can be considered
as a degeneration of K\"ahlerian K3 surfaces with the finite symplectic automorphism group $G$.
Degenerations of K\"ahlerian K3 surfaces with finite symplectic automorphism groups
were studied in our papers \cite{Nik9}, \cite{Nik10} and \cite{Nik11}.

In Section \ref{sec2:degen}, we remind to a reader our results of
\cite{Nik9}, \cite{Nik10} and \cite{Nik11} about classification of degenerations
of K\"ahlerian K3 surfaces with finite symplectic automorphism groups.

In Section \ref{sec3:piclatt}, we present their applications
to classification of Picard lattices of K3 surfaces.

In Section \ref{sec4:tables}, we give Tables 1---4 which contain
our classification of Picard lattices of K3 surfaces: all possible
cases are given in lines of the Tables 1---4 which are not marked
by $o$. Lines which are marked by $o$ are also important. They give
classification of cases when a degeneration of K3 with finite
symplectic automorphism group $G$ has larger symplectic automorphism
group $\tilde{G}$ which has $|\tilde{G}|>|G|$. The same happens for Kummer surfaces.
They can be considered as a degeneration with the Dynkin diagram $16\aaa_1$ of
K3 with trivial symplectic automorphism group.
But, their symplectic automorphism group is $(C_2)^4$. It
was shown in our paper \cite{Nik-1}.

We hope to consider similar classification for remaining groups
$D_6$, $C_4$, $(C_2)^2$, $C_3$, $C_2$ and $C_1$,
and more details in further variants of the paper and further publications.

%%%%%%%%%%%%%%%%%%%%%%%%%%%%%%%%%%%%%%%%%%%%%%%%%%%%%%%%%
%%%%%%%%%%%%%%%%%%%%%%%%%%%%%%%%%%%%%%%%%%%%%%%%%%%%%%%%%
%%%%%%%%%%%%%%%%%%%%%%%%%%%%%%%%%%%%%%%%%%%%%%%%%%%%%%%%%

\section{On classification of degenerations of\\
K\"ahlerian K3 surfaces with finite symplectic\\
automorphism groups. A reminder.}
\label{sec2:degen}

Here we remind to a reader our results from \cite{Nik9}, \cite{Nik10} and \cite{Nik11}
about classification of degenerations of K\"ahlerian K3 surfaces
with finite symplectic automorphism groups. We refer to these papers for more details.

Let $X$ be a K\"ahlerian K3 surface (e. g. see \cite{Sh}, \cite{PS},
\cite{BR}, \cite{Siu}, \cite{Tod}
about such surfaces). That is $X$ is a non-singular compact
complex surface with the trivial canonical class $K_X$, and
its irregularity $q(X)$ is equal to 0.
Then $H^2(X,\bz)$ with the intersection pairing
is an even unimodular lattice $L_{K3}$
of the signature $(3,19)$. For a non-zero holomorphic $2$-form $\omega_X\in \Omega^2[X]$,
we have $H^{2,0}(X)=\Omega^2[X]=\bc\omega_X$. The primitive sublattice
$$
S_X=H^2(X,\bz)\cap H^{1,1}(X)=\{x\in H^2(X,\bz)\ |\ x\cdot H^{2,0}(X)=0 \ \}\subset H^2(X,\bz)
$$
is the {\it Picard lattice} of $X$
generated by first Chern classes of all line bundles over $X$.
We remind that a primitive sublattice means that $H^2(X,\bz)/S_X$ has no torsion.

Let $G$ be a finite symplectic automorphism group of $X$.
Here symplectic means
that for any $g\in G$, for a non-zero holomorphic $2$-form $\omega_X\in
H^{2,0}(X)=\Omega^2[X]=\bc\omega_X$, one has $g^\ast(\omega_X)=\omega_X$.
For an $G$-invariant sublattice $M\subset H^2(X,\bz)$, we denote by
$M^G=\{x\in M\ |\ G(x)=x\}$ the {\it fixed sublattice of $M$,} and by
$M_G=(M^G)^\perp_M$ the {\it coinvariant sublattice of $M$.}
By \cite{Nik-1/2}, \cite{Nik0}, the coinvariant lattice $S_G=H^2(X,\bz)_G=(S_X)_G$ is
{\it Leech type lattice:} i. e.  it is negative definite,
it has no elements with square $(-2)$, $G$ acts trivially
on the discriminant group $A_{S_G}=(S_G)^\ast/S_G$, and $(S_G)^G=\{0\}$.
For a general pair $(X,G)$, the $S_G=S_X$, and non-general $(X,G)$
can be considered as K\"ahlerian K3 surfaces
with the condition $S_G\subset S_X$
on the Picard lattice (in terminology of \cite{Nik0}).
The dimension of their
moduli is equal to $20-\rk S_G$.

Let $E\subset X$ be a non-singular irreducible rational curve
(that is $E\cong \bp^1$).
It is equivalent to: $\alpha=cl(E)\in S_X$, $\alpha^2=-2$,
$\alpha$ is effective
and $\alpha$ is numerically effective: $\alpha\cdot D\ge 0$
for every irreducible curve
$D$ on $X$ such that $cl(D)\not=\alpha$.

Let us consider the primitive sublattice $S=[S_G, G(\alpha)]_{pr}\subset S_X$
of $S_X$ generated by the coinvariant sublattice $S_G$ and
all classes of the orbit $G(E)$.
Since $S_G$ has no elements with square $(-2)$, it follows
that $\rk S=\rk S_G+1$
and $S=[S_G,\alpha]_{pr}\subset S_X$.

Let us {\it assume that the lattice $S=[S_G,\alpha]_{pr}$
is negative definite.} Then the
elements $G(\alpha)$ define the basis of the root system $\Delta(S)$ of
all elements with square $(-2)$ of $S$. All curves $G(E)$ of $X$
can be contracted to Du Val singularities of types of
connected components of the Dynkin diagram of the basis. The group $G$
will act on the corresponding singular K3 surface
$\overline{X}$ with these Du Val
singularities. For a general such triplet $(X,G,G(E))$,
the Picard lattice $S_X=S$,
and such triplets can be considered as {\it a degeneration
of codimension $1$} of
K\"ahlerian K3 surfaces $(X,G)$ with the finite symplectic
automorphism group $G$.
Really, the dimension of moduli of K\"ahlerian K3 surfaces
with the condition $S\subset S_X$
on the Picard lattice is equal to $20-\rk S=20-\rk S_G-1$.

By Global Torelli Theorem for K3 surfaces \cite{PS},
\cite{BR}, the main invariants
of the degeneration is the {\it type of the abstract group $G$}
which is equivalent
to the isomorphism class of the coinvariant lattice $S_G$,
and the type of the degeneration which is
equivalent to the Dynkin diagram of the basis $G(\alpha)$
or the Dynkin diagram
of the rational curves $G(E)$.

We can consider only the maximal finite symplectic
automorphism group $G$ with the
same coinvariant lattice $S_G$, that is $G=Clos(G)$.
By Global Torelly Theorem for K3 surfaces, this
is equivalent to
$$
G|S_G=\{ g\in O(S_G)\ |\ g\ is\ identity\ on\ A_{S_G}=(S_G)^\ast/S_G \}.
$$
Indeed, $G$ and $Clos(G)$ have the same lattice $S_G$, the same orbits
$G(E)$ and $Clos(G)(E)$, and the same sublattice $S\subset S_X$.

In \cite{Nik9}, all types of $G=Clos(G)$ and types of degenerations
(that is Dynkin
diagrams of the orbits $G(E)$) are described. They are described
in Table 1 of Section \ref{sec4:tables} which
is the same as \cite[Table 1]{Nik10}
where ${\bf n}$ gives types of possible $G=Clos(G)$, and we show all possible
types of degenerations at the corresponding rows by their Dynkin diagrams.

In the Table 1 of Section \ref{sec4:tables},
the type of $G=Clos(G) $ and the isomorphism class of the lattice $S_G$
is marked by ${\bf n}$. We also give the genus of $S_G$ which
is defined by
the discriminant quadratic form $q_{S_G}$. They
were calculated in papers \cite{Nik0},
\cite{Muk}, \cite{Xiao}, \cite{Hash}.

In \cite{Nik10}, with two exceptions which are marked by $I$ and $II$
for ${\bf n}=10$ (group $D_8$) and the degeneration $2\aaa_1$,
and for ${\bf n}=34$ (group $\SSS_4$) and the degeneration $6\aaa_1$, it was
shown that the lattice $S=[S_G,\alpha]_{pr}$ is unique,
up to isomorphisms, and its
genus (equivalent to $\rk S$
and the discriminant quadratic form $q_S$) is shown in Table 1.
In \cite[Table 2]{Nik10},
we described possible markings of $S$, $G$ and $\alpha$ by
Niemeier lattices which
give exact lattices descriptions of $S$, and the action of
$G$ on $S$ and $G(\alpha)$.

For degenerations of arbitrary codimension $t\ge 1$ which
were considered in \cite{Nik11},
instead of one orbit $G(E)$
of a non-singular rational curve, we should consider
$t\ge 1$ different orbits
$G(E_1),\dots, G(E_t)$ of non-singular rational curves on $X$, and their
classes $G(\alpha_1), \dots, G(\alpha_t)$ in $S_X$, but we also assume that
the sublattice $S=[S_G,G(\alpha_1),\dots,G(\alpha_t)]_{pr}\subset S_X$
is negative definite. Then the codimension of the degeneration is
equal to $t=\rk S-\rk S_G$, and $S_X=S=[S_G,\alpha_1,\dots, \alpha_t]_{pr}$
in general. We remark that $\rk S=\rk S_G+t\le 19$ since $H^2(X,\bz)$ has
the signature $(3,19)$. Thus, $t\le 19-\rk S_G$.

The type of the degeneration is given by the Dynkin
diagrams and subdiagrams
$$
(Dyn(G(\alpha_1)),\dots, Dyn(G(\alpha_t)))\subset
Dyn(G(\alpha_1)\cup\dots \cup G(\alpha_t))
$$
and their types. In difficult cases, we also
consider the matrix of subdiagrams
which is defined by
$$
(Dyn (G(\alpha_i)),\, Dyn ( G(\alpha_j)))\subset Dyn(G(\alpha_i)\cup G(\alpha_j))
$$
and their types for $1\le i<j\le t$.

In Table 2 of Section \ref{sec4:tables} which is similar to \cite[Table 2]{Nik11},
we give classification of types of degenerations of arbitrary
codimension $\ge 2$ for ${\bf n}\ge 12$ (then the
codimension is less or equal to $4$).
Equivalently, either $|G|>8$ or $G\cong Q_8$ of order $8$ ($G$ is big enough).
We calculate
genuses of the lattices $S$ by $\rk S$ and $q_S$. By $\ast$,
we mark cases when we prove that
the lattice $S$ is unique up to isomorphisms for the given type.
We mark by $o$ some cases (it will be  important for classification of
Picard lattices, and we shall discuss this later).
The description of their markings by Niemeier lattices is given in \cite[Table 3]{Nik11}.
For remaining groups $D_8$ (equivalently, ${\bf n}=10$) and $(C_2)^3$ (equivalently, ${\bf n}=9$)
of order $8$, similar results are given in Tables 3 and 4 of Section \ref{sec4:tables}
which are similar to tables of \cite{Nik11}. Their markings by Niemeier lattices are given in \cite[Table 5]{Nik11}
and \cite[Table 7]{Nik11} respectively.

For remaining finite symplectic automorphism groups $G=Clos(G)$ of order $|G|<8$, which are
$D_6$ (${\bf n}=6$), $C_4$ (${\bf n}=4)$, $(C_2)^2$ (${\bf n}=3$),
$C_3$ (${\bf n}=2$), $C_2$ (${\bf n}=1$) and $C_1$, similar classification results
are unknown yet.

%%%%%%%%%%%%%%%%%%%%%%%%%%%%%%%%%%%%%%%%%%%%%%%%%
%%%%%%%%%%%%%%%%%%%%%%%%%%%%%%%%%%%%%%%%%%%%%%%%%
%%%%%%%%%%%%%%%%%%%%%%%%%%%%%%%%%%%%%%%%%%%%%%%%%

%%%%%%%%%%%%%%%%%%%%%%%%%%%%%%%%%%%%%%%%%%
%%%%%%%%%%%%%%%%%%%%%%%%%%%%%%%%%%%%%%%%%%%%
%%%%%%%%%%%%%%%%%%%%%%%%%%%%%%%%%%%%%%%%%%%

%%%%%%%%%%%%%%%%%%%%%%%%%%%%%%%%%%%%%%%%%%
%%%%%%%%%%%%%%%%%%%%%%%%%%%%%%%%%%%%%%%%%%%%
%%%%%%%%%%%%%%%%%%%%%%%%%%%%%%%%%%%%%%%%%%%

\medskip

%%%%%%%%%%%%%%%%%%%%%%%%%%%%%%%%%%
%%%%%%%%%%%%%%%%%%%%%%%%%%%%%%%%%%

\section{Classification of Picard lattices of K3 surfaces}
\label{sec3:piclatt}

Classification of degenerations in Section \ref{sec2:degen} and Tables 1---4 of
Section \ref{sec4:tables} contain important classication of Picard lattices of K3 surfaces which is
the main subject of this paper.

Let $S_X$ be the Picard lattice of K\"ahlerian K3 surface $X$
and $S_X<0$ is negative definite. For algebraic K3 surfaces, instead of $S_X$,
one can take $(h)^\perp_{S_X}$ where $h\in S_X$ is a nef element with $h^2>0$.

Let $G=Aut(X)_0$ be the full symplectic automorphism subgroup
of $X$ (it is finite) and let $E_1,\dots E_k$ are
all non-singular rational curves of $X$. Then
$$
S=[(S)_G=(S_X)_G, cl(E_1),\dots cl(E_k)]_{pr}\subset S_X
$$
is the most important part of the Picard lattice $S_X$ {(\it the main part of $S_X$
or $MS_X$)}. The larger part
$S\subset S_X$ of $S_X$ is generated by the main part $S$ and by the
{\it exotic} part $ES_X=(S)^\perp_{S_X}$ of $S_X$. Classification
of possible $S=MS_X$ is the most important.

The main part $S=MS_X$
can be described purely lattice-theoretically for negative definite $S_X$; see
\cite[Remark 1.14.7]{Nik1}. Let $H(S_X)\subset O(S_X)$ be the kernel of the
action of $O(S_X)$ on the discriminant group $(S_X)^\ast/S_X$ (equivalently,
on the discriminant quadratic form $q_{S_X}$). Then the exotic part $ES_X=(S_X)^{H(S_X)}$ is
the fixed part of $H(S_X)$ on $S_X$, and the main part $MS_X=(S_X)_{H(S_X)}$ is the
coinvariant part of $H(S_X)$. The group $H(S_X)=W^{(2)}(S_X)\rtimes Aut(X)_0$
where $W^{(2)}(S_X)$ is generated by reflections in all elements $\delta\in S_X$
with $\delta^2=-2$.

K3 surfaces with the Picard lattice $S$ give the degeneration of K3 surfaces with
the symplectic automorphism group $G=Aut (X)_0$.
{\it Thus, the lattice $S=MS_X$ is one of lattices
of Tables 1 --- 4 (if $G>D_6$ and $k\ge 1$).}
The only difference is that
{\it the group $G$ must be the maximal finite symplectic automorphism group of K3 surfaces
with the Picard lattice $S$.}

\medskip

From the point of view of abstract lattices,
$e_1=cl(E_1),\dots , e_k=cl(E_k)$ give a basis of the $(-2)$ root
system of $S$ and
$$
G|S=\{g\in O(S)\ |\ g(\{e_1,...,e_k\})=\{e_1,...,e_k\})\ and\ g|(S^\ast/S)=id\}.
$$

For some two cases $(S_1,G_1)$ and $(S,G)$ of the Tables 1 --- 4, the lattices $S_1\cong S$ can be isomorphic,
but their symplectic automorphism groups $G_1\subset G$ are only subgroups of one another with different orders,
$|G_1|<|G|$, where $G$ is the maximal symplectic automorphism group of K3 surfaces
with the Picard lattice $S$, but $G_1\subset G$ is only a proper subgroup of
the maximal symplectic automorphism group $G$ of K3 surfaces with the Picard lattice $S$.

Thus, to get from Tables 1 --- 4 classification of Picards lattices $S$ of K3 surfaces,
we should find all such pairs of degenerations $(S_1,G_1)$ and $(S,G)$ of Tables 1---4.
Classification of such pairs is given in the List 1 below
as pairs
$$
(S_1,G_1)\Longleftarrow (S,G).
$$
 {\it In Tables 1---4, we mark  the case $(S_1, G_1)$ by $o$ (old).
 For the classification of Picard lattices $S$ of K3 surfaces, we can
remove this case $(S_1,G_1)$.} But, this case is important and interesting as itself:
If a K3 surface has a symplectic automorphism group $G_1$ (not necessarily maximal)
and the degeneration $S_1$, then the full symplectic automorphism
group of $X$ is larger, it is $G$, with $|G_1|<|G|$.
For $G$,  we have less orbits and less codimension of the degeneration.

The same happens for the classical case of Kummer surfaces: $|G_1|=1$, $G=(C_2)^4$, $S=[16\aaa_1]_{pr}$:
Kummer surfaces give degeneration of the type $16\aaa_1$ (of codim. $16$) of K3 surfaces with trivial
symplectic automorphism group, but their full symplectic  automorphism group increases to $(C_2)^4$. It was
shown in  \cite{Nik-1}, that any K3 surface with 16 (sixteen) $\bp^1$ and the Dynkin diagram $16\aaa_1$ is Kummer.

Thus, we have the following result.

\begin{theorem} Classification of the main parts $S=MS_X$ of
Picard lattices $S_X<0$ of K3 surfaces with full finite symplectic
automorphism group $G$ of order $|G|\ge 8$ (equivalently, ${\bf n}\ge 9$,
or $G$ is different from $D_6$, $C_4$, $(C_2)^2$, $C_3$, $C_2$ and $C_1$)
and with at least one
non-singular rational curve is given in Tables 1---4 of
Section \ref{sec4:tables} in lines which are not marked by $o$.

Lines $(G_1,S_1)$ of Tables 1---4 of Section \ref{sec4:tables} which are marked
by $o$ give classification of degenerations $X$ of K\"ahlerian K3 surfaces
with finite symplectic automorphism groups $G_1$ such that the full
finite symplectic automorphism group $G$ of $X$ is larger than $G_1$, it contains $G_1$ as
a proper subroup. Such cases are shown in the List 1 below as
$$
(G_1,S_1)\Longleftarrow (G,S)
$$
where $(G,S)$ gives the corresponding line of Tables 1---4
with maximal finite symplectic automorphism group $G$ and $S\cong S_1$.
\label{theorem1}
\end{theorem}

\vskip2cm

{\bf List 1:}

{\it The list or cases, when a degeneration of K3 surfaces with symplectic automorphism
group $G_1$ has, actually,  the maximal finite symplectic automorphism group $G$ which
contains $G_1$ and $|G_1|<|G|$. The group $G$ has less orbits and less codimension
of the degeneration than $G_1$.}

\medskip

$({\bf n}=9,\ ((2\aaa_1,2\aaa_1)\subset 4\aaa_1)_{II})\
\Longleftarrow ({\bf n}=21,4\aaa_1)$;

\medskip

$({\bf n}=9,\ (8\aaa_1,8\aaa_1)\subset 16\aaa_1)\
\Longleftarrow ({\bf n}=21,16\aaa_1)$;

\medskip

$({\bf n}=9,\ ((2\aaa_1,2\aaa_1)_{II},4\aaa_1)\subset 8\aaa_1)\
\Longleftarrow\  ({\bf n}=21,(4\aaa_1,4\aaa_1)\subset 8\aaa_1)$;

\medskip

({\bf n}=9,\
$
\left(\begin{array}{cccc}
 2\aaa_1 & (4\aaa_1)_I & (4\aaa_1)_I  & (4\aaa_1)_I \\
         & 2\aaa_1 & (4\aaa_1)_I & (4\aaa_1)_I  \\
         &         & 2\aaa_1 &    (4\aaa_1)_I \\
         &         &         & 2\aaa_1
\end{array}\right)
\subset 8\aaa_1) \Longleftarrow ({\bf n}=40,\ 8\aaa_1)$;

\medskip

$({\bf n}=9,\ ((2\aaa_1,2\aaa_1)_{II},4\aaa_1,4\aaa_1)
\subset 12\aaa_1)\Longleftarrow ({\bf n}=49,\ 12\aaa_1)$;

\medskip

({\bf n}=9,\
$
\left(\begin{array}{cccc}
 2\aaa_1 & 2\aaa_3     & (4\aaa_1)_I & 6\aaa_1 \\
         & 4\aaa_1     & 6\aaa_1     & 8\aaa_1 \\
         &             & 2\aaa_1     & 2\aaa_3 \\
         &             &              & 4\aaa_1
\end{array}\right)
\subset 4\aaa_3)\
\Longleftarrow\ ({\bf n}=22,(4\aaa_1,8\aaa_1)\subset 4\aaa_3)$;

\medskip

$({\bf n}=9,\ (4\aaa_1,4\aaa_1,4\aaa_1,4\aaa_1)
\subset 16\aaa_1)\Longleftarrow ({\bf n}=39,\ 16\aaa_1)$;

\medskip

({\bf n}=9,\
$
\left(\begin{array}{ccccc}
 2\aaa_1 & (4\aaa_1)_I & (4\aaa_1)_I  & (4\aaa_1)_I &  10\aaa_1 \\
         & 2\aaa_1     & (4\aaa_1)_I  & (4\aaa_1)_I &  10\aaa_1 \\
         &             & 2\aaa_1      & (4\aaa_1)_I &  10\aaa_1 \\
         &             &              & 2\aaa_1     &  10\aaa_1\\
         &             &              &             &   8\aaa_1
\end{array}\right)
\subset 16\aaa_1 ) \Longleftarrow ({\bf n}=56,\ 16\aaa_1)$;

\medskip

({\bf n}=9,\
$
\left(\begin{array}{ccccc}
 2\aaa_1 & (4\aaa_1)_I & 6\aaa_1    & (4\aaa_1)_I     &  6\aaa_1 \\
         & 2\aaa_1     & 2\aaa_3  & (4\aaa_1)_I       &  6\aaa_1 \\
         &             & 4\aaa_1      & 6\aaa_1   &  8\aaa_1 \\
         &             &              & 2\aaa_1     & 2\aaa_3\\
         &             &              &             & 4\aaa_1
\end{array}\right)
\subset 2\aaa_1\amalg 4\aaa_3)
$

$\Longleftarrow\ ({\bf n}=22,\
\left(\begin{array}{rrr}
 2\aaa_1 & 6\aaa_1 & 10\aaa_1  \\
         & 4\aaa_1 & 4\aaa_3   \\
         &         & 8\aaa_1
\end{array}\right)\subset 2\aaa_1\amalg 4\aaa_3)$;

\medskip

$({\bf n}=9,\ ((2\aaa_1,2\aaa_1)_{II},4\aaa_1,4\aaa_1,4\aaa_1)
\subset 16\aaa_1)\Longleftarrow ({\bf n}=75,\ 16\aaa_1)$;

\medskip

$({\bf n}=9,
\left(\begin{array}{ccccc}
 2\aaa_1 & 2\aaa_3     & (4\aaa_1)_I  &  6\aaa_1     &  6\aaa_1 \\
         & 4\aaa_1     & 6\aaa_1     &  8\aaa_1     &  8\aaa_1 \\
         &             & 2\aaa_1      & 2\aaa_3     &  6\aaa_1 \\
         &             &              & 4\aaa_1     &  8\aaa_1\\
         &             &              &             &  4\aaa_1
\end{array}\right)
\subset 4\aaa_3\amalg  4\aaa_1)
$

$\Longleftarrow
({\bf n}=22,\
\left(\begin{array}{rrr}
 4\aaa_1 & (8\aaa_1)_{II} & 12\aaa_1  \\
         & 4\aaa_1 & 4\aaa_3   \\
         &         & 8\aaa_1
\end{array}\right)\subset 4\aaa_1\amalg 4\aaa_3)$;

\medskip

({\bf n}=10, $(\aaa_1,\aaa_1)\subset 2\aaa_1)\ \Longleftarrow\ ({\bf n}=22,2\aaa_1)$;

\medskip

({\bf n}=10,
$((4\aaa_1,4\aaa_1)\subset 8\aaa_1)_{II})
\Longleftarrow\ ({\bf n}=22,8\aaa_1)$;

\medskip

({\bf n}=10,
$(\aaa_1,\aaa_1,(2\aaa_1)_I)\subset 4\aaa_1)\ \Longleftarrow ({\bf n}=39,\ 4\aaa_1)$;

\medskip

({\bf n}=10,
$(\aaa_1,\aaa_1,4\aaa_1)\subset 6\aaa_1)\
\Longleftarrow\ ({\bf n}=22,(2\aaa_1,4\aaa_1)\subset 6\aaa_1)$;

\medskip

({\bf n}=10,
$(\aaa_1,\aaa_1,8\aaa_1)\subset 10\aaa_1)\
\Longleftarrow\ ({\bf n}=22,(2\aaa_1,8\aaa_1)\subset 10\aaa_1)$;

\medskip

({\bf n}=10,
$
\left(\begin{array}{ccc}
(2\aaa_1)_{I} & (6\aaa_1)_I & (6\aaa_1)_I  \\
         & 4\aaa_1   & (8\aaa_1)_{II}   \\
         &         & 4\aaa_1
\end{array}\right)
\subset 10\aaa_1)\
\Longleftarrow\ ({\bf n}=22,(2\aaa_1,8\aaa_1)\subset 10\aaa_1)$;

\medskip

({\bf n}=10,
$
\left(\begin{array}{ccc}
(2\aaa_1)_{II}  & 6\aaa_1 & 6\aaa_1   \\
         & 4\aaa_1 & (8\aaa_1)_{II}   \\
         &         & 4\aaa_1
\end{array}\right)
\subset 10\aaa_1)\
\Longleftarrow\ ({\bf n}=22,(2\aaa_1,8\aaa_1)\subset 10\aaa_1)$;

\medskip

({\bf n}=10,
$
\left(\begin{array}{ccc}
4\aaa_1  & (8\aaa_1)_I & (8\aaa_1)_{II}  \\
         & 4\aaa_1   & (8\aaa_1)_{I} \\
         &         & 4\aaa_1
\end{array}\right)
\subset 12\aaa_1)\
\Longleftarrow\ ({\bf n}=22,(4\aaa_1,8\aaa_1)\subset 12\aaa_1)$;

\medskip

({\bf n}=10,
$
\left(\begin{array}{ccc}
4\aaa_1  & 4\aaa_2 & (8\aaa_1)_{II}   \\
         & 4\aaa_1 & 4\aaa_2   \\
         &         & 4\aaa_1
\end{array}\right)
\subset 4\aaa_3)
\Longleftarrow\ ({\bf n}=22,(4\aaa_1,8\aaa_1)\subset 4\aaa_3)$;

\medskip

({\bf n}=10,\
$
\left(\begin{array}{ccc}
4\aaa_1 & (8\aaa_1)_{II} & 12\aaa_1   \\
         & 4\aaa_1 & 12\aaa_1   \\
         &         & 8\aaa_1
\end{array}\right)
\subset 16\aaa_1) \Longleftarrow ({\bf n}=39,\ 16\aaa_1)$;

\medskip

$({\bf n}=10,\
(4\aaa_1,4\aaa_1,2\aaa_2)\subset 6\aaa_2)
\Longleftarrow ({\bf n}=34,\ 6\aaa_2)$;

\medskip

({\bf n}=10,\
$
\left(\begin{array}{cccc}
\aaa_1 & 2\aaa_1 & 3\aaa_1      &  5\aaa_1       \\
       & \aaa_1  & 3\aaa_1      &  5\aaa_1       \\
       &         & (2\aaa_1)_I  & (6\aaa_1)_I   \\
       &         &              &  4\aaa_1
\end{array}\right)
\subset 8\aaa_1) \Longleftarrow ({\bf n}=56,\ 8\aaa_1)$;

\medskip

$({\bf n}=10,\ (\aaa_1,\aaa_1,(2\aaa_1)_{I},8\aaa_1)
\subset 12\aaa_1)\Longleftarrow ({\bf n}=65,\ 12\aaa_1)$;

\medskip

({\bf n}=10,\
$
\left(\begin{array}{cccc}
\aaa_1 & 2\aaa_1 & 5\aaa_1      &  5\aaa_1       \\
       & \aaa_1  & 5\aaa_1      &  5\aaa_1       \\
       &         & 4\aaa_1      & (8\aaa_1)_I    \\
       &         &              &  4\aaa_1
\end{array}\right)
\subset 10\aaa_1)\
\Longleftarrow ({\bf n}=22,\ (2\aaa_1, (4\aaa_1,4\aaa_1)_{II})\subset 10\aaa_1)$;

\medskip

$({\bf n}=10,\  (\aaa_1,\,\aaa_1,\,4\aaa_1,\,8\aaa_1)\subset 14\aaa_1)\
\Longleftarrow\ ({\bf n}=22,\ (2\aaa_1,4\aaa_1,8\aaa_1)\subset 14\aaa_1)$;

\medskip

$({\bf n}=10,\ (\aaa_1,\,\aaa_1,\,4\aaa_1,\, 8\aaa_1)\subset 2\aaa_1\amalg 4\aaa_3)$

$\Longleftarrow\ ({\bf n}=22,\
\left(\begin{array}{rrr}
 2\aaa_1 & 6\aaa_1 & 10\aaa_1  \\
         & 4\aaa_1 & 4\aaa_3   \\
         &         & 8\aaa_1
\end{array}\right)\subset 2\aaa_1\amalg 4\aaa_3)$;

\medskip

({\bf n}=10,\
$
\left(\begin{array}{cccc}
\aaa_1 & \aaa_3              & 3\aaa_1          &  5\aaa_1     \\
       & (2\aaa_1)_{II}      &  4\aaa_1         &  6\aaa_1     \\
       &                     & (2\aaa_1)_I      &  2\aaa_3     \\
       &                     &                  &  4\aaa_1
\end{array}\right)
\subset 3\aaa_3)\
\Longleftarrow ({\bf n}=34,\ (3\aaa_1,(6\aaa_1)_{II})\subset 3\aaa_3)$;

\medskip

({\bf n}=10,\
$
\left(\begin{array}{cccc}
\aaa_1 & 3\aaa_1            & 5\aaa_1          &  9\aaa_1     \\
       & (2\aaa_1)_{I}      &  (6\aaa_1)_I     &  10\aaa_1     \\
       &                     & 4\aaa_1         &  12\aaa_1     \\
       &                     &                  &  8\aaa_1
\end{array}\right)
\subset 15\aaa_1)\
\Longleftarrow ({\bf n}=34,\ (3\aaa_1,12\aaa_1))\subset 15\aaa_1))$;

\medskip

({\bf n}=10,\
$(\aaa_1,\,4\aaa_1,\,4\aaa_1,\,2\aaa_2)\subset \aaa_1\amalg 6\aaa_2)
\Longleftarrow ({\bf n}=34,\ (\aaa_1,6\aaa_2)\subset \aaa_1\amalg 6\aaa_2))$;

\medskip

({\bf n}=10,\
$
\left(\begin{array}{cccc}
(2\aaa_1)_I  & 4\aaa_1        & (6\aaa_1)_I  &  (6\aaa_1)_{I}     \\
             & (2\aaa_1)_{II} &  6\aaa_1     &  6\aaa_1     \\
       &                      & 4\aaa_1      &  (8\aaa_1)_{II}     \\
       &                      &              &  4\aaa_1
\end{array}\right)
\subset 12\aaa_1)\  \Longleftarrow ({\bf n}=65,\ 12\aaa_1)$;

\medskip

({\bf n}=10,\
$
\left(\begin{array}{cccc}
(2\aaa_1)_I & (6\aaa_1)_{II}   & (6\aaa_1)_I     &  (6\aaa_1)_I     \\
            & 4\aaa_1          & (8\aaa_1)_{I}   &  (8\aaa_1)_I     \\
       &                       & 4\aaa_1         &   (8\aaa_1)_{II} \\
       &                       &                 &  4\aaa_1
\end{array}\right)
\subset 14\aaa_1)$

$\Longleftarrow\ ({\bf n}=22,\ (2\aaa_1,4\aaa_1,8\aaa_1)\subset 14\aaa_1)$;

\medskip

({\bf n}=10,\
$
\left(\begin{array}{cccc}
(2\aaa_1)_{II} & 6\aaa_1     & 6\aaa_1          &  6\aaa_1        \\
               & 4\aaa_1     &(8\aaa_1)_I       &  (8\aaa_1)_I    \\
               &             & 4\aaa_1          &  (8\aaa_1)_{II} \\
       &                     &                  &  4\aaa_1
\end{array}\right)
\subset 14\aaa_1)
$

$\Longleftarrow\ ({\bf n}=22,\ (2\aaa_1,4\aaa_1,8\aaa_1)\subset 14\aaa_1)$;

\medskip

({\bf n}=10,\
$
\left(\begin{array}{cccc}
(2\aaa_1)_{I} & (6\aaa_1)_I  & (6\aaa_1)_I   & (6\aaa_1)_I      \\
               & 4\aaa_1     & 4\aaa_2       &  (8\aaa_1)_{II}  \\
               &             & 4\aaa_1       &  4\aaa_2         \\
       &                     &                  &  4\aaa_1
\end{array}\right)
\subset 2\aaa_1\amalg 4\aaa_3)
$,\

({\bf n}=10,\
$
\left(\begin{array}{cccc}
(2\aaa_1)_{II} & 6\aaa_1     & 6\aaa_1       & 6\aaa_1      \\
               & 4\aaa_1     & 4\aaa_2       &  (8\aaa_1)_{II}  \\
               &             & 4\aaa_1       &  4\aaa_2         \\
       &                     &                  &  4\aaa_1
\end{array}\right)
\subset 2\aaa_1\amalg 4\aaa_3)
$

$\Longleftarrow\ ({\bf n}=22,\
\left(\begin{array}{rrr}
 2\aaa_1 & 6\aaa_1 & 10\aaa_1  \\
         & 4\aaa_1 & 4\aaa_3   \\
         &         & 8\aaa_1
\end{array}\right)\subset 2\aaa_1\amalg 4\aaa_3)$;

\medskip

({\bf n}=10,\
$
\left(\begin{array}{cccc}
(2\aaa_1)_I & 2\aaa_3          & (6\aaa_1)_I     &  10\aaa_1  \\
            & 4\aaa_1          & (8\aaa_1)_{I}   &  12\aaa_1  \\
       &                       & 4\aaa_1         &  4\aaa_3   \\
       &                       &                 &  8\aaa_1
\end{array}\right)
\subset 6\aaa_3) \
\Longleftarrow ({\bf n}=34,\ ((6\aaa_1)_I,12\aaa_1)\subset 6\aaa_3)$;

\medskip

({\bf n}=10,\
$
\left(\begin{array}{cccc}
4\aaa_1 & (8\aaa_1)_{II}     & (8\aaa_1)_I        &  (8\aaa_1)_I  \\
              & 4\aaa_1            & (8\aaa_1)_I  &  (8\aaa_1)_I  \\
               &                   & 4\aaa_1      &  (8\aaa_1)_{II}\\
       &                     &                    &  4\aaa_1
\end{array}\right)
\subset 16\aaa_1) \Longleftarrow ({\bf n}=56,\ 16\aaa_1)$;

\medskip

({\bf n}=10,\
$
\left(\begin{array}{cccc}
4\aaa_1 & (8\aaa_1)_I        & (8\aaa_1)_I      &  (8\aaa_1)_I   \\
               & 4\aaa_1     & 4\aaa_2          &  (8\aaa_1)_{II}\\
               &             & 4\aaa_1          &  4\aaa_2       \\
       &                     &                  &  4\aaa_1
\end{array}\right)
\subset 4\aaa_1\amalg 4\aaa_3)
$

$\Longleftarrow
({\bf n}=22,\
\left(\begin{array}{rrr}
 4\aaa_1 & (8\aaa_1)_{II} & 12\aaa_1  \\
         & 4\aaa_1 & 4\aaa_3   \\
         &         & 8\aaa_1
\end{array}\right)\subset 4\aaa_1\amalg 4\aaa_3)$;

\medskip

({\bf n}=10,\
$
\left(\begin{array}{cccc}
4\aaa_1 & (8\aaa_1)_{I}     & (8\aaa_1)_I    & 4\aaa_1\amalg 2\aaa_2   \\
               & 4\aaa_1     & 4\aaa_2       &  4\aaa_1\amalg 2\aaa_2  \\
               &             & 4\aaa_1       &  4\aaa_1\amalg 2\aaa_2  \\
       &                     &               &  2\aaa_2
\end{array}\right)
\subset 4\aaa_1\amalg 6\aaa_2)$\

$\Longleftarrow ({\bf n}=34,\ (4\aaa_1,6\aaa_2)\subset 4\aaa_1\amalg 6\aaa_2)$;

\medskip

$({\bf n}=12,\ (8\aaa_1,8\aaa_1)\subset 16\aaa_1)\Longleftarrow ({\bf n}=75,\ 16\aaa_1)$;

\medskip

$({\bf n}=12,\ (\aaa_2,\aaa_2)\subset 2\aaa_2)\Longleftarrow ({\bf n}=26,\ 2\aaa_2)$;

\medskip

$({\bf n}=16,\
(5\aaa_1,5\aaa_1,5\aaa_1)\subset 15\aaa_1)
\Longleftarrow ({\bf n}=55,\ 15\aaa_1)$;

\medskip

$({\bf n}=17,\ (\aaa_1,\aaa_1)
\subset 2\aaa_1)\Longleftarrow ({\bf n}=34,\ 2\aaa_1)$;

\medskip

$({\bf n}=17,\ (\aaa_1,3\aaa_1)
\subset 4\aaa_1)\Longleftarrow ({\bf n}=49,\ 4\aaa_1)$;

\medskip

$({\bf n}=17,\ (4\aaa_1,4\aaa_1)
\subset 8\aaa_1)\Longleftarrow ({\bf n}=34,\ 8\aaa_1)$;

\medskip

$({\bf n}=17,\ (4\aaa_1,12\aaa_1)
\subset 16\aaa_1)\Longleftarrow ({\bf n}=49,\ 16\aaa_1)$;

\medskip

$({\bf n}=17,\ (6\aaa_1,6\aaa_1)
\subset 12\aaa_1)\Longleftarrow ({\bf n}=49,\ 12\aaa_1)$;

\medskip

$({\bf n}=17,\ (6\aaa_1,6\aaa_1)
\subset 6\aaa_2)\Longleftarrow ({\bf n}=34,\ 6\aaa_2)$;

\medskip

$({\bf n}=17,\ (\aaa_1,\aaa_1,\aaa_1)
\subset 3\aaa_1)\Longleftarrow ({\bf n}=61,\ 3\aaa_1)$;

\medskip

$({\bf n}=17,\ (\aaa_1,\aaa_1,4\aaa_1)
\subset 6\aaa_1)\Longleftarrow ({\bf n}=34,\ (2\aaa_1,4\aaa_1)\subset 6\aaa_1)$;

\medskip

$({\bf n}=17,\ (\aaa_1,\aaa_1,6\aaa_1)
\subset 8\aaa_1)\Longleftarrow ({\bf n}=34,\ (2\aaa_1,(6\aaa_1)_{II})\subset 8\aaa_1)$;

\medskip

$({\bf n}=17,\ (\aaa_1,\aaa_1,12\aaa_1)\subset 14\aaa_1)
\Longleftarrow ({\bf n}=34,\ (2\aaa_1,12\aaa_1)\subset 14\aaa_1)$;

\medskip

$({\bf n}=17,\ (\aaa_1,3\aaa_1,4\aaa_1)
\subset 8\aaa_1)\Longleftarrow ({\bf n}=65,\ 8\aaa_1)$;

\medskip

$({\bf n}=17,\ (\aaa_1,3\aaa_1,12\aaa_1)
\subset 16\aaa_1)\Longleftarrow ({\bf n}=75,\ 16\aaa_1)$;

\medskip

$({\bf n}=17,\ (\aaa_1,4\aaa_1,4\aaa_1)
\subset 9\aaa_1)\Longleftarrow ({\bf n}=34,\ (\aaa_1,8\aaa_1)\subset 9\aaa_1)$;

\medskip

$({\bf n}=17,\ (\aaa_1,6\aaa_1,6\aaa_1)
\subset \aaa_1\amalg 6\aaa_2)
\Longleftarrow ({\bf n}=34,\ (\aaa_1,6\aaa_2)\subset \aaa_1\amalg 6\aaa_2)$;

\medskip

$({\bf n}=17,\  (3\aaa_1,4\aaa_1,4\aaa_1)\subset 11\aaa_1)
\Longleftarrow ({\bf n}=34,\ (3\aaa_1,8\aaa_1)\subset 11\aaa_1)$;

\medskip

$({\bf n}=17,\ (4\aaa_1,4\aaa_1,4\aaa_1)\subset 12\aaa_1)
\Longleftarrow ({\bf n}=61,\ 12\aaa_1)$;

\medskip

({\bf n}=17,\
$
\left(\begin{array}{rrr}
4\aaa_1 & 4\aaa_2 & 8\aaa_1  \\
         & 4\aaa_1 & 4\aaa_2   \\
         &         & 4\aaa_1
\end{array}\right)
\subset 4\aaa_3)\
\Longleftarrow ({\bf n}=34,\ (4\aaa_1,8\aaa_1)\subset 4\aaa_3)$;

\medskip

$({\bf n}=17,\ (4\aaa_1,4\aaa_1,6\aaa_1)\subset 14\aaa_1)\
\Longleftarrow ({\bf n}=34,\ ((6\aaa_1)_I,8\aaa_1)\subset 14\aaa_1)$;

\medskip

$({\bf n}=17,\ (4\aaa_1,6\aaa_1,6\aaa_1)
\subset 16\aaa_1)\Longleftarrow ({\bf n}=75,\ 16\aaa_1)$;

\medskip

$({\bf n}=17,\ (4\aaa_1,6\aaa_1,6\aaa_1)\subset 4\aaa_1\amalg 6\aaa_2)
\Longleftarrow ({\bf n}=34,\ (4\aaa_1,6\aaa_2)\subset 4\aaa_1\amalg 6\aaa_2)$;

\medskip

$({\bf n}=21,\ (4\aaa_1, 4\aaa_1,4\aaa_1)
\subset 12\aaa_1)\Longleftarrow ({\bf n}=49,\ 12\aaa_1)$;

\medskip

$({\bf n}=21,\ (4\aaa_1, 4\aaa_1,4\aaa_1,4\aaa_1)
\subset 16\aaa_1)\Longleftarrow ({\bf n}=75,\ 16\aaa_1)$;

\medskip

$({\bf n}=22,\ (2\aaa_1,2\aaa_1)
\subset 4\aaa_1)\Longleftarrow ({\bf n}=39,\ 4\aaa_1)$;

\medskip

$({\bf n}=22,\ ((4\aaa_1,4\aaa_1)
\subset 8\aaa_1)_I)\Longleftarrow ({\bf n}=40,\ 8\aaa_1)$;

\medskip

$({\bf n}=22,\ (8\aaa_1,8\aaa_1)
\subset 16\aaa_1)\ \Longleftarrow ({\bf n}=39,\ 16\aaa_1)$;

\medskip

$({\bf n}=22,\ (2\aaa_1, 2\aaa_1,4\aaa_1)
\subset 8\aaa_1)\Longleftarrow ({\bf n}=56,\ 8\aaa_1)$;

\medskip

$({\bf n}=22,\ (2\aaa_1, 2\aaa_1,8\aaa_1)
\subset 12\aaa_1)\Longleftarrow ({\bf n}=65,\ 12\aaa_1)$;

\medskip

$({\bf n}=22,\ ((4\aaa_1,4\aaa_1)_I,8\aaa_1)
\subset 16\aaa_1)\Longleftarrow ({\bf n}=56,\ 16\aaa_1)$;

\medskip

$({\bf n}=34,\ (\aaa_1,\aaa_1)\subset 2\aaa_1)\Longleftarrow ({\bf n}=51,\ 2\aaa_1)$;

\medskip

$({\bf n}=34,\ (\aaa_1,2\aaa_1)\subset 3\aaa_1)\Longleftarrow ({\bf n}=61,\ 3\aaa_1)$;

\medskip

$({\bf n}=34,\ (\aaa_1,3\aaa_1)\subset 4\aaa_1)\Longleftarrow ({\bf n}=65,\ 4\aaa_1)$;

\medskip

$({\bf n}=34,\ (4\aaa_1,4\aaa_1)\subset 8\aaa_1)\Longleftarrow ({\bf n}=51,\ 8\aaa_1)$;

\medskip

$({\bf n}=34,\ (4\aaa_1,8\aaa_1)\subset 12\aaa_1)\Longleftarrow ({\bf n}=61,\ 12\aaa_1)$;

\medskip

$({\bf n}=34,\ (4\aaa_1,12\aaa_1)\subset 16\aaa_1)\Longleftarrow ({\bf n}=65,\ 16\aaa_1)$;

\medskip

$({\bf n}=34,\ ((6\aaa_1)_{I},(6\aaa_1)_{II})\subset 12\aaa_1)\Longleftarrow ({\bf n}=65,\ 12\aaa_1)$;

\medskip

$({\bf n}=39,\ (4\aaa_1,4\aaa_1)\subset 8\aaa_1)\Longleftarrow ({\bf n}=56,\ 8\aaa_1)$;

\medskip

$({\bf n}=39,\ (4\aaa_1,8\aaa_1)\subset 12\aaa_1)\Longleftarrow ({\bf n}=65,\ 12\aaa_1)$;

\medskip

$({\bf n}=39,\ (8\aaa_1,8\aaa_1)\subset 16\aaa_1)\Longleftarrow ({\bf n}=75,\ 16\aaa_1)$;

\medskip

$({\bf n}=40,\ (8\aaa_1,8\aaa_1)\subset 16\aaa_1)\Longleftarrow ({\bf n}=56,\ 16\aaa_1)$;

\medskip

$({\bf n}=49,\ (4\aaa_1,4\aaa_1)\subset 8\aaa_1)\Longleftarrow ({\bf n}=65,\ 8\aaa_1)$;

\medskip

$({\bf n}=49,\ (4\aaa_1,12\aaa_1)\subset 16\aaa_1)\Longleftarrow ({\bf n}=75,\ 16\aaa_1)$.

\medskip

%%%%%%%%%%%%%%%%%%%%%%%%%%%%%%%%%%%%%
%%%%%%%%%%%%%%%%%%%%%%%%%%%%%%%%%%%%%%

\begin{proof}
Assume that for two degenerations given by two lines
$({\bf n}, Deg, S)$ and $({\bf n}_1, Deg_1, S_1)$ of
Tables 1---4 we have
\begin{equation}
({\bf n_1}, Deg_1, S_1)\Longleftarrow ({\bf n}, Deg, S).
\label{pairslist1}
\end{equation}
Then the group $G$ (marked by ${\bf n}$) is maximal
symplectic automorphism group for $S$, or shortly
\begin{equation}
the\ degeneration\ ({\bf n},Deg,S)\ is\ maximal;
\label{A}
\end{equation}
group $G_1$ (marked by ${\bf n}_1$)
is isomorphic to a proper subgroup of $G$ (marked by ${\bf n}$), or shortly,
\begin{equation}
{\bf n}_1\subset {\bf n},
\label{B}
\end{equation}
full Dynkin diagrams of $Deg$ and $Deg_1$ are isomorphic:
\begin{equation}
Dyn(Deg_1)\cong Dyn(Deg);
\label{C}
\end{equation}
genuses of lattices $S_1$ and $S$ are equal, equivalently,
\begin{equation}
\rk S_1=\rk S\ and\ q_{S_1}\cong q_S;
\label{D}
\end{equation}
finally, each marking by Niemeier lattices of the degeneration
$({\bf n}_1, Deg_1, S_1)$ is isomorphic to a restriction on
an appropriate subgroup $G_1\subset G$
of a marking of the degeneration $({\bf n}, Deg, S)$
by Niemeier lattices, or shortly,
\begin{equation}
Niemeier\ markings\ of\ ({\bf n}_1, Deg_1, S_1)\ and\ ({\bf n}, Deg, S)\ agree.
\label{E}
\end{equation}

Using calculations by Hashimoto in \cite{Hash} (which can be easily done by
the Program GAP, see \cite{GAP}), we can find all possible pairs
${\bf n}$, ${\bf n_1}$ from Table 1 such that ${\bf n}\supset {\bf n}_1$.
They are given in the list of subgroups \eqref{listsubgroups} below.
\medskip

\begin{equation}
{\bf List\ of\ subgroups}
\label{listsubgroups}
\end{equation}

$$
{\bf n}=75\supset {\bf n}_1=1,2,3,4,9,10,12,17,21,22,39,49;
$$
$$
{\bf n}=65\supset {\bf n}_1=1,2,3,4,6,9,10,17,21,22,34,39,49;
$$
$$
{\bf n}=61\supset {\bf n}_1=1,2,3,4,6,10,17,18,30,34;
$$
$$
{\bf n}=56 \supset {\bf n}_1=1,3,4,9,10,21,22,39,40;
$$
$$
{\bf n}=55 \supset {\bf n}_1=1,2,3,6,16,17;
$$
$$
{\bf n}=51 \supset {\bf n}_1=1,2,3,4,6,9,10,17,18,22,34;
$$
$$
{\bf n}=49 \supset {\bf n}_1=1,2,3,9,17,21;
$$
$$
{\bf n}=48 \supset {\bf n}_1=1,2,3,6,18,30;
$$
$$
{\bf n}=46 \supset {\bf n}_1=1,2,4,6,30;
$$
$$
{\bf n}=40 \supset {\bf n}_1=1,3,4,9,10,22;
$$
$$
{\bf n}=39 \supset {\bf n}_1=1,3,4,9,10,21,22;
$$
$$
{\bf n}=34 \supset {\bf n}_1=1,2,3,4,6,10,17;
$$
$$
{\bf n}=33 \supset {\bf n}_1=2;
$$
$$
{\bf n}=32 \supset {\bf n}_1=1,4,16;
$$
$$
{\bf n}=30 \supset {\bf n}_1=1,2,6;
$$
$$
{\bf n}=26 \supset {\bf n}_1=1,3,4,10,12;
$$
$$
{\bf n}=22 \supset {\bf n}_1=1,3,4,9,10;
$$
$$
{\bf n}=21 \supset {\bf n}_1=1,3,9;
$$
$$
{\bf n}=18 \supset {\bf n}_1=1,2,3,6;
$$
$$
{\bf n}=17 \supset {\bf n}_1=1,2,3;
$$
$$
{\bf n}=16 \supset {\bf n}_1=1;
$$
$$
{\bf n}=12 \supset {\bf n}_1=1,4;
$$
$$
{\bf n}=10 \supset {\bf n}_1=1,3,4;
$$
$$
{\bf n}=9 \supset {\bf n}_1=1,3;
$$
$$
{\bf n}=6 \supset {\bf n}_1=1,2;
$$
$$
{\bf n}=4 \supset {\bf n}_1=1;
$$
$$
{\bf n}=3 \supset {\bf n}_1=1.
$$

All degenerations $({\bf n}, Deg, S)$ with one orbit of $(-2)$-curves
(equivalently, of the codimension $1$) are, obviously, maximal.
They are classified in Table 1. Using the list of subgroups \eqref{listsubgroups},
Tables 1---4 and our description of markings by Niemeier lattices in \cite{Nik10} and \cite{Nik11},
firstly, we find all pairs
$({\bf n_1}, Deg_1, S_1)\Longleftarrow ({\bf n}, Deg, S)$ where the degeneration
$({\bf n}, Deg, S)$ has one orbit. They give pairs of List 1
where $({\bf n}, Deg, S)$ belongs to Table 1. The corresponding cases
$({\bf n_1}, Deg_1, S_1)$ we mark by $o$ in Tables 1---4.

Then all remaining non-marked by $o$
degenerations (lines of Tables 2---4) with two orbits (equivalently, of the codimension $2$)
are maximal. For them, ${\bf n}\le 34$.
They are given in Tables 2---4 as degenerations with two orbits which are not marked by $o$.
Let us denote the set of them by $D_2$.

Next, in 3rd step, using the list of subgroups \eqref{listsubgroups}, Tables 2---4 and our description of
markings by Niemeier lattices in \cite{Nik10} and \cite{Nik11}, we find all
pairs $({\bf n_1}, Deg_1, S_1)\Longleftarrow ({\bf n}, Deg, S)$ where the degeneration
$({\bf n}, Deg, S)$ belongs to the set $D_2$. They are given in List 1 as
pairs $({\bf n_1}, Deg_1, S_1)\Longleftarrow ({\bf n}, Deg, S)$ where
the degeneration $({\bf n}, Deg, S)$ has two orbits. The corresponding
degenerations $({\bf n_1}, Deg_1, S_1)$ we mark by $o$.

Let us denote the set of remaining (not marked by $o$)
degenerations with $3$ orbits of Tables 2---4 as $D_3$.
For them, ${\bf n}\le 30$.

Next, in 4th step, using the list of subgroups \eqref{listsubgroups},
Tables 2---4 and our description of markings by Niemeier lattices in \cite{Nik10}
and \cite{Nik11}, we find all pairs $({\bf n_1}, Deg_1, S_1)\Longleftarrow ({\bf n}, Deg, S)$
where the degeneration $({\bf n}, Deg, S)$ belongs to $D_3$. They are given in List 1 as
pairs $({\bf n_1}, Deg_1, S_1)\Longleftarrow ({\bf n}, Deg, S)$ where the
degeneration $({\bf n}, Deg, S)$ has three orbits (equivalently, it has the codimension $3$).
The corresponding degenerations $({\bf n_1}, Deg_1, S_1)$ of Tables 2---4 we mark by $o$.

Let us denote the set of remaining (nor marked by $o$) degenerations with $4$ orbits (equivalently,
of the codimension $4$) of Tables 2---4 as $D_4$. For them, ${\bf n}\le 10$.

Next, in 5th step, using the list of subgroups \eqref{listsubgroups},
Tables 2---4, we find all pairs \linebreak $({\bf n_1}, Deg_1, S_1)\Longleftarrow ({\bf n}, Deg, S)$
where the degeneration $({\bf n}, Deg, S)$ belongs to $D_4$ and ${\bf n}_1\ge 8$.
We see that there are no such pairs satisfying conditions \eqref{A}, \eqref{B}.
This gives the List 1 of the Theorem.

Below, for all cases of the List 1, we give details of the check of the most complicated
condition \eqref{E}. They are also important as itself because they describe exact isomorphisms
between lattices $S$ and $S_1$ of the List 1.

%%%%%%%%%%%%%%%%%%%%%%%%%%%%%%%%%%%%%%%%%%
%%%%%%%%%%%%%%%%%%%%%%%%%%%%%%%%%%%%%%%%%%%
%%%%%%%%%%%%%%%%%%%%%%%%%%%%%%%%%%%%%%%%%%%

\medskip

Case 1:
$({\bf n}=9,\ ((2\aaa_1,2\aaa_1)\subset 4\aaa_1)_{II})\
\Longleftarrow ({\bf n}=21,4\aaa_1)$.

For ${\bf n}=21$, the group $G\cong (C_2)^4$. By \cite{Nik10},
the case $({\bf n}=21,4\aaa_1)$ is marked by the Niemeier
lattice $N_{23}$ with $G=H_{21,1}\subset A_{23}$
and by any of $4$-elements orbits
$\{\alpha_1,\alpha_{16},\alpha_{14},\alpha_{18}\}$,
$\{\alpha_2,\alpha_{20},\alpha_{19},\alpha_{24}\}$,
$\{\alpha_3,\alpha_{10},\alpha_{5},\alpha_6\}$,
$\{\alpha_8,\alpha_{11},\alpha_{9},\alpha_{21}\}$,
$\{\alpha_{12},\alpha_{22},\alpha_{23},\alpha_{17}\}$
of $H_{21,1}$ in the set of simple $(-2)$-roots $\alpha_1, \dots \alpha_{24}$
of $N_{23}$. All other orbits of $H_{21,1}$ have one element. Exactly (see \cite{Nik10}),
the group $G=H_{21,1}$ is equal to
$$
G=H_{21,1}=
$$
$$
[(\alpha_2\alpha_{20})(\alpha_{3}\alpha_{10})(\alpha_{5}\alpha_{6})(\alpha_{8}\alpha_{11})
(\alpha_{9}\alpha_{21})(\alpha_{12}\alpha_{22})(\alpha_{17}\alpha_{23})(\alpha_{19}\alpha_{24}),
$$
$$
(\alpha_{2}\alpha_{19})(\alpha_{3}\alpha_{5})(\alpha_{6}\alpha_{10})(\alpha_{8}\alpha_{9})
(\alpha_{11}\alpha_{21})(\alpha_{12}\alpha_{23})(\alpha_{17}\alpha_{22})(\alpha_{20}\alpha_{24}),
$$
$$
(\alpha_{1}\alpha_{16})(\alpha_{2}\alpha_{20})(\alpha_{3}\alpha_{6})(\alpha_{5}\alpha_{10})
(\alpha_{12}\alpha_{23})(\alpha_{14}\alpha_{18})(\alpha_{17}\alpha_{22})(\alpha_{19}\alpha_{24}),
$$
$$
(\alpha_{1}\alpha_{14})(\alpha_{2}\alpha_{24})(\alpha_{3}\alpha_{5})(\alpha_{6}\alpha_{10})
(\alpha_{12}\alpha_{22})(\alpha_{16}\alpha_{18})(\alpha_{17}\alpha_{23})(\alpha_{19}\alpha_{20})].
$$
Here, we denote the corresponding subgroup of the symmetric group $\SSS_{24}$ on
$\alpha_{1},\dots, \alpha_{24}$ by $G$ too.

For ${\bf n}=9$, the group $G_1\cong (C_2)^3$ has the order $8$ and the identificator
$i=5$ (see Hashimoto \cite{Hash}) for the program GAP (see \cite{GAP}). The group $G_1$ is isomorphic
to ($G1:={\rm SmallGroup}(8,5);$) for the program GAP.

Using the program GAP, we find all conjugacy classes \newline (${\rm F}:={\rm IsomorphicSubgroups}(G,G1)$;) of
$G1$ in $G$. There are $15$ conjugacy classes. For the first of them
$F[1]$ (for other conjugacy classes, the result is similar),
we find the corresponding subgroup ($G_1:={\rm Image}(F[1]);)\subset
G=H_{21,1}$ and its orbits (${\rm ORB}:={\rm Orbits}(G_1);$). These orbits of order $>1$ are
$\{\alpha_1,\alpha_{14}\}$, $\{\alpha_{16},\alpha_{18}\}$,
$\{\alpha_2,\alpha_{20},\alpha_{19},\alpha_{24}\}$,
$\{\alpha_3,\alpha_{10},\alpha_{5},\alpha_6\}$,
$\{\alpha_8,\alpha_{11},\alpha_{9},\alpha_{21}\}$,
$\{\alpha_{12},\alpha_{22},\alpha_{23},\alpha_{17}\}$.

Thus, for $({\bf n}=21,4\aaa_1)$ marked by $H_{21,1}$ with the orbit
$\{\alpha_1,\alpha_{16},\alpha_{14},\alpha_{18}\}$, for the subgroup
$G_1\subset G=H_{21,1}$, the orbit splits in two suborbits
$\{\alpha_1,\alpha_{14}\}$, $\{\alpha_{16},\alpha_{18}\}$ of the order two,
and it gives
$({\bf n}=9,\ ((2\aaa_1,2\aaa_1)\subset 4\aaa_1)_{II})$ with the same
lattice.

\medskip

More exactly (by the program GAP), the corresponding subgroup $G_1\subset G$ is equal to
%gap> T1:=Image(F[1]);
%Group([ (2,19)(3,5)(6,10)(8,9)(11,21)(12,23)(17,22)(20,24), (2,20)(3,10)(5,6)(8,11)(9,21)(12,22)(17,23)(19,24), (1,14)
%(3,10)(5,6)(8,21)(9,11)(12,23)(16,18)(17,22) ])
%gap> Orbits(T1);
%[ [ 1, 14 ], [ 2, 19, 20, 24 ], [ 3, 5, 10, 6 ], [ 8, 9, 11, 21 ], [ 12, 23, 22, 17 ], [ 16, 18 ] ]
$$
G=H_{21,2}\supset G_1=
$$
$$
[(\alpha_{2}\alpha_{19})(\alpha_{3}\alpha_{5})(\alpha_{6}\alpha_{10})(\alpha_{8}\alpha_{9})
(\alpha_{11}\alpha_{21})(\alpha_{12}\alpha_{23})(\alpha_{17}\alpha_{22})(\alpha_{20}\alpha_{24}),
$$
$$
(\alpha_{2}\alpha_{20})(\alpha_{3}\alpha_{10})(\alpha_{5}\alpha_{6})(\alpha_{8}\alpha_{11})
(\alpha_{9}\alpha_{21})(\alpha_{12}\alpha_{22})(\alpha_{17}\alpha_{23})(\alpha_{19}\alpha_{24}),
$$
$$
(\alpha_{1}\alpha_{14})(\alpha_{3}\alpha_{10})(\alpha_{5}\alpha_{6})(\alpha_{8}\alpha_{21})
(\alpha_{9}\alpha_{11})(\alpha_{12}\alpha_{23})(\alpha_{16}\alpha_{18})(\alpha_{17}\alpha_{22})].
$$
Since both cases
$({\bf n}=9,\ ((2\aaa_1,2\aaa_1)\subset 4\aaa_1)_{II})$ and $({\bf n}=21,4\aaa_1)$ have a unique
lattices $S_1$ and $S$, up to isomorphism (they are marked by $\ast$ in Tables 1---4; for Table 1,
all lattices are marked by $\ast$ by definition)), they  have isomorphic lattices $S_1\cong S$ at any case: it is enough to find just one case when they are isomorphic.

\medskip

Case 2:
$({\bf n}=9,\ (8\aaa_1,8\aaa_1)\subset 16\aaa_1)\
\Longleftarrow ({\bf n}=21,16\aaa_1)$.

Similar to Case 1. By \cite{Nik10}, the $G\cong (C_2)^4$ is
marked by $N_{23}$ and
$$
G=H_{21,2}=
$$
$$
%%G2:=Group([ (1,3)(2,23)(5,14)(6,16)(10,18)(12,20)(17,24)(19,22),
%(1,2)(3,23)(5,17)(6,12)(10,22)(14,24)(16,20)(18,19),
%(1,16)(2,20)(3,6)(5,10)(12,23)(14,18)(17,22)(19,24),
%(1,14)(2,24)(3,5)(6,10)(12,22)(16,18)(17,23)(19,20) ]);
[(\alpha_{1}\alpha_{3})(\alpha_{2}\alpha_{23})
(\alpha_{5}\alpha_{14})(\alpha_{6}\alpha_{16})
(\alpha_{10}\alpha_{18})(\alpha_{12}\alpha_{20})
(\alpha_{17}\alpha_{24})(\alpha_{19}\alpha_{22}),
$$
$$
(\alpha_{1}\alpha_{2})(\alpha_{3}\alpha_{23})
(\alpha_{5}\alpha_{17})(\alpha_{6}\alpha_{12})
(\alpha_{10}\alpha_{22})(\alpha_{14}\alpha_{24})
(\alpha_{16}\alpha_{20})(\alpha_{18}\alpha_{19}),
$$
$$
(\alpha_{1}\alpha_{16})(\alpha_{2}\alpha_{20})
(\alpha_{3}\alpha_{6})(\alpha_{5}\alpha_{10})
(\alpha_{12}\alpha_{23})(\alpha_{14}\alpha_{18})
(\alpha_{17}\alpha_{22})(\alpha_{19}\alpha_{24}),
$$
$$
(\alpha_{1}\alpha_{14})(\alpha_{2}\alpha_{24})
(\alpha_{3}\alpha_{5})(\alpha_{6}\alpha_{10})
(\alpha_{12}\alpha_{22})(\alpha_{16}\alpha_{18})
(\alpha_{17}\alpha_{23})(\alpha_{19}\alpha_{20})]
$$
with the orbit
%%[ [ 1, 3, 2, 16, 14, 23, 6, 5, 20, 24, 18, 12, 17, 10, 19, 22 ] ]
$
\{\alpha_{1},\alpha_{3},\alpha_{2},\alpha_{16},\alpha_{14},
\alpha_{23},\alpha_{6},\alpha_{5},
\alpha_{20},\alpha_{24},\alpha_{18},\alpha_{12},\alpha_{17},
\alpha_{10},\alpha_{19},\alpha_{22}\}
$. The $G_1\cong (C_2)^3$ is marked by
$$
G=H_{21,2}\supset G_1=
$$
$$
[(\alpha_{1}\alpha_{3})(\alpha_{2}\alpha_{23})
(\alpha_{5}\alpha_{14})(\alpha_{6}\alpha_{16})
(\alpha_{10}\alpha_{18})(\alpha_{12}\alpha_{20})
(\alpha_{17}\alpha_{24})(\alpha_{19}\alpha_{22}),
$$
$$
(\alpha_{1}\alpha_{2})(\alpha_{3}\alpha_{23})
(\alpha_{5}\alpha_{17})(\alpha_{6}\alpha_{12})
(\alpha_{10}\alpha_{22})(\alpha_{14}\alpha_{24})
(\alpha_{16}\alpha_{20})(\alpha_{18}\alpha_{19}),
$$
$$
(\alpha_{1}\alpha_{16})(\alpha_{2}\alpha_{20})
(\alpha_{3}\alpha_{6})(\alpha_{5}\alpha_{10})
(\alpha_{12}\alpha_{23})(\alpha_{14}\alpha_{18})
(\alpha_{17}\alpha_{22})(\alpha_{19}\alpha_{24})]
$$
with suborbits
$\{\alpha_{1},\alpha_{3},\alpha_{2},\alpha_{16},\alpha_{23},\alpha_{6},\alpha_{20},\alpha_{12}\}$,
$\{\alpha_{5},\alpha_{14},\alpha_{17},\alpha_{10},\alpha_{24},\alpha_{18},\alpha_{22},\alpha_{19}\}$
(of the orbit of $G$ above). Both $G$ and $G_1$ are marked by $\ast$ in Tables 1---4.

\medskip

Case 3:
$({\bf n}=9,\ ((2\aaa_1,2\aaa_1)_{II},4\aaa_1)\subset 8\aaa_1)\
\Longleftarrow\  ({\bf n}=21,(4\aaa_1,4\aaa_1)\subset 8\aaa_1)$.

Similar to Case 1. By \cite{Nik11}, the $G\cong (C_2)^4$
is marked by $N_{23}$ and $G=H_{21,2}$ of Case 1 with orbits
$\{\alpha_1,\alpha_{16},\alpha_{14},\alpha_{18}\}$,
$\{\alpha_2,\alpha_{20},\alpha_{19},\alpha_{24}\}$. The
$G_1\cong (C_2)^3$ is marked by $G_1\subset H_{21,2}$
of Case 1 with suborbits
$\{\alpha_1,\alpha_{14}\}$, $\{\alpha_{16},\alpha_{18}\}$,
$\{\alpha_2,\alpha_{20},\alpha_{19},\alpha_{24}\}$.
Both $G$ and $G_1$ are marked by $\ast$ in the tables 1---4.

\medskip

Case 4:
({\bf n}=9,\
$
\left(\begin{array}{cccc}
 2\aaa_1 & (4\aaa_1)_I & (4\aaa_1)_I  & (4\aaa_1)_I \\
         & 2\aaa_1 & (4\aaa_1)_I & (4\aaa_1)_I  \\
         &         & 2\aaa_1 &    (4\aaa_1)_I \\
         &         &         & 2\aaa_1
\end{array}\right)
\subset 8\aaa_1) \Longleftarrow ({\bf n}=40,\ 8\aaa_1)$.

Similar to Case 1. By \cite{Nik10}, the $G\cong Q_8\ast Q_8$ is marked
by $N_{23}$ and
$$
G=H_{40,1}=
%%G:=Group([ (1,16)(2,18)(4,9)(5,23)(10,12)(14,20)(15,21)(19,24),
%(3,5)(6,12)(8,13)(10,22)(15,21)(16,19)(17,23)(18,20),
%(1,16)(2,20)(3,6)(5,10)(12,23)(14,18)(17,22)(19,24),
%(1,14)(2,24)(3,5)(6,10)(12,22)(16,18)(17,23)(19,20) ]);
$$
$$
[(\alpha_{1}\alpha_{16})(\alpha_{2}\alpha_{18})
(\alpha_{4}\alpha_{9})(\alpha_{5}\alpha_{23})
(\alpha_{10}\alpha_{12})(\alpha_{14}\alpha_{20})
(\alpha_{15}\alpha_{21})(\alpha_{19}\alpha_{24}),
$$
$$
(\alpha_{3}\alpha_{5})(\alpha_{6}\alpha_{12})
(\alpha_{8}\alpha_{13})(\alpha_{10}\alpha_{22})
(\alpha_{15}\alpha_{21})(\alpha_{16}\alpha_{19})
(\alpha_{17}\alpha_{23})(\alpha_{18}\alpha_{20}),
$$
$$
(\alpha_{1}\alpha_{16})(\alpha_{2}\alpha_{20})
(\alpha_{3}\alpha_{6})(\alpha_{5}\alpha_{10})
(\alpha_{12}\alpha_{23})(\alpha_{14}\alpha_{18})
(\alpha_{17}\alpha_{22})(\alpha_{19}\alpha_{24}),
$$
$$
(\alpha_{1}\alpha_{14})(\alpha_{2}\alpha_{24})
(\alpha_{3}\alpha_{5})(\alpha_{6}\alpha_{10})
(\alpha_{12}\alpha_{22})(\alpha_{16}\alpha_{18})
(\alpha_{17}\alpha_{23})(\alpha_{19}\alpha_{20})]
$$
with the orbit
%%[ [ 1, 16, 14, 19, 18, 20, 24, 2 ], [ 3, 5, 6, 23, 10, 12, 17, 22 ],
%%[ 4, 9 ], [ 8, 13 ], [ 15, 21 ] ]
$\{\alpha_{1},\alpha_{16},\alpha_{14},\alpha_{19},\alpha_{18},\alpha_{20},
\alpha_{24},\alpha_{2}\}$. The $G_1\cong (C_2)^3$ is marked by
$$
G=H_{40,1}\supset G_1=
$$
$$
[(\alpha_{1}\alpha_{24})(\alpha_{2}\alpha_{14})(\alpha_{3}\alpha_{23})(\alpha_{5}\alpha_{17})
(\alpha_{6}\alpha_{10})(\alpha_{8}\alpha_{13})(\alpha_{12}\alpha_{22})(\alpha_{15}\alpha_{21}),
$$
$$
(\alpha_{3}\alpha_{5})(\alpha_{6}\alpha_{12})(\alpha_{8}\alpha_{13})(\alpha_{10}\alpha_{22})
(\alpha_{15}\alpha_{21})(\alpha_{16}\alpha_{19})(\alpha_{17}\alpha_{23})(\alpha_{18}\alpha_{20}),
$$
$$
(\alpha_{2}\alpha_{14})(\alpha_{3}\alpha_{6})(\alpha_{4}\alpha_{9})(\alpha_{5}\alpha_{12})
(\alpha_{10}\alpha_{23})(\alpha_{15}\alpha_{21})(\alpha_{17}\alpha_{22})(\alpha_{18}\alpha_{20})]
$$
%[ (1,24)(2,14)(3,23)(5,17)(6,10)(8,13)(12,22)(15,21), (3,5)(6,12)(8,13)(10,22)(15,21)(16,19)(17,23)(18,20),
%  (2,14)(3,6)(4,9)(5,12)(10,23)(15,21)(17,22)(18,20) ]
with suborbits $\{\alpha_{1},\alpha_{24}\}$, $\{\alpha_{2},\alpha_{14}\}$,
$\{\alpha_{16},\alpha_{19}\}$, $\{\alpha_{18},\alpha_{20}\}$. Both $G$ and $G_1$
are marked by $\ast$ in Tables 1---4.

\medskip

Case 5:
$({\bf n}=9,\ ((2\aaa_1,2\aaa_1)_{II},4\aaa_1,4\aaa_1)
\subset 12\aaa_1)\Longleftarrow ({\bf n}=49,\ 12\aaa_1)$.

Similar to Case 1. By \cite{Nik10}, the $G\cong 2^4C_3$ is marked
by $N_{23}$ and
$$
G=H_{49,1}=
$$
%%G1:=Group([(1,22,19)(3,16,17)(4,20,9)(7,10,8)(12,13,23)(14,18,21),
%(2,12)(3,8)(4,20)(7,16)(9,11)(13,23)(14,22)(18,19),
%(2,13)(3,22)(4,9)(7,18)(8,14)(11,20)(12,23)(16,19)]);
$$
[(\alpha_{1}\alpha_{22}\alpha_{19})
(\alpha_{3}\alpha_{16}\alpha_{17})
(\alpha_{4}\alpha_{20}\alpha_{9})
(\alpha_{7}\alpha_{10}\alpha_{8})
(\alpha_{12}\alpha_{13}\alpha_{23})
(\alpha_{14}\alpha_{18}\alpha_{21}),
$$
$$
(\alpha_{2}\alpha_{12})(\alpha_{3}\alpha_{8})
(\alpha_{4}\alpha_{20})(\alpha_{7}\alpha_{16})
(\alpha_{9}\alpha_{11})(\alpha_{13}\alpha_{23})
(\alpha_{14}\alpha_{22})(\alpha_{18}\alpha_{19}),
$$
$$
(\alpha_{2}\alpha_{13})(\alpha_{3}\alpha_{22})
(\alpha_{4}\alpha_{9})(\alpha_{7}\alpha_{18})
(\alpha_{8}\alpha_{14})(\alpha_{11}\alpha_{20})
(\alpha_{12}\alpha_{23})(\alpha_{16}\alpha_{19})]
$$
with the orbit
$\{\alpha_{1},\alpha_{22},\alpha_{21},\alpha_{10},\alpha_{19},\alpha_{14},
\alpha_{8},\alpha_{17},\alpha_{18},\alpha_{7},\alpha_{3},\alpha_{16}\}$.
%gap> H1:=Image(F[1]);
%Group([ (1,21)(2,23)(3,8)(4,9)(10,17)(11,20)(12,13)(14,22), (2,12)(3,8)(4,20)(7,16)(9,11)(13,23)(14,22)(18,19), (2,13)
%(3,22)(4,9)(7,18)(8,14)(11,20)(12,23)(16,19) ])
%gap> Orbits(G);
%[ [ 1, 22, 19, 14, 3, 18, 16, 8, 21, 7, 17, 10 ], [ 2, 12, 13, 23 ], [ 4, 20, 9, 11 ] ]
%gap>
%gap> Orbits(H1);
%[ [ 1, 21 ], [ 2, 23, 12, 13 ], [ 3, 8, 22, 14 ], [ 4, 9, 20, 11 ], [ 7, 16, 18, 19 ], [ 10, 17 ] ]
The $G_1\cong (C_2)^3$ is marked by
$$
G=H_{49,1}\supset G_1=
$$
$$
[(\alpha_{1}\alpha_{21})(\alpha_{2}\alpha_{23})(\alpha_{3}\alpha_{8})(\alpha_{4}\alpha_{9})
(\alpha_{10}\alpha_{17})(\alpha_{11}\alpha_{20})(\alpha_{12}\alpha_{13})(\alpha_{14}\alpha_{22}),
$$
$$
(\alpha_{2}\alpha_{12})(\alpha_{3}\alpha_{8})(\alpha_{4}\alpha_{20})(\alpha_{7}\alpha_{16})
(\alpha_{9}\alpha_{11})(\alpha_{13}\alpha_{23})(\alpha_{14}\alpha_{22})(\alpha_{18}\alpha_{19}),
$$
$$
(\alpha_{2}\alpha_{13})(\alpha_{3}\alpha_{22})(\alpha_{4}\alpha_{9})(\alpha_{7}\alpha_{18})
(\alpha_{8}\alpha_{14})(\alpha_{11}\alpha_{20})(\alpha_{12}\alpha_{23})(\alpha_{16}\alpha_{19})]
$$
with suborbits
$\{\alpha_{1},\alpha_{21}\}$, $\{\alpha_{10},\alpha_{17}\}$,
$\{\alpha_{3},\alpha_{8},\alpha_{22},\alpha_{14}\}$, $\{\alpha_{7},\alpha_{16},\alpha_{18},\alpha_{19}\}$.
Both $G$ and $G_1$ are marked by $\ast$ in Tables 1---4.

\medskip

Case 6:
({\bf n}=9,\
$
\left(\begin{array}{cccc}
 2\aaa_1 & 2\aaa_3     & (4\aaa_1)_I & 6\aaa_1 \\
         & 4\aaa_1     & 6\aaa_1     & 8\aaa_1 \\
         &             & 2\aaa_1     & 2\aaa_3 \\
         &             &              & 4\aaa_1
\end{array}\right)
\subset 4\aaa_3)\
\Longleftarrow\ ({\bf n}=22,(4\aaa_1,8\aaa_1)\subset 4\aaa_3)$.

Similar to Case 1. By \cite{Nik11}, the $G\cong C_2\times D_8$ is marked
by $N_{21}$ and
%G:=Group([ (7,9)(10,15)(11,14)(12,13)(16,22)(17,23)(18,24)(19,21),(4,6)(7,18)(8,17)(9,16)(10,12)(19,22)(20,23)(21,24),
%(1,3)(4,6)(7,9)(10,12)(13,15)(16,18)(19,21)(22,24) ]);
$$
G=H_{22,1}=
$$
$$
[(\alpha_{1,3}\alpha_{3,3})(\alpha_{1,4}\alpha_{3,5})
(\alpha_{2,4}\alpha_{2,5})(\alpha_{3,4}\alpha_{1,5})
(\alpha_{1,6}\alpha_{1,8})
(\alpha_{2,6}\alpha_{2,8})
(\alpha_{3,6}\alpha_{3,8})(\alpha_{1,7}\alpha_{3,7}),
$$
$$
(\alpha_{1,2}\alpha_{3,2})(\alpha_{1,3}\alpha_{3,6})
(\alpha_{2,3}\alpha_{2,6})(\alpha_{3,3}\alpha_{1,6})
(\alpha_{1,4}\alpha_{3,4})(\alpha_{1,7}\alpha_{1,8})
(\alpha_{2,7}\alpha_{2,8})(\alpha_{3,7}\alpha_{3,8}),
$$
$$
(\alpha_{1,1}\alpha_{3,1})(\alpha_{1,2}\alpha_{3,2})
(\alpha_{1,3}\alpha_{3,3})(\alpha_{1,4}\alpha_{3,4})
(\alpha_{1,5}\alpha_{3,5})(\alpha_{1,6}\alpha_{3,6})
(\alpha_{1,7}\alpha_{3,7})(\alpha_{1,8}\alpha_{3,8}) ]
$$
with orbits
%gap> Orbits(H4);
%[ 7, 9, 18, 16, 24, 22, 21, 19 ], [ 8, 17, 23, 20 ],
$\{\alpha_{2,3},\alpha_{2,6},\alpha_{2,8},\alpha_{2,7}\}$,
$\{\alpha_{1,3},\alpha_{3,3},\alpha_{3,6},\alpha_{1,6},\alpha_{3,8},
\alpha_{1,8},\alpha_{3,7},\alpha_{1,7}\}$.
The $G_1\cong (C_2)^3$ is marked by
$$
G=H_{22,1}\supset G_1=
$$
%gap> H1:=Image(F[1]);
%Group([ (7,21)(8,20)(9,19)(10,13)(11,14)(12,15)(16,18)(22,24), (7,9)(10,15)(11,14)(12,13)(16,22)(17,23)(18,24)
%(19,21), (1,3)(4,6)(10,13)(11,14)(12,15)(16,24)(17,23)(18,22) ])
$$
[(\alpha_{1,3}\alpha_{3,7})(\alpha_{2,3}\alpha_{2,7})(\alpha_{3,3}\alpha_{1,7})(\alpha_{1,4}\alpha_{1,5})
(\alpha_{2,4}\alpha_{2,5})(\alpha_{3,4}\alpha_{3,5})(\alpha_{1,6}\alpha_{3,6})(\alpha_{1,8}\alpha_{3,8}),
$$
$$
(\alpha_{1,3}\alpha_{3,3})(\alpha_{1,4}\alpha_{3,5})(\alpha_{2,4}\alpha_{2,5})(\alpha_{3,4}\alpha_{1,5})
(\alpha_{1,6}\alpha_{1,8})(\alpha_{2,6}\alpha_{2,8})(\alpha_{3,6}\alpha_{3,8})(\alpha_{1,7}\alpha_{3,7}),
$$
$$
(\alpha_{1,1}\alpha_{3,1})(\alpha_{1,2}\alpha_{3,2})(\alpha_{1,4}\alpha_{1,5})(\alpha_{2,4}\alpha_{2,5})
(\alpha_{3,4}\alpha_{3,5})(\alpha_{1,6}\alpha_{3,8})(\alpha_{2,6}\alpha_{2,8})(\alpha_{3,6}\alpha_{1,8})]
$$
with suborbits
$\{\alpha_{2,3},\alpha_{2,7}\}$, $\{\alpha_{1,3},\alpha_{3,7},\alpha_{3,3},\alpha_{1,7}\}$,
$\{\alpha_{2,6},\alpha_{2,8}\}$,
$\{\alpha_{1,6},\alpha_{3,6},\alpha_{1,8},\alpha_{3,8}\}$.
Both $G$ and $G_1$ are marked by $\ast$ in Tables 1---4.

\medskip

Case 7:
$({\bf n}=9,\ (4\aaa_1,4\aaa_1,4\aaa_1,4\aaa_1)
\subset 16\aaa_1)\Longleftarrow ({\bf n}=39,\ 16\aaa_1)$.

Similar to Case 1. By \cite{Nik10}, the $G\cong 2^4C_2$ is marked
by $N_{23}$ and
%%G1:=Group([ (5,10)(7,13)(8,12)(9,22)(11,23)(14,16)(17,21)(19,20),
%(2,12)(3,8)(4,20)(7,16)(9,11)(13,23)(14,22)(18,19),
%(2,13)(3,22)(4,9)(7,18)(8,14)(11,20)(12,23)(16,19) ]);
$$
G=H_{39,1}=
$$
$$
[(\alpha_{5}\alpha_{10})(\alpha_{7}\alpha_{13})
(\alpha_{8}\alpha_{12})(\alpha_{9}\alpha_{22})
(\alpha_{11}\alpha_{23})(\alpha_{14}\alpha_{16})
(\alpha_{17}\alpha_{21})(\alpha_{19}\alpha_{20}),
$$
$$
(\alpha_{2}\alpha_{12})(\alpha_{3}\alpha_{8})
(\alpha_{4}\alpha_{20})(\alpha_{7}\alpha_{16})
(\alpha_{9}\alpha_{11})(\alpha_{13}\alpha_{23})
(\alpha_{14}\alpha_{22})(\alpha_{18}\alpha_{19}),
$$
$$
(\alpha_{2}\alpha_{13})(\alpha_{3}\alpha_{22})
(\alpha_{4}\alpha_{9})(\alpha_{7}\alpha_{18})
(\alpha_{8}\alpha_{14})(\alpha_{11}\alpha_{20})
(\alpha_{12}\alpha_{23})(\alpha_{16}\alpha_{19})]
$$
with the orbit
$\{\alpha_{2},\alpha_{12},\alpha_{13},\alpha_{3},\alpha_{18},\alpha_{8},
\alpha_{23},\alpha_{19},\alpha_{7},\alpha_{22},\alpha_{4},
\alpha_{14},\alpha_{20},\alpha_{11},\alpha_{16},\alpha_{9}\}$.
The $G_1\cong (C_2)^3$ is marked by
$$
G=H_{39,1}\supset G_1=
$$
%gap> H1:=Image(F[1]);
%Group([ (2,3)(4,18)(5,10)(7,22)(9,13)(11,14)(16,23)(17,21), (2,18)(3,4)(5,10)(8,19)(11,16)(12,20)(14,23)
%(17,21), (5,10)(7,13)(8,12)(9,22)(11,23)(14,16)(17,21)(19,20) ])
$$
[(\alpha_{2}\alpha_{3})(\alpha_{4}\alpha_{18})(\alpha_{5}\alpha_{10})(\alpha_{7}\alpha_{22})
(\alpha_{9}\alpha_{13})(\alpha_{11}\alpha_{14})(\alpha_{16}\alpha_{23})(\alpha_{17}\alpha_{21}),
$$
$$
(\alpha_{2}\alpha_{18})(\alpha_{3}\alpha_{4})(\alpha_{5}\alpha_{10})(\alpha_{8}\alpha_{19})
(\alpha_{11}\alpha_{16})(\alpha_{12}\alpha_{20})(\alpha_{14}\alpha_{23})(\alpha_{17}\alpha_{21}),
$$
$$
(\alpha_{5}\alpha_{10})(\alpha_{7}\alpha_{13})(\alpha_{8}\alpha_{12})(\alpha_{9}\alpha_{22})
(\alpha_{11}\alpha_{23})(\alpha_{14}\alpha_{16})(\alpha_{17}\alpha_{21})(\alpha_{19}\alpha_{20})]
$$
with suborbits
$\{\alpha_{2},\alpha_{3},\alpha_{18},\alpha_{4}\}$, $\{\alpha_{7},\alpha_{22},\alpha_{13},\alpha_{9}\}$,
$\{\alpha_{8},\alpha_{19},\alpha_{12},\alpha_{20}\}$, $\{\alpha_{11},\alpha_{14},\alpha_{16},\alpha_{23}\}$.
\newline
%%gap> Orbits(H1);
%%[ [ 2, 3, 18, 4 ], [ 5, 10 ], [ 7, 22, 13, 9 ], [ 8, 19, 12, 20 ], [ 11, 14, 16, 23 ], [ 17, 21 ] ]
Both $G$ and $G_1$ are marked by $\ast$ in Tables 1---4.

\medskip

Case 8:
({\bf n}=9,\
$
\left(\begin{array}{ccccc}
 2\aaa_1 & (4\aaa_1)_I & (4\aaa_1)_I  & (4\aaa_1)_I &  10\aaa_1 \\
         & 2\aaa_1     & (4\aaa_1)_I  & (4\aaa_1)_I &  10\aaa_1 \\
         &             & 2\aaa_1      & (4\aaa_1)_I &  10\aaa_1 \\
         &             &              & 2\aaa_1     &  10\aaa_1\\
         &             &              &             &   8\aaa_1
\end{array}\right)
\subset 16\aaa_1 ) \Longleftarrow ({\bf n}=56,\ 16\aaa_1)$.

Similar to Case 1. By \cite{Nik10}, the $G\cong \Gamma_{25}a_1$ is marked
by $N_{23}$ and
$$
G=H_{56,2}=
$$
%G2:=Group([(2,3)(5,6)(7,18)(8,23)(10,20)(11,17)(15,16)(19,24),
%(1,14)(2,23)(3,5)(6,18)(7,8)(9,13)(10,16)(17,24),
%(1,14)(2,24)(3,5)(6,10)(12,22)(16,18)(17,23)(19,20)]);
$$
[(\alpha_{2}\alpha_{3})(\alpha_{5}\alpha_{6})(\alpha_{7}\alpha_{18})
(\alpha_{8}\alpha_{23})(\alpha_{10}\alpha_{20})(\alpha_{11}\alpha_{17})
(\alpha_{15}\alpha_{16})(\alpha_{19}\alpha_{24}),
$$
$$
(\alpha_{1}\alpha_{14})(\alpha_{2}\alpha_{23})(\alpha_{3}\alpha_{5})
(\alpha_{6}\alpha_{18})(\alpha_{7}\alpha_{8})(\alpha_{9}\alpha_{13})
(\alpha_{10}\alpha_{16})(\alpha_{17}\alpha_{24}),
$$
$$
(\alpha_{1}\alpha_{14})(\alpha_{2}\alpha_{24})(\alpha_{3}\alpha_{5})
(\alpha_{6}\alpha_{10})(\alpha_{12}\alpha_{22})(\alpha_{16}\alpha_{18})
(\alpha_{17}\alpha_{23})(\alpha_{19}\alpha_{20})]
$$
with the orbit
$\{\alpha_{2}, \alpha_{3}, \alpha_{23},
\alpha_{24}, \alpha_{5},\alpha_{8}, \alpha_{17}, \alpha_{19}, \alpha_{6},
\alpha_{7}, \alpha_{11}, \alpha_{20}, \alpha_{18},
\alpha_{10},\alpha_{16},\alpha_{15}\}$.
The $G_1\cong (C_2)^3$ is marked by
$$
G=H_{56,2}\supset G_1=
$$
%gap> H2:=Image(F[2]);
%Group([ (2,6)(3,15)(5,11)(7,19)(8,20)(9,13)(12,22)(16,17), (3,11)(5,15)(7,20)(8,19)(9,13)(10,24)(12,22)(18,23), (1,14)
%(3,7)(5,8)(9,13)(10,24)(11,20)(15,19)(16,17) ])
$$
[(\alpha_{2}\alpha_{6})(\alpha_{3}\alpha_{15})(\alpha_{5}\alpha_{11})(\alpha_{7}\alpha_{19})
(\alpha_{8}\alpha_{20})(\alpha_{9}\alpha_{13})(\alpha_{12}\alpha_{22})(\alpha_{16}\alpha_{17}),
$$
$$
(\alpha_{3}\alpha_{11})(\alpha_{5}\alpha_{15})(\alpha_{7}\alpha_{20})(\alpha_{8}\alpha_{19})
(\alpha_{9}\alpha_{13})(\alpha_{10}\alpha_{24})(\alpha_{12}\alpha_{22})(\alpha_{18}\alpha_{23}),
$$
$$
(\alpha_{1}\alpha_{14})(\alpha_{3}\alpha_{7})(\alpha_{5}\alpha_{8})(\alpha_{9}\alpha_{13})
(\alpha_{10}\alpha_{24})(\alpha_{11}\alpha_{20})(\alpha_{15}\alpha_{19})(\alpha_{16}\alpha_{17})]
$$
with suborbits
%gap> Orbits(H2);
%[ [ 1, 14 ], [ 2, 6 ], [ 3, 15, 11, 7, 5, 19, 20, 8 ], [ 9, 13 ], [ 10, 24 ], [ 12, 22 ], [ 16, 17 ], [ 18, 23 ] ]
$\{\alpha_{2},\alpha_{6}\}$, $\{\alpha_{10},\alpha_{24}\}$, $\{\alpha_{16},\alpha_{17}\}$, $\{\alpha_{18},\alpha_{23}\}$,
$\{\alpha_{3},\alpha_{15},\alpha_{11},\alpha_{7},\alpha_{5},\alpha_{19},\alpha_{20},\alpha_{8}\}$.
Both $G$ and $G_1$ are marked by $\ast$ Tables 1---4.

\medskip

Case 9:
({\bf n}=9,\
$
\left(\begin{array}{ccccc}
 2\aaa_1 & (4\aaa_1)_I & 6\aaa_1    & (4\aaa_1)_I     &  6\aaa_1 \\
         & 2\aaa_1     & 2\aaa_3  & (4\aaa_1)_I       &  6\aaa_1 \\
         &             & 4\aaa_1      & 6\aaa_1   &  8\aaa_1 \\
         &             &              & 2\aaa_1     & 2\aaa_3\\
         &             &              &             & 4\aaa_1
\end{array}\right)
\subset 2\aaa_1\amalg 4\aaa_3)
$

$\Longleftarrow\ ({\bf n}=22,\
\left(\begin{array}{rrr}
 2\aaa_1 & 6\aaa_1 & 10\aaa_1  \\
         & 4\aaa_1 & 4\aaa_3   \\
         &         & 8\aaa_1
\end{array}\right)\subset 2\aaa_1\amalg 4\aaa_3)$.

Similar to Cases 1 and 6. By \cite{Nik11}, the $G\cong C_2\times D_8$
%G:=Group([ (7,9)(10,15)(11,14)(12,13)(16,22)(17,23)(18,24)(19,21),(4,6)(7,18)(8,17)(9,16)(10,12)(19,22)(20,23)(21,24),
%(1,3)(4,6)(7,9)(10,12)(13,15)(16,18)(19,21)(22,24) ]);
%gap> Orbits(G);
%[ [ 1, 3 ], [ 4, 6 ], [ 7, 9, 18, 16, 24, 22, 21, 19 ], [ 8, 17, 23, 20 ], [ 10, 15, 12, 13 ], [ 11, 14 ] ]
is marked by $N_{21}$ and $G=H_{22,1}$ of Case 6 with orbits
$\{\alpha_{1,1},\alpha_{3,1}\}$,
$\{\alpha_{2,3},\alpha_{2,6},\alpha_{2,8},\alpha_{2,7}\}$,
$\{\alpha_{1,3},\alpha_{3,3},\alpha_{3,6},\alpha_{1,6},
\alpha_{3,8}, \alpha_{1,8},\alpha_{3,7},\alpha_{1,7}\}$.
%gap> H1:=Image(F[1]);
%Group([ (7,21)(8,20)(9,19)(10,13)(11,14)(12,15)(16,18)(22,24), (7,9)(10,15)(11,14)(12,13)(16,22)(17,23)(18,24)
%(19,21), (1,3)(4,6)(10,13)(11,14)(12,15)(16,24)(17,23)(18,22) ])
%gap> Orbits(H1);
%[ [ 1, 3 ], [ 4, 6 ], [ 7, 21, 9, 19 ], [ 8, 20 ], [ 10, 13, 15, 12 ], [ 11, 14 ], [ 16, 18, 22, 24 ], [ 17, 23 ] ]
The $G_1\cong (C_2)^3$ is marked by
$$
G=H_{22,1}\supset G_1=
$$
$$
[(\alpha_{1,3}\alpha_{3,7})(\alpha_{2,3}\alpha_{2,7})(\alpha_{3,3}\alpha_{1,7})(\alpha_{1,4}\alpha_{1,5})
(\alpha_{2,4}\alpha_{2,5})(\alpha_{3,4}\alpha_{3,5})(\alpha_{1,6}\alpha_{3,6})(\alpha_{1,8}\alpha_{3,8}),
$$
$$
(\alpha_{1,3}\alpha_{3,3})(\alpha_{1,4}\alpha_{3,5})(\alpha_{2,4}\alpha_{2,5})(\alpha_{3,4}\alpha_{1,5})
(\alpha_{1,6}\alpha_{1,8})(\alpha_{2,6}\alpha_{2,8})(\alpha_{3,6}\alpha_{3,8})(\alpha_{1,7}\alpha_{3,7}),
$$
$$
(\alpha_{1,1}\alpha_{3,1})(\alpha_{1,2}\alpha_{3,2})(\alpha_{1,4}\alpha_{1,5})(\alpha_{2,4}\alpha_{2,5})
(\alpha_{3,4}\alpha_{3,5})(\alpha_{1,6}\alpha_{3,8})(\alpha_{2,6}\alpha_{2,8})(\alpha_{3,6}\alpha_{1,8})]
$$
with suborbits
$\{\alpha_{1,1},\alpha_{3,1}\}$,
$\{\alpha_{2,3},\alpha_{2,7}\}$, $\{\alpha_{1,3},\alpha_{3,3}\alpha_{1,7},\alpha_{3,7}\}$,
$\{\alpha_{2,6},\alpha_{2,8}\}$, $\{\alpha_{1,6},\alpha_{3,6},\alpha_{1,8},\alpha_{3,8}\}$.
Both $G$ and $G_1$
are marked by $\ast$ in Tables 1---4.

\medskip

Case 10:
$({\bf n}=9,\ ((2\aaa_1,2\aaa_1)_{II},4\aaa_1,4\aaa_1,4\aaa_1)
\subset 16\aaa_1)\Longleftarrow ({\bf n}=75,\ 16\aaa_1)$.

Similar to Case 1. By \cite{Nik10}, the $G\cong 4^2\AAA_4$ is marked
by $N_{23}$ and
$$
G=H_{75,1}=
$$
%G:=Group([(3,16)(4,5)(6,21)(10,20)(11,12)(13,17)(14,22)(23,24),
%(1,9,7)(3,21,22)(4,14,24)(5,10,16)(6,20,23)(12,17,13),
%(1,6,22)(3,24,16)(4,20,9)(5,10,7)(11,13,12)(14,18,21)]);
$$
[(\alpha_{3}\alpha_{16})(\alpha_{4}\alpha_{5})
(\alpha_{6}\alpha_{21})(\alpha_{10}\alpha_{20})
(\alpha_{11}\alpha_{12})(\alpha_{13}\alpha_{17})
(\alpha_{14}\alpha_{22})(\alpha_{23}\alpha_{24}),\
$$
$$
(\alpha_{1}\alpha_{9}\alpha_{7})
(\alpha_{3}\alpha_{21}\alpha_{22})
(\alpha_{4}\alpha_{14}\alpha_{24})
(\alpha_{5}\alpha_{10}\alpha_{16})
(\alpha_{6}\alpha_{20}\alpha_{23})
(\alpha_{12}\alpha_{17}\alpha_{13}),\
$$
$$
(\alpha_{1}\alpha_{6}\alpha_{22})
(\alpha_{3}\alpha_{24}\alpha_{16})
(\alpha_{4}\alpha_{20}\alpha_{9})
(\alpha_{5}\alpha_{10}\alpha_{7})
(\alpha_{11}\alpha_{13}\alpha_{12})
(\alpha_{14}\alpha_{18}\alpha_{21})]
$$
with the orbit
$\{\alpha_{1},\alpha_{9},\alpha_{6},\alpha_{7},\alpha_{4},
\alpha_{21},\alpha_{20},
\alpha_{22},\alpha_{5},\alpha_{14},\alpha_{10},\alpha_{23},
\alpha_{3},\alpha_{24},\alpha_{18},\alpha_{16}\}$.
The $G_1\cong (C_2)^3$ is marked by
$$
G=H_{75,1}\supset G_1=
$$
%gap> H1:=Image(F[1]);
%Group([ (1,7)(3,24)(4,21)(5,6)(9,18)(11,12)(13,17)(16,23), (3,24)(4,6)(5,21)(10,14)(11,17)(12,13)(16,23)
%(20,22), (3,16)(4,5)(6,21)(10,20)(11,12)(13,17)(14,22)(23,24) ])
$$
[(\alpha_{1}\alpha_{7})(\alpha_{3}\alpha_{24})(\alpha_{4}\alpha_{21})(\alpha_{5}\alpha_{6})
(\alpha_{9}\alpha_{18})(\alpha_{11}\alpha_{12})(\alpha_{13}\alpha_{17})(\alpha_{16}\alpha_{23}),
$$
$$
(\alpha_{3}\alpha_{24})(\alpha_{4}\alpha_{6})(\alpha_{5}\alpha_{21})(\alpha_{10}\alpha_{14})
(\alpha_{11}\alpha_{17})(\alpha_{12}\alpha_{13})(\alpha_{16}\alpha_{23})(\alpha_{20}\alpha_{22}),
$$
$$
(\alpha_{3}\alpha_{16})(\alpha_{4}\alpha_{5})(\alpha_{6}\alpha_{21})(\alpha_{10}\alpha_{20})
(\alpha_{11}\alpha_{12})(\alpha_{13}\alpha_{17})(\alpha_{14}\alpha_{22})(\alpha_{23}\alpha_{24})]
$$
with suborbits
%gap> Orbits(H1);
%[ [ 1, 7 ], [ 9, 18 ], [ 3, 24, 16, 23 ], [ 4, 21, 6, 5 ], [ 10, 14, 20, 22 ]]
$\{\alpha_{1},\alpha_{7}\}$, $\{\alpha_{9},\alpha_{18}\}$,
$\{\alpha_{3},\alpha_{24},\alpha_{16},\alpha_{23}\}$, $\{\alpha_{4},\alpha_{21},\alpha_{6},\alpha_{5}\}$,
$\{\alpha_{10},\alpha_{14},\alpha_{20},\alpha_{22}\}$.
Both $G$ and $G_1$ are marked by $\ast$ in Tables 1---4.

\medskip

Case 11:
$({\bf n}=9,
\left(\begin{array}{ccccc}
 2\aaa_1 & 2\aaa_3     & (4\aaa_1)_I  &  6\aaa_1     &  6\aaa_1 \\
         & 4\aaa_1     & 6\aaa_1     &  8\aaa_1     &  8\aaa_1 \\
         &             & 2\aaa_1      & 2\aaa_3     &  6\aaa_1 \\
         &             &              & 4\aaa_1     &  8\aaa_1\\
         &             &              &             &  4\aaa_1
\end{array}\right)
\subset 4\aaa_3\amalg  4\aaa_1)
$

$\Longleftarrow
({\bf n}=22,\
\left(\begin{array}{rrr}
 4\aaa_1 & (8\aaa_1)_{II} & 12\aaa_1  \\
         & 4\aaa_1 & 4\aaa_3   \\
         &         & 8\aaa_1
\end{array}\right)\subset 4\aaa_1\amalg 4\aaa_3)$.

%G:=Group([ (7,9)(10,15)(11,14)(12,13)(16,22)(17,23)(18,24)(19,21),(4,6)(7,18)(8,17)(9,16)(10,12)(19,22)(20,23)(21,24),
%(1,3)(4,6)(7,9)(10,12)(13,15)(16,18)(19,21)(22,24) ]);
%gap> Orbits(H4);
%[ [ 1, 3 ], [ 4, 6 ], [ 7, 9, 18, 16, 24, 22, 21, 19 ], [ 8, 17, 23, 20 ], [ 10, 15, 12, 13 ], [ 11, 14 ] ]

Similar to Cases 1 and 6. By \cite{Nik11}, the $G\cong C_2\times D_8$ is marked
by $N_{21}$ and $G=H_{22,1}$ from Case 6
with orbits
%[ 10, 15, 12, 13 ], [ 8, 17, 23, 20 ],[ 7, 9, 18, 16, 24, 22, 21, 19 ],
$\{\alpha_{1,4},\alpha_{3,5},\alpha_{3,4},\alpha_{1,5}\}$,
$\{\alpha_{2,3},\alpha_{2,6},\alpha_{2,8},\alpha_{2,7}\}$,
$\{\alpha_{1,3},\alpha_{3,3},\alpha_{3,6},\alpha_{1,6},\alpha_{3,8},
\linebreak \alpha_{1,8},\alpha_{3,7},\alpha_{1,7}\}$.
The $G_1\cong (C_2)^3$ is marked by
$G=H_{22,1}\supset G_1$ of Case 6
%gap> H1:=Image(F[1]);
%Group([ (7,21)(8,20)(9,19)(10,13)(11,14)(12,15)(16,18)(22,24), (7,9)(10,15)(11,14)(12,13)(16,22)(17,23)(18,24)
%(19,21), (1,3)(4,6)(10,13)(11,14)(12,15)(16,24)(17,23)(18,22) ])
%gap> Orbits(H1);
%[ [ 1, 3 ], [ 4, 6 ], [ 7, 21, 9, 19 ], [ 8, 20 ], [ 10, 13, 15, 12 ], [ 11, 14 ], [ 16, 18, 22, 24 ], [ 17, 23 ] ]
with suborbits
$\{\alpha_{2,3},\alpha_{2,7}\}$, $\{\alpha_{1,3},\alpha_{3,7},\alpha_{3,3},\alpha_{1,7}\}$,
$\{\alpha_{2,6},\alpha_{2,8}\}$, $\{\alpha_{1,6},\alpha_{3,6},\alpha_{1,8},\alpha_{3,8}\}$, $\{\alpha_{1,4},\alpha_{1,5},\alpha_{3,5},\alpha_{3,4}\}$.
Both $G$ and $G_1$ are marked by $\ast$ in Tables 1---4.

\medskip

Case 12: ({\bf n}=10, $(\aaa_1,\aaa_1)\subset 2\aaa_1)\ \Longleftarrow\ ({\bf n}=22,2\aaa_1)$.

Similar to Case 1. By \cite{Nik10}, the $G\cong C_2\times D_8$ is marked
by $N_{23}$ and
$$
G=H_{22,2}=
$$
%G2:=Group([(3,4)(6,15)(8,11)(9,22)(12,23)(14,20)(16,19)(17,21),
%(2,22)(3,13)(4,18)(7,9)(10,24)(11,20)(15,17)(16,19),
%(2,13)(3,22)(4,9)(7,18)(8,14)(11,20)(12,23)(16,19)]);
$$
[
(\alpha_{3}\alpha_{4})(\alpha_{6}\alpha_{15})
(\alpha_{8}\alpha_{11})(\alpha_{9}\alpha_{22})
(\alpha_{12}\alpha_{23})(\alpha_{14}\alpha_{20})
(\alpha_{16}\alpha_{19})(\alpha_{17}\alpha_{21}),
$$
$$
(\alpha_{2}\alpha_{22})(\alpha_{3}\alpha_{13})
(\alpha_{4}\alpha_{18})(\alpha_{7}\alpha_{9})
(\alpha_{10}\alpha_{24})(\alpha_{11}\alpha_{20})
(\alpha_{15}\alpha_{17})(\alpha_{16}\alpha_{19}),
$$
$$
(\alpha_{2}\alpha_{13})(\alpha_{3}\alpha_{22})
(\alpha_{4}\alpha_{9})(\alpha_{7}\alpha_{18})
(\alpha_{8}\alpha_{14})(\alpha_{11}\alpha_{20})
(\alpha_{12}\alpha_{23})(\alpha_{16}\alpha_{19})]
$$
with the orbit
$\{\alpha_{12},\alpha_{23}\}$.
%gap> H4:=Image(F[4]);
%Group([ (2,9)(3,18)(4,13)(6,21)(7,22)(8,14)(10,24)(16,19), (2,13)(3,9)(4,22)(6,15)(7,18)(8,20)(11,14)(17,21) ])
%gap> Orbits(H4);
%[ [ 2, 9, 13, 3, 4, 18, 22, 7 ], [ 6, 21, 15, 17 ], [ 8, 14, 20, 11 ], [ 10, 24 ], [ 16, 19 ] ]
The $G_1\cong D_8$ is marked by
$$
G=H_{22,2}\supset G_1=
$$
$$
[(\alpha_{2}\alpha_{9})(\alpha_{3}\alpha_{18})(\alpha_{4}\alpha_{13})(\alpha_{6}\alpha_{21})
(\alpha_{7}\alpha_{22})(\alpha_{8}\alpha_{14})(\alpha_{10}\alpha_{24})(\alpha_{16}\alpha_{19}),
$$
$$
(\alpha_{2}\alpha_{13})(\alpha_{3}\alpha_{9})(\alpha_{4}\alpha_{22})(\alpha_{6}\alpha_{15})
(\alpha_{7}\alpha_{18})(\alpha_{8}\alpha_{20})(\alpha_{11}\alpha_{14})(\alpha_{17}\alpha_{21})]
$$
with suborbits
$\{\alpha_{12}\}$, $\{\alpha_{23}\}$.
Both $G$ and $G_1$
are marked by $\ast$ in Tables 1---4.

\medskip

Case 13:
({\bf n}=10,
$((4\aaa_1,4\aaa_1)\subset 8\aaa_1)_{II})
\Longleftarrow\ ({\bf n}=22,8\aaa_1)$.

Similar to Cases 1 and 12. By \cite{Nik10}, the $G\cong C_2\times D_8$ is marked
by $N_{23}$ and $G=H_{22,2}$ from Case 12 with the orbit
$\{\alpha_{2},\alpha_{22},\alpha_{13},\alpha_{7},\alpha_{9},
\alpha_{3},\alpha_{18},\alpha_{4}\}$.
The $G_1\cong D_8$ is marked by
$$
G=H_{22,2}\supset G_1=
$$
%gap> H1:=Image(F[1]);
%Group([ (3,4)(6,15)(8,11)(9,22)(12,23)(14,20)(16,19)(17,21), (2,3)(4,7)(8,14)(9,18)(10,24)(12,23)(13,22)(15,17) ])
%gap> Orbits(H1);
%[ [ 2, 3, 4, 7 ], [ 6, 15, 17, 21 ], [ 8, 11, 14, 20 ], [ 9, 22, 18, 13 ], [ 10, 24 ], [ 12, 23 ], [ 16, 19 ] ]
$$
[(\alpha_{3}\alpha_{4})(\alpha_{6}\alpha_{15})(\alpha_{8}\alpha_{11})(\alpha_{9}\alpha_{22})
(\alpha_{12}\alpha_{23})(\alpha_{14}\alpha_{20})(\alpha_{16}\alpha_{19})(\alpha_{17}\alpha_{21}),
$$
$$
(\alpha_{2}\alpha_{3})(\alpha_{4}\alpha_{7})(\alpha_{8}\alpha_{14})(\alpha_{9}\alpha_{18})
(\alpha_{10}\alpha_{24})(\alpha_{12}\alpha_{23})(\alpha_{13}\alpha_{22})(\alpha_{15}\alpha_{17})]
$$
with suborbits
$\{\alpha_{2},\alpha_{3},\alpha_{4},\alpha_{7}\}$,
$\{\alpha_{9},\alpha_{22},\alpha_{18},\alpha_{13}\}$.
Both $G$ and $G_1$
are marked by $\ast$ in Tables 1---4.

\medskip

Case 14:
({\bf n}=10,
$(\aaa_1,\aaa_1,(2\aaa_1)_I)\subset 4\aaa_1)\ \Longleftarrow ({\bf n}=39,\ 4\aaa_1)$.

Similar to Case 1. By \cite{Nik10}, the $G\cong 2^4C_2$ is marked
by $N_{23}$ and
%%G2:=Group([(5,18)(6,7)(10,17)(11,20)(12,13)(14,22)(15,16)(19,24),
%(2,12)(3,8)(4,20)(7,16)(9,11)(13,23)(14,22)(18,19),
%(2,13)(3,22)(4,9)(7,18)(8,14)(11,20)(12,23)(16,19)]);
$$
G=H_{39,2}=
$$
$$
[(\alpha_{5}\alpha_{18})(\alpha_{6}\alpha_{7})
(\alpha_{10}\alpha_{17})(\alpha_{11}\alpha_{20})
(\alpha_{12}\alpha_{13})(\alpha_{14}\alpha_{22})
(\alpha_{15}\alpha_{16})(\alpha_{19}\alpha_{24}),
$$
$$
(\alpha_{2}\alpha_{12})(\alpha_{3}\alpha_{8})
(\alpha_{4}\alpha_{20})(\alpha_{7}\alpha_{16})
(\alpha_{9}\alpha_{11})(\alpha_{13}\alpha_{23})
(\alpha_{14}\alpha_{22})(\alpha_{18}\alpha_{19}),
$$
$$
(\alpha_{2}\alpha_{13})(\alpha_{3}\alpha_{22})
(\alpha_{4}\alpha_{9})(\alpha_{7}\alpha_{18})
(\alpha_{8}\alpha_{14})(\alpha_{11}\alpha_{20})
(\alpha_{12}\alpha_{23})(\alpha_{16}\alpha_{19})]
$$
with the orbit
$\{\alpha_{2},\alpha_{12},\alpha_{13},\alpha_{23}\}$.
%gap> H1:=Image(F[1]);
%Group([ (3,14,8,22)(4,20,9,11)(5,7,15,19)(6,18,24,16)(10,17)(12,13), (5,18)(6,7)(10,17)(11,20)(12,13)(14,22)(15,16)
%(19,24) ])
%gap> Orbits(H1);
%[ [ 3, 14, 8, 22 ], [ 4, 20, 9, 11 ], [ 5, 7, 18, 15, 6, 24, 19, 16 ], [ 10, 17 ], [ 12, 13 ] ]
The $G_1\cong D_8$ is marked by
$$
G=H_{39,2}\supset G_1=
$$
$$
[(\alpha_{3}\alpha_{14}\alpha_{8}\alpha_{22})(\alpha_{4}\alpha_{20}\alpha_{9}\alpha_{11})
(\alpha_{5}\alpha_{7}\alpha_{15}\alpha_{19})(\alpha_{6}\alpha_{18}\alpha_{24}\alpha_{16})
(\alpha_{10}\alpha_{17})(\alpha_{12}\alpha_{13}),
$$
$$
(\alpha_{5}\alpha_{18})(\alpha_{6}\alpha_{7})(\alpha_{10}\alpha_{17})(\alpha_{11}\alpha_{20})
(\alpha_{12}\alpha_{13})(\alpha_{14}\alpha_{22})(\alpha_{15}\alpha_{16})(\alpha_{19}\alpha_{24})]
$$
with suborbits
$\{\alpha_{2}\}$, $\{\alpha_{23}\}$, $\{\alpha_{12},\alpha_{13}\}$.
Both $G$ and $G_1$ are marked by $\ast$ in Tables 1---4.

\medskip

Case 15:
({\bf n}=10,
$(\aaa_1,\aaa_1,4\aaa_1)\subset 6\aaa_1)\
\Longleftarrow\ ({\bf n}=22,(2\aaa_1,4\aaa_1)\subset 6\aaa_1)$.

Similar to Cases 1 and 12. By \cite{Nik11}, the $G\cong C_2\times D_8$
is marked by $N_{23}$ and $G=H_{22,2}$ of Case 12 with orbits
$\{\alpha_{12},\alpha_{23}\}$, $\{\alpha_{6},\alpha_{15},\alpha_{21},\alpha_{17}\}$.
%gap> H4:=Image(F[4]);
%Group([ (2,9)(3,18)(4,13)(6,21)(7,22)(8,14)(10,24)(16,19), (2,13)(3,9)(4,22)(6,15)(7,18)(8,20)(11,14)(17,21) ])
%gap> Orbits(H4);
%[ [ 2, 9, 13, 3, 4, 18, 22, 7 ], [ 6, 21, 15, 17 ], [ 8, 14, 20, 11 ], [ 10, 24 ], [ 16, 19 ] ]
The $G_1\cong D_8$ is marked by
$$
G=H_{22,2}\supset G_1=
$$
$$
[(\alpha_{2}\alpha_{9})(\alpha_{3}\alpha_{18})(\alpha_{4}\alpha_{13})(\alpha_{6}\alpha_{21})
(\alpha_{7}\alpha_{22})(\alpha_{8}\alpha_{14})(\alpha_{10}\alpha_{24})(\alpha_{16}\alpha_{19}),
$$
$$
(\alpha_{2}\alpha_{13})(\alpha_{3}\alpha_{9})(\alpha_{4}\alpha_{22})(\alpha_{6}\alpha_{15})
(\alpha_{7}\alpha_{18})(\alpha_{8}\alpha_{20})(\alpha_{11}\alpha_{14})(\alpha_{17}\alpha_{21})]
$$
with suborbits
$\{\alpha_{12}\}$, $\{\alpha_{23}\}$, $\{\alpha_{6},\alpha_{15},\alpha_{21},\alpha_{17}\}$.
Both $G$ and $G_1$
are marked by $\ast$ in Tables 1---4.

\medskip

Case 16: ({\bf n}=10,
$(\aaa_1,\aaa_1,8\aaa_1)\subset 10\aaa_1)\
\Longleftarrow\ ({\bf n}=22,(2\aaa_1,8\aaa_1)\subset 10\aaa_1)$.

\medskip

By \cite{Nik11} which describes markings of these cases by Niemeier lattices,
we should consider two cases, when $({\bf n}=22,(2\aaa_1,8\aaa_1)\subset 10\aaa_1)$ is marked
by $N_{23}$ and by $N_{21}$.

If $({\bf n}=22,(2\aaa_1,8\aaa_1)\subset 10\aaa_1)$ is marked by $N_{23}$, then
$G$ is marked by $G=H_{22,2}$ from Case 12 with orbits
$\{\alpha_{12},\alpha_{23}\}$, $\{\alpha_{2},\alpha_{22},\alpha_{13},\alpha_{7},\alpha_{9},
\alpha_{3},\alpha_{18},\alpha_{4}\}$, and  $G_1\cong D_8$ is marked by
$G_1\subset G=H_{22,2}$ from Case 15 with suborbits
$\{\alpha_{12}\}$, $\{\alpha_{23}\}$, $\{\alpha_{2},\alpha_{22},\alpha_{13},\alpha_{7},\alpha_{9},
\alpha_{3},\alpha_{18},\alpha_{4}\}$.

%N_{21}
%gap> H4:=Image(HH[4]);
%G:=Group([ (7,9)(10,15)(11,14)(12,13)(16,22)(17,23)(18,24)(19,21),(4,6)(7,18)(8,17)(9,16)(10,12)(19,22)(20,23)(21,24),
%(1,3)(4,6)(7,9)(10,12)(13,15)(16,18)(19,21)(22,24) ]);
%gap> Orbits(H4);
%[ [ 1, 3 ], [ 4, 6 ], [ 7, 9, 18, 16, 24, 22, 21, 19 ], [ 8, 17, 23, 20 ], [ 10, 15, 12, 13 ], [ 11, 14 ] ]

If $({\bf n}=22,(2\aaa_1,8\aaa_1)\subset 10\aaa_1)$ is marked by $N_{21}$, then
$G$ is marked by $G=H_{22,1}$ from Case 6 with orbits
%$
%\{\alpha_{1,1},\alpha_{3,1}\},\
%\{\alpha_{1,2},\alpha_{3,2}\},\
%\{\alpha_{1,3},\alpha_{3,3},\alpha_{3,6},\alpha_{1,6},\alpha_{3,8},
%\alpha_{1,8},\linebreak
%\alpha_{3,7},\alpha_{1,7}\},\
%\{\alpha_{2,3},\alpha_{2,6},\alpha_{2,8},\alpha_{2,7}\},\
%\{\alpha_{1,4},\alpha_{3,5},\alpha_{3,4},\alpha_{1,5}\},\
%\{\alpha_{2,4},\alpha_{2,5}\}
%$.
$\{\alpha_{1,1},\alpha_{3,1}\}$,
$\{\alpha_{1,3},\alpha_{3,3},\alpha_{3,6},\alpha_{1,6},\alpha_{3,8},
\alpha_{1,8},\alpha_{3,7},\alpha_{1,7}\}$,
%gap> H1:=Image(F[1]);
%Group([ (7,9)(10,15)(11,14)(12,13)(16,22)(17,23)(18,24)(19,21), (4,6)(7,18)(8,17)(9,16)(10,12)(19,22)(20,23)(21,24) ])
%gap> Orbits(H1);
%[ [ 4, 6 ], [ 7, 9, 18, 16, 24, 22, 21, 19 ], [ 8, 17, 23, 20 ], [ 10, 15, 12, 13 ], [ 11, 14 ] ]
and  $G_1\cong D_8$ is marked by
$$
G=H_{22,1}\supset G_1=
$$
$$
[(\alpha_{1,3}\alpha_{3,3})(\alpha_{1,4}\alpha_{3,5})(\alpha_{2,4}\alpha_{2,5})(\alpha_{3,4}\alpha_{1,5})
(\alpha_{1,6}\alpha_{1,8})(\alpha_{2,6}\alpha_{2,8})(\alpha_{3,6}\alpha_{3,8})(\alpha_{1,7}\alpha_{3,7}),
$$
$$
(\alpha_{1,2}\alpha_{3,2})(\alpha_{1,3}\alpha_{3,6})(\alpha_{2,3}\alpha_{2,6})(\alpha_{3,3}\alpha_{1,6})
(\alpha_{1,4}\alpha_{3,4})(\alpha_{1,7}\alpha_{1,8})(\alpha_{2,7}\alpha_{2,8})(\alpha_{3,7}\alpha_{3,8})]
$$
with suborbits
$\{\alpha_{1,1}\}$, $\{\alpha_{3,1}\}$,
$\{\alpha_{1,3},\alpha_{3,3},\alpha_{3,6},\alpha_{1,6},\alpha_{3,8},
\alpha_{1,8},\alpha_{3,7},\alpha_{1,7}\}$.

Case 17:

({\bf n}=10,
$
\left(\begin{array}{ccc}
(2\aaa_1)_{I} & (6\aaa_1)_I & (6\aaa_1)_I  \\
         & 4\aaa_1   & (8\aaa_1)_{II}   \\
         &         & 4\aaa_1
\end{array}\right)
\subset 10\aaa_1)\
\Longleftarrow\ ({\bf n}=22,(2\aaa_1,8\aaa_1)\subset 10\aaa_1)$.

Similar to Case 16.

If $({\bf n}=22,(2\aaa_1,8\aaa_1)\subset 10\aaa_1)$ is marked by $N_{23}$, then $G=H_{22,2}$ from Case 12 with orbits
$\{\alpha_{12},\alpha_{23}\}$,
$\{\alpha_{2},\alpha_{22},\alpha_{13},\alpha_{7},\alpha_{9},\alpha_{3},\alpha_{18},\alpha_{4}\}$,
and  $G_1\cong D_8$ is marked by
%N23
%gap> H2;
%Group([ (3,4)(6,15)(8,11)(9,22)(12,23)(14,20)(16,19)(17,21), (2,9)(3,18)(4,13)(6,21)(7,22)(8,14)(10,24)(16,19) ])
%gap> Orbits(H2);
%[ [ 2, 9, 22, 7 ], [ 3, 4, 18, 13 ], [ 6, 15, 21, 17 ], [ 8, 11, 14, 20 ], [ 10, 24 ], [ 12, 23 ], [ 16, 19 ] ]
$$
G=H_{22,2}\supset G_1=
$$
$$
[(\alpha_{3}\alpha_{4})(\alpha_{6}\alpha_{15})(\alpha_{8}\alpha_{11})(\alpha_{9}\alpha_{22})
(\alpha_{12}\alpha_{23})(\alpha_{14}\alpha_{20})(\alpha_{16}\alpha_{19})(\alpha_{17}\alpha_{21}),
$$
$$
(\alpha_{2}\alpha_{9})(\alpha_{3}\alpha_{18})(\alpha_{4}\alpha_{13})(\alpha_{6}\alpha_{21})
(\alpha_{7}\alpha_{22})(\alpha_{8}\alpha_{14})(\alpha_{10}\alpha_{24})(\alpha_{16}\alpha_{19})]
$$
with suborbits
$\{\alpha_{12},\alpha_{23}\}$, $\{\alpha_{2},\alpha_{9},\alpha_{22},\alpha_{7}\}$,
$\{\alpha_{3},\alpha_{4},\alpha_{18},\alpha_{13}\}$.

If $({\bf n}=22,(2\aaa_1,8\aaa_1)\subset 10\aaa_1)$ is marked by $N_{21}$, then
$G=H_{22,1}$ from Case 6 with orbits
%$
%\{\alpha_{1,1},\alpha_{3,1}\},\
%\{\alpha_{1,2},\alpha_{3,2}\},\
%\{\alpha_{1,3},\alpha_{3,3},\alpha_{3,6},\alpha_{1,6},\alpha_{3,8},
%\alpha_{1,8},\linebreak
%\alpha_{3,7},\alpha_{1,7}\},\
%\{\alpha_{2,3},\alpha_{2,6},\alpha_{2,8},\alpha_{2,7}\},\
%\{\alpha_{1,4},\alpha_{3,5},\alpha_{3,4},\alpha_{1,5}\},\
%\{\alpha_{2,4},\alpha_{2,5}\}
%$.
$\{\alpha_{1,1},\alpha_{3,1}\}$,
$\{\alpha_{1,3},\alpha_{3,3},\alpha_{3,6},\alpha_{1,6},\alpha_{3,8},
\alpha_{1,8},\alpha_{3,7},\alpha_{1,7}\}$,
%gap> H3:=Image(F[3]);
%Group([ (4,6)(7,18)(8,17)(9,16)(10,12)(19,22)(20,23)(21,24), (1,3)(4,6)(10,13)(11,14)(12,15)(16,24)(17,23)(18,22) ])
%gap> Orbits(H3);
%[ [ 1, 3 ], [ 4, 6 ], [ 7, 18, 22, 19 ], [ 8, 17, 23, 20 ], [ 9, 16, 24, 21 ], [ 10, 12, 13, 15 ], [ 11, 14 ] ]
and  $G_1\cong D_8$ is marked by
$$
G=H_{22,1}\supset G_1=
$$
$$
[(\alpha_{1,2}\alpha_{3,2})(\alpha_{1,3}\alpha_{3,6})(\alpha_{2,3}\alpha_{2,6})(\alpha_{3,3}\alpha_{1,6})
(\alpha_{1,4}\alpha_{3,4})(\alpha_{1,7}\alpha_{1,8})(\alpha_{2,7}\alpha_{2,8})(\alpha_{3,7}\alpha_{3,8}),
$$
$$
(\alpha_{1,1}\alpha_{3,1})(\alpha_{1,2}\alpha_{3,2})(\alpha_{1,4}\alpha_{1,5})(\alpha_{2,4}\alpha_{2,5})
(\alpha_{3,4}\alpha_{3,5})(\alpha_{1,6}\alpha_{3,8})(\alpha_{2,6}\alpha_{2,8})(\alpha_{3,6}\alpha_{1,8})]
$$
with suborbits
$\{\alpha_{1,1},\alpha_{3,1}\}$, $\{\alpha_{1,3},\alpha_{3,6},\alpha_{1,8},\alpha_{1,7}\}$,
$\{\alpha_{3,3},\alpha_{1,6},\alpha_{3,8},\alpha_{3,7}\}$.

\medskip

Case 18:

({\bf n}=10,
$
\left(\begin{array}{ccc}
(2\aaa_1)_{II}  & 6\aaa_1 & 6\aaa_1   \\
         & 4\aaa_1 & (8\aaa_1)_{II}   \\
         &         & 4\aaa_1
\end{array}\right)
\subset 10\aaa_1)\
\Longleftarrow\ ({\bf n}=22,(2\aaa_1,8\aaa_1)\subset 10\aaa_1)$.

Similar to Case 17.

If $({\bf n}=22,(2\aaa_1,8\aaa_1)\subset 10\aaa_1)$ is marked by $N_{23}$, then
$G=H_{22,2}$ from Case 12 with orbits
$\{\alpha_{12},\alpha_{23}\}$,
$\{\alpha_{2},\alpha_{22},\alpha_{13},\alpha_{7},\alpha_{9},\alpha_{3},\alpha_{18},\alpha_{4}\}$,
and $D_8\cong G_1$ is marked by
%N23
%gap> H1;
%Group([ (3,4)(6,15)(8,11)(9,22)(12,23)(14,20)(16,19)(17,21), (2,3)(4,7)(8,14)(9,18)(10,24)(12,23)(13,22)(15,17) ])
%gap> Orbits(H1);
%[ [ 2, 3, 4, 7 ], [ 6, 15, 17, 21 ], [ 8, 11, 14, 20 ], [ 9, 22, 18, 13 ], [ 10, 24 ], [ 12, 23 ], [ 16, 19 ] ]
$$
G=H_{22,2}\supset G_1=
$$
$$
[(\alpha_{3}\alpha_{4})(\alpha_{6}\alpha_{15})(\alpha_{8}\alpha_{11})(\alpha_{9}\alpha_{22})
(\alpha_{12}\alpha_{23})(\alpha_{14}\alpha_{20})(\alpha_{16}\alpha_{19})(\alpha_{17}\alpha_{21}),
$$
$$
(\alpha_{2}\alpha_{3})(\alpha_{4}\alpha_{7})(\alpha_{8}\alpha_{14})(\alpha_{9}\alpha_{18})
(\alpha_{10}\alpha_{24})(\alpha_{12}\alpha_{23})(\alpha_{13}\alpha_{22})(\alpha_{15}\alpha_{17})]
$$
with suborbits
$\{\alpha_{12},\alpha_{23}\}$, $\{\alpha_{2},\alpha_{3},\alpha_{4},\alpha_{7}\}$,
$\{\alpha_{9},\alpha_{22},\alpha_{18},\alpha_{13}\}$.

If $({\bf n}=22,(2\aaa_1,8\aaa_1)\subset 10\aaa_1)$ is marked by $N_{21}$, then
$G=H_{22,1}$ from Case 6 with orbits
%$
%\{\alpha_{1,1},\alpha_{3,1}\},\
%\{\alpha_{1,2},\alpha_{3,2}\},\
%\{\alpha_{1,3},\alpha_{3,3},\alpha_{3,6},\alpha_{1,6},\alpha_{3,8},
%\alpha_{1,8},\linebreak
%\alpha_{3,7},\alpha_{1,7}\},\
%\{\alpha_{2,3},\alpha_{2,6},\alpha_{2,8},\alpha_{2,7}\},\
%\{\alpha_{1,4},\alpha_{3,5},\alpha_{3,4},\alpha_{1,5}\},\
%\{\alpha_{2,4},\alpha_{2,5}\}
%$.
$\{\alpha_{1,1},\alpha_{3,1}\}$,
$\{\alpha_{1,3},\alpha_{3,3},\alpha_{3,6},\alpha_{1,6},\alpha_{3,8},
\alpha_{1,8},\alpha_{3,7},\alpha_{1,7}\}$,
%gap> H4:=Image(F[4]);
%Group([ (1,3)(7,16)(8,17)(9,18)(13,15)(19,24)(20,23)(21,22), (1,3)(4,6)(10,13)(11,14)(12,15)(16,24)(17,23)(18,22) ])
%gap> Orbits(H4);
%[ [ 1, 3 ], [ 4, 6 ], [ 7, 16, 24, 19 ], [ 8, 17, 23, 20 ], [ 9, 18, 22, 21 ], [ 10, 13, 15, 12 ], [ 11, 14 ] ]
and $G_1\cong D_8$ is marked by
$$
G=H_{22,1}\supset G_1=
$$
$$
[(\alpha_{1,1}\alpha_{3,1})(\alpha_{1,3}\alpha_{1,6})(\alpha_{2,3}\alpha_{2,6})(\alpha_{3,3}\alpha_{3,6})
(\alpha_{1,5}\alpha_{3,5})(\alpha_{1,7}\alpha_{3,8})(\alpha_{2,7}\alpha_{2,8})(\alpha_{3,7}\alpha_{1,8}),
$$
$$
(\alpha_{1,1}\alpha_{3,1})(\alpha_{1,2}\alpha_{3,2})(\alpha_{1,4}\alpha_{1,5})(\alpha_{2,4}\alpha_{2,5})
(\alpha_{3,4}\alpha_{3,5})(\alpha_{1,6}\alpha_{3,8})(\alpha_{2,6}\alpha_{2,8})(\alpha_{3,6}\alpha_{1,8})]
$$
with suborbits
$\{\alpha_{1,1},\alpha_{3,1}\}$, $\{\alpha_{1,3},\alpha_{1,6},\alpha_{3,8},\alpha_{1,7}\}$,
$\{\alpha_{3,3},\alpha_{3,6},\alpha_{1,8},\alpha_{3,7}\}$.

\medskip

Case 19:
({\bf n}=10,
$
\left(\begin{array}{ccc}
4\aaa_1  & (8\aaa_1)_I & (8\aaa_1)_{II}  \\
         & 4\aaa_1   & (8\aaa_1)_{I} \\
         &         & 4\aaa_1
\end{array}\right)
\subset 12\aaa_1)\
\Longleftarrow\ ({\bf n}=22,(4\aaa_1,8\aaa_1)\subset 12\aaa_1)$.

Similar to Cases 1, 12 and 18. By \cite{Nik11}, the $G\cong C_2\times D_8$
is marked by $N_{23}$ and $G=H_{22,2}$ from Case 12 with orbits
%$$
%G=H_{22,2}=
%$$
%G2:=Group([(3,4)(6,15)(8,11)(9,22)(12,23)(14,20)(16,19)(17,21),
%(2,22)(3,13)(4,18)(7,9)(10,24)(11,20)(15,17)(16,19),
%(2,13)(3,22)(4,9)(7,18)(8,14)(11,20)(12,23)(16,19)]);
%with orbits
%[ [ 2, 22, 13, 9, 3, 7, 4, 18 ], [ 6, 15, 17, 21 ], [ 8, 11, 14, 20 ], [ 10, 24 ], [ 12, 23 ],  [ 16, 19 ] ]
$\{\alpha_{6},\alpha_{15},\alpha_{17},\alpha_{21}\}$,
$\{\alpha_{2},\alpha_{22},\alpha_{13},\alpha_{9},\alpha_{3},\alpha_{7},\alpha_{4},\alpha_{8}\}$.
The $G_1\cong D_8$ is marked by
%gap> H1:=Image(F[1]);
%Group([ (3,4)(6,15)(8,11)(9,22)(12,23)(14,20)(16,19)(17,21), (2,3)(4,7)(8,14)(9,18)(10,24)(12,23)
%(13,22)(15,17) ])
%gap> Orbits(H1);
%[ [ 2, 3, 4, 7 ], [ 6, 15, 17, 21 ], [ 8, 11, 14, 20 ], [ 9, 22, 18, 13 ], [ 10, 24 ], [ 12, 23 ],
%  [ 16, 19 ] ]
$G_1\subset G=H_{22,2}$ of Case 18 with suborbits
$\{\alpha_{6}, \alpha_{15},\alpha_{17},\alpha_{21}\}$
$\{\alpha_{2},\alpha_{3},\alpha_{4},\alpha_{7}\}$,
$\{\alpha_{9},\alpha_{22},\alpha_{18},\alpha_{13}\}$.
Both $G$ and $G_1$ are marked by $\ast$ in Tables 1---4.

\medskip

Case 20:
({\bf n}=10,
$
\left(\begin{array}{ccc}
4\aaa_1  & 4\aaa_2 & (8\aaa_1)_{II}   \\
         & 4\aaa_1 & 4\aaa_2   \\
         &         & 4\aaa_1
\end{array}\right)
\subset 4\aaa_3)
\Longleftarrow\ ({\bf n}=22,(4\aaa_1,8\aaa_1)\subset 4\aaa_3)$.

Similar to Cases 1 and 6. By \cite{Nik11}, the
$G\cong C_2\times D_8$
is marked by $N_{21}$ and
$G=H_{22,1}$ from Case 6 with orbits
$\{\alpha_{2,3},\alpha_{2,6},\alpha_{2,8},\alpha_{2,7}\}$,
$\{\alpha_{1,3},\alpha_{3,3},\alpha_{3,6},\alpha_{1,6},\alpha_{3,8},
\alpha_{1,8},\alpha_{3,7},\alpha_{1,7}\}$.
The $G_1\cong D_8$ is marked by
$$
G=H_{22,1}\supset G_1=
$$
%gap> H3:=Image(F[3]);
%Group([ (4,6)(7,18)(8,17)(9,16) (10,12)(19,22)(20,23)(21,24), (1,3)(7,18,19,22)(8,17,20,23)(9,16,21,24)(10,13,12,15)
%(11,14) ])
$$
[(\alpha_{1,2}\alpha_{3,2})(\alpha_{1,3}\alpha_{3,6})(\alpha_{2,3}\alpha_{2,6})(\alpha_{3,3}\alpha_{1,6})
(\alpha_{1,4}\alpha_{3,4})(\alpha_{1,7}\alpha_{1,8})(\alpha_{2,7}\alpha_{2,8})(\alpha_{3,7}\alpha_{3,8}),
$$
$$
(\alpha_{1,1}\alpha_{3,1}) (\alpha_{1,3}\alpha_{3,6}\alpha_{1,7}\alpha_{1,8})
(\alpha_{2,3}\alpha_{2,6}\alpha_{2,7}\alpha_{2,8})
(\alpha_{3,3}\alpha_{1,6}\alpha_{3,7}\alpha_{3,8}) (\alpha_{1,4}\alpha_{1,5}\alpha_{3,4}\alpha_{3,5})
(\alpha_{2,4}\alpha_{2,5})]
$$
with suborbits
%[ 8, 17, 20, 23 ], [ 7, 18, 19, 22 ], [ 9, 16, 21, 24 ],
$\{ \alpha_{2,3},\alpha_{2,6},\alpha_{2,7},\alpha_{2,8} \}$,
$\{ \alpha_{1,3}, \alpha_{3,6},\alpha_{1,7}, \alpha_{1,8}\}$,
$\{ \alpha_{3,3}, \alpha_{1,6},\alpha_{3,7}, \alpha_{3,8} \}$.
%gap> Orbits(H3);
%[ [ 1, 3 ], [ 4, 6 ], [ 7, 18, 19, 22 ], [ 8, 17, 20, 23 ], [ 9, 16, 21, 24 ], [ 10, 12, 13, 15 ], [ 11, 14 ] ]
Both $G$ and $G_1$ are marked by $\ast$ in Tables 1---4.

\medskip

Case 21:
({\bf n}=10,\
$
\left(\begin{array}{ccc}
4\aaa_1 & (8\aaa_1)_{II} & 12\aaa_1   \\
         & 4\aaa_1 & 12\aaa_1   \\
         &         & 8\aaa_1
\end{array}\right)
\subset 16\aaa_1) \Longleftarrow ({\bf n}=39,\ 16\aaa_1)$.

Similar to Cases 1 and 7. By \cite{Nik10}, the $G\cong 2^4C_2$
is marked by $N_{23}$ and
$G=H_{39,1}$ from Case 7 with the orbit
$\{\alpha_{2},\alpha_{12},\alpha_{13},\alpha_{3},\alpha_{18},\alpha_{8},
\alpha_{23},\alpha_{19},\alpha_{7},\alpha_{22},\alpha_{4},
\alpha_{14},\alpha_{20},\alpha_{11},\alpha_{16},\alpha_{9}\}$.
The $G_1\cong D_8$ is marked by
$$
G=H_{39,1}\supset G_1=
$$
%gap> H1:=Image(F[1]);
%Group([ (2,7,18,13)(3,9,4,22)(5,10)(8,16,20,23)(11,19,14,12)(17,21), (5,10)(7,13)(8,12)(9,22)(11,23)(14,16)(17,21)
%(19,20) ])
$$
[(\alpha_{2}\alpha_{7}\alpha_{18}\alpha_{13})(\alpha_{3}\alpha_{9}\alpha_{4}\alpha_{22})
(\alpha_{5}\alpha_{10})
(\alpha_{8}\alpha_{16}\alpha_{20}\alpha_{23})(\alpha_{11}\alpha_{19}\alpha_{14}\alpha_{12})
(\alpha_{17}\alpha_{21}),
$$
$$
(\alpha_{5}\alpha_{10})(\alpha_{7}\alpha_{13})(\alpha_{8}\alpha_{12})(\alpha_{9}\alpha_{22})
(\alpha_{11}\alpha_{23})(\alpha_{14}\alpha_{16})(\alpha_{17}\alpha_{21})(\alpha_{19}\alpha_{20})]
$$
with suborbits
%gap> Orbits(H1);
%[ [ 2, 7, 18, 13 ], [ 3, 9, 4, 22 ], [ 5, 10 ], [ 8, 16, 12, 20, 14, 11, 23, 19 ], [ 17, 21 ] ]
$\{\alpha_{2},\alpha_{7},\alpha_{18},\alpha_{13}\}$, $\{\alpha_{3},\alpha_{9},\alpha_{4},\alpha_{22}\}$,
$\{\alpha_{8},\alpha_{16},\alpha_{12},\alpha_{20},\alpha_{14},\alpha_{11},\alpha_{23},\alpha_{19}\}$.
Both $G$ and $G_1$ are marked by $\ast$ in Tables 1---4.

\medskip

Case 22:
$({\bf n}=10,\
(4\aaa_1,4\aaa_1,2\aaa_2)\subset 6\aaa_2)
\Longleftarrow ({\bf n}=34,\ 6\aaa_2)$.

Similar to Case 1. By \cite{Nik10}, the $G\cong \SSS_4$
is marked by $N_{22}$ and
$$
G=H_{34,1}=
$$
%%gap> h34c1;
%%G1:=Group([ (3,10)(4,9)(7,16)(8,15)(13,18)(14,17)(21,23)(22,24),
%%(1,5,8)(2,6,7)(3,9,22)(4,10,21)(13,17,24)(14,18,23) ]);
$$
[(\alpha_{1,2}\alpha_{2,5})(\alpha_{2,2}\alpha_{1,5})
(\alpha_{1,4}\alpha_{2,8})(\alpha_{2,4}\alpha_{1,8})
(\alpha_{1,7}\alpha_{2,9})(\alpha_{2,7}\alpha_{1,9})
(\alpha_{1,11}\alpha_{1,12})(\alpha_{2,11}\alpha_{2,12}),
$$
$$
(\alpha_{1,1}\alpha_{1,3}\alpha_{2,4})
(\alpha_{2,1}\alpha_{2,3}\alpha_{1,4})
(\alpha_{1,2}\alpha_{1,5}\alpha_{2,11})
(\alpha_{2,2}\alpha_{2,5}\alpha_{1,11})
(\alpha_{1,7}\alpha_{1,9}\alpha_{2,12})
(\alpha_{2,7}\alpha_{2,9}\alpha_{1,12})]
$$
with the orbit
%%Orbits(h34c1);
%%[ [ 1, 5, 8, 15 ], [ 2, 6, 7, 16 ],
%%[ 3, 10, 9, 21, 4, 22, 23, 24, 14, 13, 17, 18 ] ]
$\{\alpha_{1,2},\alpha_{2,5},
\alpha_{1,5},\alpha_{1,11},\alpha_{2,2},\alpha_{2,11},
\alpha_{1,12},\alpha_{2,12},\alpha_{2,7},\alpha_{1,7},\alpha_{1,9},\alpha_{2,9}\}$.
The $G_1\cong D_8$ is marked by
$$
G=H_{34,1}\supset G_1=
$$
%gap> H1:=Image(F[1]);
%Group([ (3,21)(4,22)(5,15)(6,16)(9,17)(10,18)(13,23)(14,24), (1,5)(2,6)(3,14)(4,13)(7,16)(8,15)(9,18)(10,17) ])
%gap>
%gap> Orbits(H1);
%[ [ 1, 5, 15, 8 ], [ 2, 6, 16, 7 ], [ 3, 21, 14, 24 ], [ 4, 22, 13, 23 ], [ 9, 17, 18, 10 ] ]
$$
[(\alpha_{1,2}\alpha_{1,11})(\alpha_{2,2}\alpha_{2,11})(\alpha_{1,3}\alpha_{1,8})(\alpha_{2,3}\alpha_{2,8})
(\alpha_{1,5}\alpha_{1,9})(\alpha_{2,5}\alpha_{2,9})(\alpha_{1,7}\alpha_{1,12})(\alpha_{2,7}\alpha_{2,12}),
$$
$$
(\alpha_{1,1}\alpha_{1,3})(\alpha_{2,1}\alpha_{2,3})(\alpha_{1,2}\alpha_{2,7})(\alpha_{2,2}\alpha_{1,7})
(\alpha_{1,4}\alpha_{2,8})(\alpha_{2,4}\alpha_{1,8})(\alpha_{1,5}\alpha_{2,9})(\alpha_{2,5}\alpha_{1,9})]
$$
with suborbits
$\{\alpha_{1,2},\alpha_{1,11},\alpha_{2,7},\alpha_{2,12}\}$, $\{\alpha_{2,2},\alpha_{2,11},\alpha_{1,7},\alpha_{1,12}\}$,
$\{\alpha_{1,5},\alpha_{1,9},\alpha_{2,9},\alpha_{2,5}\}$.
Both $G$ and $G_1$ are marked by $\ast$ in Tables 1---4.

\medskip

Case 23:
({\bf n}=10,\
$
\left(\begin{array}{cccc}
\aaa_1 & 2\aaa_1 & 3\aaa_1      &  5\aaa_1       \\
       & \aaa_1  & 3\aaa_1      &  5\aaa_1       \\
       &         & (2\aaa_1)_I  & (6\aaa_1)_I   \\
       &         &              &  4\aaa_1
\end{array}\right)
\subset 8\aaa_1) \Longleftarrow ({\bf n}=56,\ 8\aaa_1)$.

Similar to Case 1. By \cite{Nik10},
the $G\cong \Gamma_{15}a_1$ is marked by $N_{23}$ and
$$
G=H_{56,1}=
$$
%G1:=Group([(5,23)(6,10)(7,9)(8,16)(13,18)(14,24)(15,20)(19,21),
%(1,14)(2,23)(3,5)(6,18)(7,8)(9,13)(10,16)(17,24),
%(1,14)(2,24)(3,5)(6,10)(12,22)(16,18)(17,23)(19,20)]);
$$
[(\alpha_{5}\alpha_{23})(\alpha_{6}\alpha_{10})(\alpha_{7}\alpha_{9})
(\alpha_{8}\alpha_{16})(\alpha_{13}\alpha_{18})(\alpha_{14}\alpha_{24})
(\alpha_{15}\alpha_{20})(\alpha_{19}\alpha_{21}),
$$
$$
(\alpha_{1}\alpha_{14})(\alpha_{2}\alpha_{23})(\alpha_{3}\alpha_{5})
(\alpha_{6}\alpha_{18})(\alpha_{7},\alpha_{8})(\alpha_{9}\alpha_{13})
(\alpha_{10}\alpha_{16})(\alpha_{17}\alpha_{24}),
$$
$$
(\alpha_{1}\alpha_{14})(\alpha_{2}\alpha_{24})(\alpha_{3}\alpha_{5})
(\alpha_{6}\alpha_{10})(\alpha_{12}\alpha_{22})(\alpha_{16}\alpha_{18})
(\alpha_{17}\alpha_{23})(\alpha_{19}\alpha_{20})]
$$
with the orbit
%$
% \{\alpha_{1}, \alpha_{14}, \alpha_{24}, \alpha_{17}, \alpha_{2},
%\alpha_{23}, \alpha_{5}, \alpha_{3} \},\  \{ \alpha_{6}, \alpha_{10}, \alpha_{18},
%\alpha_{16}, \alpha_{13}, \alpha_{8},\linebreak
%\alpha_{9}, \alpha_{7}\},\
%\{ \alpha_{12}, \alpha_{22}\},\
%\{ \alpha_{15}, \alpha_{20}, \alpha_{19}, \alpha_{21}\}$;
$\{\alpha_{1}, \alpha_{14}, \alpha_{24}, \alpha_{17}, \alpha_{2},\alpha_{23}, \alpha_{5}, \alpha_{3} \}$.
%gap> H3:=Image(F[3]);
%Group([ (5,23)(6,10)(7,9)(8,16)(13,18)(14,24)(15,20)(19,21), (2,17)(6,16)(7,8)(9,13)(10,18)(12,22)(19,20)(23,24) ])
%gap> Orbits(H3);
%[ [ 2, 17 ], [ 5, 23, 24, 14 ], [ 6, 10, 16, 18, 8, 13, 7, 9 ], [ 12, 22 ], [ 15, 20, 19, 21 ] ]
The $G_1\cong D_8$ is marked by
$$
G=H_{56,1}\supset G_1=
$$
$$
[(\alpha_{5}\alpha_{23})(\alpha_{6}\alpha_{10})(\alpha_{7}\alpha_{9})(\alpha_{8}\alpha_{16})
(\alpha_{13}\alpha_{18})(\alpha_{14}\alpha_{24})(\alpha_{15}\alpha_{20})(\alpha_{19}\alpha_{21}),
$$
$$
(\alpha_{2}\alpha_{17})(\alpha_{6}\alpha_{16})(\alpha_{7}\alpha_{8})(\alpha_{9}\alpha_{13})
(\alpha_{10}\alpha_{18})(\alpha_{12}\alpha_{22})(\alpha_{19}\alpha_{20})(\alpha_{23}\alpha_{24})]
$$
with suborbits
$\{\alpha_{1}\}$, $\{\alpha_{3}\}$, $\{\alpha_{2},\alpha_{17}\}$,
$\{\alpha_{5},\alpha_{23},\alpha_{24},\alpha_{14}\}$.
Both $G$ and $G_1$ are marked by $\ast$ in Tables 1---4.

\medskip

Case 24:
$({\bf n}=10,\ (\aaa_1,\aaa_1,(2\aaa_1)_{I},8\aaa_1)
\subset 12\aaa_1)\Longleftarrow ({\bf n}=65,\ 12\aaa_1)$.

Similar to Case 1. By \cite{Nik10},
the $G\cong 2^4D_6$ is marked by $N_{23}$ and
$$
G=H_{65,3}=
$$
%G3:=Group([(1,7,17)(3,22,5)(4,23,20)(6,14,8)(9,13,11)(10,19,21),
%(1,11,23,8)(4,6,13,20)(5,16)(7,17,9,14)(10,21,19,15)(12,24)]);
$$
[(\alpha_{1}\alpha_{7}\alpha_{17})(\alpha_{3}\alpha_{22}\alpha_{5})
(\alpha_{4}\alpha_{23}\alpha_{20})(\alpha_{6}\alpha_{14}\alpha_{8})
(\alpha_{9}\alpha_{13}\alpha_{11})(\alpha_{10}\alpha_{19}\alpha_{21}),
$$
$$
(\alpha_{1}\alpha_{11}\alpha_{23}\alpha_{8})
(\alpha_{4}\alpha_{6}\alpha_{13}\alpha_{20})
(\alpha_{5}\alpha_{16})
(\alpha_{7}\alpha_{17}\alpha_{9}\alpha_{14})
(\alpha_{10}\alpha_{21}\alpha_{19}\alpha_{15})
(\alpha_{12}\alpha_{24})];
$$
with the orbit
%$
%\{\alpha_{1},\alpha_{7},\alpha_{23},\alpha_{8},\alpha_{11},\alpha_{17},\alpha_{20},\alpha_{9},
%\alpha_{6},
%\alpha_{4},\alpha_{14},\alpha_{13}\},\
%\{\alpha_{3},\alpha_{22},\linebreak
%\alpha_{16},\alpha_{5}\},\
%\{\alpha_{10},\alpha_{19},
%\alpha_{15},\alpha_{21}\},\
%\{\alpha_{12},\alpha_{24}\}
$\{\alpha_{1},\alpha_{7},\alpha_{23},\alpha_{8},\alpha_{11},\alpha_{17},\alpha_{20},\alpha_{9},
\alpha_{6},\alpha_{4},\alpha_{14},\alpha_{13}\}$.
%gap> H1:=Image(F[1]);
%Group([ (4,7)(5,16)(6,14)(8,11)(9,13)(12,24)(15,21)(17,20), (3,5)(4,13)(6,7)(9,20)(10,15)(14,17)(16,22)(19,21) ])
%gap> Orbits(H1);
%[ [ 3, 5, 16, 22 ], [ 4, 7, 13, 6, 9, 14, 20, 17 ], [ 8, 11 ], [ 10, 15, 21, 19 ], [ 12, 24 ] ]
The $G_1\cong D_8$ is marked by
$$
G=H_{65,3}\supset G_1=
$$
$$
[(\alpha_{4}\alpha_{7})(\alpha_{5}\alpha_{16})(\alpha_{6}\alpha_{14})(\alpha_{8}\alpha_{11})
(\alpha_{9}\alpha_{13})(\alpha_{12}\alpha_{24})(\alpha_{15}\alpha_{21})(\alpha_{17}\alpha_{20}),
$$
$$
(\alpha_{3}\alpha_{5})(\alpha_{4}\alpha_{13})(\alpha_{6}\alpha_{7})(\alpha_{9}\alpha_{20})
(\alpha_{10}\alpha_{15})(\alpha_{14}\alpha_{17})(\alpha_{16}\alpha_{22})(\alpha_{19}\alpha_{21})]
$$
with suborbits
$\{\alpha_{1}\}$, $\{\alpha_{23}\}$, $\{\alpha_{8},\alpha_{11}\}$,
$\{\alpha_{4},\alpha_{7},\alpha_{13},\alpha_{6},\alpha_{9},\alpha_{14},\alpha_{20},\alpha_{17}\}$.
Both $G$ and $G_1$ are marked by $\ast$ in Tables 1---4.

\medskip

Case 25:

({\bf n}=10,\
$
\left(\begin{array}{cccc}
\aaa_1 & 2\aaa_1 & 5\aaa_1      &  5\aaa_1       \\
       & \aaa_1  & 5\aaa_1      &  5\aaa_1       \\
       &         & 4\aaa_1      & (8\aaa_1)_I    \\
       &         &              &  4\aaa_1
\end{array}\right)
\subset 10\aaa_1)\
\Longleftarrow ({\bf n}=22,\ (2\aaa_1, (4\aaa_1,4\aaa_1)_{II})\subset 10\aaa_1)$.

Similar to Case 1. By \cite{Nik11},
the $G\cong C_2\times D_8$ is marked by $N_{23}$ and
%%H3:=Group([(1,10)(3,22)(4,9)(5,24)(6,15)(8,14)(11,20)(17,21),
%(2,22)(3,13)(4,18)(7,9)(10,24)(11,20)(15,17)(16,19),
%(2,13)(3,22)(4,9)(7,18)(8,14)(11,20)(12,23)(16,19)]);
$$
G=H_{22,3}=
$$
$$
[(\alpha_{1}\alpha_{10})(\alpha_{3}\alpha_{22})
(\alpha_{4}\alpha_{9})(\alpha_{5}\alpha_{24})
(\alpha_{6}\alpha_{15})(\alpha_{8}\alpha_{14})
(\alpha_{11}\alpha_{20})(\alpha_{17}\alpha_{21}),
$$
$$
(\alpha_{2}\alpha_{22})(\alpha_{3}\alpha_{13})
(\alpha_{4}\alpha_{18})(\alpha_{7}\alpha_{9})
(\alpha_{10}\alpha_{24})(\alpha_{11}\alpha_{20})
(\alpha_{15}\alpha_{17})(\alpha_{16}\alpha_{19}),
$$
$$
(\alpha_{2}\alpha_{13})(\alpha_{3}\alpha_{22})
(\alpha_{4}\alpha_{9})(\alpha_{7}\alpha_{18})
(\alpha_{8}\alpha_{14})(\alpha_{11}\alpha_{20})
(\alpha_{12}\alpha_{23})(\alpha_{16}\alpha_{19})]
$$
with orbits
%$
%\{\alpha_{1},\alpha_{10},\alpha_{5},\alpha_{24}\},\
%\{\alpha_{2},\alpha_{22},\alpha_{13},\alpha_{3}\},\
%\{\alpha_{4},\alpha_{9},\alpha_{18},\alpha_{7}\},\
%\linebreak
%\{\alpha_{6},\alpha_{15},\alpha_{21},\alpha_{17}\},\
%\{\alpha_{8},\alpha_{14}\},\
%\{\alpha_{11},\alpha_{20}\},\
%\{\alpha_{12},\alpha_{23}\},\ \{\alpha_{16},\alpha_{19}\}$.
$\{\alpha_{8},\alpha_{14}\}$,
$\{\alpha_{1},\alpha_{10},\alpha_{5},\alpha_{24}\}$,
$\{\alpha_{2},\alpha_{22},\alpha_{13},\alpha_{3}\}$.
%gap> H4:=Image(F[4]);
%Group([ (2,22)(3,13)(4,18)(7,9)(10,24)(11,20)(15,17)(16,19), (1,10)(2,13)(5,24)(6,15)(7,18)(12,23)(16,19)(17,21) ])
%gap> Orbits(H4);
%[ [ 1, 10, 24, 5 ], [ 2, 22, 13, 3 ], [ 4, 18, 7, 9 ], [ 6, 15, 17, 21 ], [ 11, 20 ], [ 12, 23 ], [ 16, 19 ] ]
The $G_1\cong D_8$ is marked by
$$
G=H_{22,3}\supset G_1=
$$
$$
[(\alpha_{2}\alpha_{22})(\alpha_{3}\alpha_{13})(\alpha_{4}\alpha_{18})(\alpha_{7}\alpha_{9})
(\alpha_{10}\alpha_{24})(\alpha_{11}\alpha_{20})(\alpha_{15}\alpha_{17})(\alpha_{16}\alpha_{19}),
$$
$$
(\alpha_{1}\alpha_{10})(\alpha_{2}\alpha_{13})(\alpha_{5}\alpha_{24})(\alpha_{6}\alpha_{15})
(\alpha_{7}\alpha_{18})(\alpha_{12}\alpha_{23})(\alpha_{16}\alpha_{19})(\alpha_{17}\alpha_{21})]
$$
with suborbits
$\{\alpha_{8}\}$, $\{\alpha_{14}\}$,
$\{\alpha_{1},\alpha_{10},\alpha_{24},\alpha_{5}\}$,
$\{\alpha_{2},\alpha_{22},\alpha_{13},\alpha_{3}\}$.
Both $G$ and $G_1$ are marked by $\ast$ in Tables 1---4.

\medskip

Case 26:
$({\bf n}=10,\  (\aaa_1,\,\aaa_1,\,4\aaa_1,\,8\aaa_1)\subset 14\aaa_1)\
\Longleftarrow\ ({\bf n}=22,\ (2\aaa_1,4\aaa_1,8\aaa_1)\subset 14\aaa_1)$.

Similar to Case 16.
By \cite{Nik11},
we should consider two cases, when \newline
$({\bf n}=22,(2\aaa_1,4\aaa_1,8\aaa_1)\subset 14\aaa_1)$ is marked
by $N_{23}$ and by $N_{21}$.

If $({\bf n}=22,(2\aaa_1,4\aaa_1,8\aaa_1)\subset 14\aaa_1)$ is marked by $N_{23}$, then
$G$ is marked by $G=H_{22,2}$ from Case 12 with orbits
%G:=Group([(3,4)(6,15)(8,11)(9,22)(12,23)(14,20)(16,19)(17,21),
%(2,22)(3,13)(4,18)(7,9)(10,24)(11,20)(15,17)(16,19),
%(2,13)(3,22)(4,9)(7,18)(8,14)(11,20)(12,23)(16,19)]);
%gap> Orbits(G);
%[ [ 2, 22, 13, 9, 3, 7, 4, 18 ], [ 6, 15, 17, 21 ], [ 8, 11, 14, 20 ], [ 10, 24 ], [ 12, 23 ], [ 16, 19 ] ]
%
%$
%\{\alpha_{2},\alpha_{22},\alpha_{13},\alpha_{7},\alpha_{9},
%\alpha_{3},\alpha_{18},\alpha_{4}\},\
%\{\alpha_{6},\alpha_{15},\alpha_{21},\alpha_{17}\},\
%\{\alpha_{8},
%\linebreak
%\alpha_{11},\alpha_{14},\alpha_{20}\},\
%\{\alpha_{10},\alpha_{24}\},\
%\{\alpha_{12},\alpha_{23}\},\  \{\alpha_{16},\alpha_{19}\}$;
$\{\alpha_{12},\alpha_{23}\}$,
$\{\alpha_{6},\alpha_{15},\alpha_{21},\alpha_{17}\}$,
$\{\alpha_{2},\alpha_{22},\alpha_{13},\alpha_{7},\alpha_{9},
\alpha_{3},\alpha_{18},\alpha_{4}\}$.
The $G_1\cong D_8$ is marked by
%gap> H4:=Image(F[4]);
%Group([ (2,9)(3,18)(4,13)(6,21)(7,22)(8,14)(10,24)(16,19), (2,13)(3,9)(4,22)(6,15)(7,18)(8,20)(11,14)(17,21) ])
%gap> Orbits(H4);
%[ [ 2, 9, 13, 3, 4, 18, 22, 7 ], [ 6, 21, 15, 17 ], [ 8, 14, 20, 11 ], [ 10, 24 ], [ 16, 19 ] ]
$G_1\subset G=H_{22,2}$ from Case 15 with suborbits
$\{\alpha_{12}\}$, $\{\alpha_{23}\}$,
$\{\alpha_{6},\alpha_{15},\alpha_{21},\alpha_{17}\}$,
$\{\alpha_{2},\alpha_{22},\alpha_{13},\alpha_{7},\alpha_{9},
\alpha_{3},\alpha_{18},\alpha_{4}\}$.

%N_{21}
%gap> H4:=Image(HH[4]);
%G:=Group([ (7,9)(10,15)(11,14)(12,13)(16,22)(17,23)(18,24)(19,21),(4,6)(7,18)(8,17)(9,16)(10,12)(19,22)(20,23)(21,24),
%(1,3)(4,6)(7,9)(10,12)(13,15)(16,18)(19,21)(22,24) ]);
%gap> Orbits(H4);
%[ [ 1, 3 ], [ 4, 6 ], [ 7, 9, 18, 16, 24, 22, 21, 19 ], [ 8, 17, 23, 20 ], [ 10, 15, 12, 13 ], [ 11, 14 ] ]

If $({\bf n}=22,(2\aaa_1,4\aaa_1,8\aaa_1)\subset 14\aaa_1)$ is marked by $N_{21}$, then
$G$ is marked by $G=H_{22,1}$ from Case 6 with orbits
%$
%\{\alpha_{1,1},\alpha_{3,1}\},\
%\{\alpha_{1,2},\alpha_{3,2}\},\
%\{\alpha_{1,3},\alpha_{3,3},\alpha_{3,6},\alpha_{1,6},\alpha_{3,8},
%\alpha_{1,8},\linebreak
%\alpha_{3,7},\alpha_{1,7}\},\
%\{\alpha_{2,3},\alpha_{2,6},\alpha_{2,8},\alpha_{2,7}\},\
%\{\alpha_{1,4},\alpha_{3,5},\alpha_{3,4},\alpha_{1,5}\},\
%\{\alpha_{2,4},\alpha_{2,5}\}
%$.
$\{\alpha_{1,1},\alpha_{3,1}\}$,
$\{\alpha_{1,4},\alpha_{3,5},\alpha_{3,4},\alpha_{1,5}$\},
$\{\alpha_{1,3},\alpha_{3,3},\alpha_{3,6},\alpha_{1,6},\alpha_{3,8},
\alpha_{1,8},\alpha_{3,7},\alpha_{1,7}\}$.
%gap> H1:=Image(F[1]);
%Group([ (7,9)(10,15)(11,14)(12,13)(16,22)(17,23)(18,24)(19,21), (4,6)(7,18)(8,17)(9,16)(10,12)(19,22)(20,23)(21,24) ])
%gap> Orbits(H1);
%[ [ 4, 6 ], [ 7, 9, 18, 16, 24, 22, 21, 19 ], [ 8, 17, 23, 20 ], [ 10, 15, 12, 13 ], [ 11, 14 ] ]
The $G_1\cong D_8$ is marked by
$G_1\subset G=H_{22,1}$ of Case 16 with suborbits
$\{\alpha_{1,1}\}$, $\{\alpha_{3,1}\}$,
$\{\alpha_{1,4},\alpha_{3,5}, \alpha_{3,4}, \alpha_{1,5}\}$ ,
$\{\alpha_{1,3},\alpha_{3,3},\alpha_{3,6},\alpha_{1,6},\alpha_{3,8},
\alpha_{1,8},\alpha_{3,7},\alpha_{1,7}\}$.

\medskip

Case 27:
$({\bf n}=10,\ (\aaa_1,\,\aaa_1,\,4\aaa_1,\, 8\aaa_1)\subset 2\aaa_1\amalg 4\aaa_3)$

$\Longleftarrow\ ({\bf n}=22,\
\left(\begin{array}{rrr}
 2\aaa_1 & 6\aaa_1 & 10\aaa_1  \\
         & 4\aaa_1 & 4\aaa_3   \\
         &         & 8\aaa_1
\end{array}\right)\subset 2\aaa_1\amalg 4\aaa_3)$.

%gap> H1:=Image(F[1]);
%Group([ (7,9)(10,15)(11,14)(12,13)(16,22)(17,23)(18,24)(19,21),
%(4,6)(7,16,19,24)(8,17,20,23)(9,18,21,22)(10,15,12,13)(11,14) ])
%gap> Orbits(H1);
%[ [ 4, 6 ], [ 7, 9, 16, 18, 22, 19, 24, 21 ], [ 8, 17, 23, 20 ], [ 10, 15, 12, 13 ], [ 11, 14 ] ]
Similar to Cases 1 and 9. The $G\cong C_2\times D_8$
is marked by $N_{21}$ and
$G=H_{22,1}$ from Case 9 with orbits
$\{\alpha_{1,1},\alpha_{3,1}\}$,
$\{\alpha_{2,3},\alpha_{2,6},\alpha_{2,8},\alpha_{2,7}\}$,
$\{\alpha_{1,3},\alpha_{3,3},\alpha_{3,6},\alpha_{1,6},
\alpha_{3,8}, \alpha_{1,8},\alpha_{3,7},\alpha_{1,7}\}$.
The $G_1\cong D_8$ is marked
by
$$
G=H_{22,1}\supset G_1=
$$
$$
[(\alpha_{1,3}\alpha_{3,3})(\alpha_{1,4}\alpha_{3,5})(\alpha_{2,4}\alpha_{2,5})(\alpha_{3,4}\alpha_{1,5})
(\alpha_{1,6}\alpha_{1,8})(\alpha_{2,6}\alpha_{2,8})(\alpha_{3,6}\alpha_{3,8})(\alpha_{1,7}\alpha_{3,7}),
$$
$$
(\alpha_{1,2}\alpha_{3,2})(\alpha_{1,3}\alpha_{1,6}\alpha_{1,7}\alpha_{3,8})
(\alpha_{2,3}\alpha_{2,6}\alpha_{2,7}\alpha_{2,8})
(\alpha_{3,3}\alpha_{3,6}\alpha_{3,7}\alpha_{1,8})
(\alpha_{1,4}\alpha_{3,5}\alpha_{3,4}\alpha_{1,5})(\alpha_{2,4}\alpha_{2,5})]
$$
with suborbits
$\{\alpha_{1,1}\}$, $\{\alpha_{1,3}\}$,
$\{\alpha_{2,3},\alpha_{2,7},\alpha_{2,6},\alpha_{2,8}\}$,
$\{\alpha_{1,3},\alpha_{1,7},\alpha_{3,6},\alpha_{1,6},\alpha_{1,8},\alpha_{3,8},\alpha_{3,3},\alpha_{3,7}\}$.
Both $G$ and $G_1$ are marked by $\ast$ in Tables 1---4.

\medskip

Case 28:

({\bf n}=10,\
$
\left(\begin{array}{cccc}
\aaa_1 & \aaa_3              & 3\aaa_1          &  5\aaa_1     \\
       & (2\aaa_1)_{II}      &  4\aaa_1         &  6\aaa_1     \\
       &                     & (2\aaa_1)_I      &  2\aaa_3     \\
       &                     &                  &  4\aaa_1
\end{array}\right)
\subset 3\aaa_3)\
\Longleftarrow ({\bf n}=34,\ (3\aaa_1,(6\aaa_1)_{II})\subset 3\aaa_3)$.

Similar to Case 1. By \cite{Nik11}, the $G\cong \SSS_4$
is marked by $N_{21}$ and
%gap> Image(HH[3]);
%G:=Group([ (4,7,21)(5,8,20)(6,9,19)(10,24,15)(11,23,14)(12,22,13), %(1,3)(7,21,9,19)(8,20)(10,18,15,24)(11,17,14,23)(12,16,13,22)]);
%gap> Orbits(Image(HH[3]));
%[ [ 1, 3 ], [ 4, 7, 21, 9, 19, 6 ], [ 5, 8, 20 ], [ 10, 24, 18, 15 ], [ 11, 23, 17, 14 ], [ 12, 22, 16, 13 ] ]
$$
G=H_{34,2}=
$$
$$
[(\alpha_{1,2}\alpha_{1,3}\alpha_{3,7})(\alpha_{2,2}\alpha_{2,3}\alpha_{2,7})
(\alpha_{3,2}\alpha_{3,3}\alpha_{1,7})(\alpha_{1,4}\alpha_{3,8}\alpha_{3,5})
(\alpha_{2,4}\alpha_{2,8}\alpha_{2,5})
(\alpha_{3,4}\alpha_{1,8}\alpha_{1,5}),
$$
$$
 (\alpha_{1,1}\alpha_{3,1})
 (\alpha_{1,3}\alpha_{3,7}\alpha_{3,3}\alpha_{1,7})
 (\alpha_{2,3}\alpha_{2,7})
 (\alpha_{1,4}\alpha_{3,6}\alpha_{3,5}\alpha_{3,8})
 (\alpha_{2,4}\alpha_{2,6}\alpha_{2,5}\alpha_{2,8})
 (\alpha_{3,4}\alpha_{1,6}\alpha_{1,5}\alpha_{1,8})]
$$
with orbits
%$
%\{\alpha_{1,1},\alpha_{3,1}\},\
%\{\alpha_{1,2},\alpha_{1,3},\alpha_{3,7},
%\alpha_{3,3},\alpha_{1,7},\alpha_{3,2}\},\
%\{\alpha_{2,2},\alpha_{2,3},\linebreak
%\alpha_{2,7}\},\
%\{\alpha_{1,4},\alpha_{3,8},\alpha_{3,6},\alpha_{3,5}\},\
%\{\alpha_{2,4},\alpha_{2,8},\alpha_{2,6},\alpha_{2,5}\},\
%\{\alpha_{3,4},\alpha_{1,8},\alpha_{1,6},\alpha_{1,5}\};
%$
$\{\alpha_{2,2},\alpha_{2,3},\alpha_{2,7}\}$,
$\{\alpha_{1,2},\alpha_{1,3},\alpha_{3,7},\alpha_{3,3},\alpha_{1,7},\alpha_{3,2}\}$.
The $G_1\cong D_8$ is marked by
%gap> H1;
%Group([ (4,6)(7,9)(10,18)(11,17)(12,16)(13,22)(14,23)(15,24), (1,3)(4,6)(7,19)(8,20)(9,21)(10,15)(11,14)(12,13) ])
%gap> Orbits(H1);
%[ [ 1, 3 ], [ 4, 6 ], [ 7, 9, 19, 21 ], [ 8, 20 ], [ 10, 18, 15, 24 ], [ 11, 17, 14, 23 ], [ 12, 16, 13, 22 ] ]
$$
G=H_{34,2}\supset G_1=
$$
$$
[(\alpha_{1,2}\alpha_{3,2})(\alpha_{1,3}\alpha_{3,3})(\alpha_{1,4}\alpha_{3,6})(\alpha_{2,4}\alpha_{2,6})
(\alpha_{3,4}\alpha_{1,6})(\alpha_{1,5}\alpha_{1,8})(\alpha_{2,5}\alpha_{2,8})(\alpha_{3,5}\alpha_{3,8}),
$$
$$
(\alpha_{1,1}\alpha_{3,1})(\alpha_{1,2}\alpha_{3,2})(\alpha_{1,3}\alpha_{1,7})(\alpha_{2,3}\alpha_{2,7})
(\alpha_{3,3}\alpha_{3,7})(\alpha_{1,4}\alpha_{3,5})(\alpha_{2,4}\alpha_{2,5})(\alpha_{3,4}\alpha_{1,5})]
$$
with suborbits
$\{\alpha_{2,2}\}$, $\{\alpha_{1,2},\alpha_{3,2}\}$,
$\{\alpha_{2,3},\alpha_{2,7}\}$,
$\{\alpha_{1,3},\alpha_{3,3},\alpha_{1,7},\alpha_{3,7}\}$.
Both $G$ and $G_1$ are marked by $\ast$ in Tables 1---4.

\medskip

Case 29:

({\bf n}=10,\
$
\left(\begin{array}{cccc}
\aaa_1 & 3\aaa_1            & 5\aaa_1          &  9\aaa_1     \\
       & (2\aaa_1)_{I}      &  (6\aaa_1)_I     &  10\aaa_1     \\
       &                     & 4\aaa_1         &  12\aaa_1     \\
       &                     &                  &  8\aaa_1
\end{array}\right)
\subset 15\aaa_1)\
\Longleftarrow ({\bf n}=34,\ (3\aaa_1,12\aaa_1)\subset 15\aaa_1)$.

Similar to Case 1. By \cite{Nik11}, the $G\cong \SSS_4$
is marked by $N_{23}$ and
%%G:=Group([(1,16,19,15)(3,5,12,14)(4,9)(7,20,23,22)(8,11,17,18)(10,13),
%(2,13)(3,22)(4,9)(7,18)(8,14)(11,20)(12,23)(16,19)]);
$$
G=H_{34,2}=
$$
$$
[(\alpha_{1}\alpha_{16}\alpha_{19}\alpha_{15})
(\alpha_{3}\alpha_{5}\alpha_{12}\alpha_{14})
(\alpha_{4}\alpha_{9})
(\alpha_{7}\alpha_{20}\alpha_{23}\alpha_{22})
(\alpha_{8}\alpha_{11}\alpha_{17}\alpha_{18})
(\alpha_{10}\alpha_{13}),
$$
$$
(\alpha_{2}\alpha_{13})(\alpha_{3}\alpha_{22})
(\alpha_{4}\alpha_{9})(\alpha_{7}\alpha_{18})
(\alpha_{8}\alpha_{14})(\alpha_{11}\alpha_{20})
(\alpha_{12}\alpha_{23})(\alpha_{16}\alpha_{19})]
$$
with orbits
%$
%\{\alpha_{1},\alpha_{16},\alpha_{19},\alpha_{15}\},\
%\{\alpha_{2},\alpha_{10},\alpha_{13}\},\
%\{\alpha_{3},\alpha_{22},\alpha_{12},\alpha_{18},\alpha_{17},
%\linebreak
%\alpha_{20},\alpha_{7},\alpha_{11},\alpha_{23},\alpha_{8},
%\alpha_{5},\alpha_{14}\},\
%\{\alpha_{4},\alpha_{9}\}$;
$\{\alpha_{2},\alpha_{10},\alpha_{13}\}$,\
$\{\alpha_{3},\alpha_{22},\alpha_{12},\alpha_{18},\alpha_{17},
\alpha_{20},\alpha_{7},\alpha_{11},\alpha_{23},\alpha_{8},
\alpha_{5},\alpha_{14}\}$.
The $G_1\cong D_8$ is marked by
%gap> H1:=Image(F[1]);
%Group([ (3,8)(4,9)(5,18)(10,13)(11,14)(12,17)(15,16)(20,22), (1,15)(3,18)(5,17)(7,22)(8,14)(11,12)(16,19)(20,23) ])
%gap>
%gap> Orbits(H1);
%[ [ 1, 15, 16, 19 ], [ 3, 8, 18, 14, 5, 11, 17, 12 ], [ 4, 9 ], [ 7, 22, 20, 23 ], [ 10, 13 ] ]
$$
G=H_{34,2}\supset G_1=
$$
$$
[(\alpha_{3}\alpha_{8})(\alpha_{4}\alpha_{9})(\alpha_{5}\alpha_{18})(\alpha_{10}\alpha_{13})
(\alpha_{11}\alpha_{14})(\alpha_{12}\alpha_{17})(\alpha_{15}\alpha_{16})(\alpha_{20}\alpha_{22}),
$$
$$
(\alpha_{1}\alpha_{15})(\alpha_{3}\alpha_{18})(\alpha_{5}\alpha_{17})(\alpha_{7}\alpha_{22})
(\alpha_{8}\alpha_{14})(\alpha_{11}\alpha_{12})(\alpha_{16}\alpha_{19})(\alpha_{20}\alpha_{23})]
$$
with suborbits
$\{\alpha_{2}\}$, $\{\alpha_{10},\alpha_{13}\}$,
$\{\alpha_{7},\alpha_{22},\alpha_{20},\alpha_{23}\}$,
$\{\alpha_{3},\alpha_{8},\alpha_{18},\alpha_{14},\alpha_{5},\alpha_{11},\alpha_{17},\alpha_{12}\}$.
Both $G$ and $G_1$ are marked by $\ast$ in Tables 1---4.

\medskip

Case 30:
({\bf n}=10,\
$(\aaa_1,\,4\aaa_1,\,4\aaa_1,\,2\aaa_2)\subset \aaa_1\amalg 6\aaa_2)
\Longleftarrow ({\bf n}=34,\ (\aaa_1,6\aaa_2)\subset \aaa_1\amalg 6\aaa_2))$.

Similar to Cases 1 and 22. By \cite{Nik11}, the $G\cong \SSS_4$ is
marked by by $N_{22}$ and
$G=H_{34,1}$ of Case 22 with orbits
$\{\alpha_{1,6}\}$, $\{\alpha_{1,2},\alpha_{2,5},
\alpha_{1,5},\alpha_{1,11},\alpha_{2,2},\alpha_{2,11},
\alpha_{1,12},\alpha_{2,12},\alpha_{2,7},\alpha_{1,7},\alpha_{1,9},\alpha_{2,9}\}$.
The $G_1\cong D_8$ is marked by
$G_1\subset G=H_{34,1}$ of Case 22
with suborbits
$\{\alpha_{1,6}\}$, $\{\alpha_{1,2},\alpha_{1,11},\alpha_{2,7},\linebreak
\alpha_{2,12}\}$, $\{\alpha_{2,2},\alpha_{2,11},\alpha_{1,7},\alpha_{1,12}\}$,
$\{\alpha_{1,5},\alpha_{1,9},\alpha_{2,9},\alpha_{2,5}\}$.
Both $G$ and $G_1$ are marked by $\ast$ in Tables 1---4.

\medskip
\medskip

Case 31:
({\bf n}=10,\
$
\left(\begin{array}{cccc}
(2\aaa_1)_I  & 4\aaa_1        & (6\aaa_1)_I  &  (6\aaa_1)_{I}     \\
             & (2\aaa_1)_{II} &  6\aaa_1     &  6\aaa_1     \\
       &                      & 4\aaa_1      &  (8\aaa_1)_{II}     \\
       &                      &              &  4\aaa_1
\end{array}\right)
\subset 12\aaa_1)\  \Longleftarrow ({\bf n}=65,\ 12\aaa_1)$.

Similar to Cases 1 and 24. By \cite{Nik10},
the $G\cong 2^4D_6$ is marked by $N_{23}$ and
$G=H_{65,3}$ of Case 24 with the orbit
$\{\alpha_{1},\alpha_{7},\alpha_{23},\alpha_{8},\alpha_{11},\alpha_{17},\alpha_{20},
\alpha_{9}, \alpha_{6},\alpha_{4},\alpha_{14},\alpha_{13}\}$.
%gap> H4:=Image(F[4]);
%Group([ (1,13)(4,11)(5,22)(7,20)(8,17)(10,15)(12,24)(14,23), (3,22)(4,17)(5,16)(6,9)(7,20)(10,19)(13,14)(15,21) ])
%gap> Orbits(H4);
%[ [ 1, 13, 14, 23 ], [ 3, 22, 5, 16 ], [ 4, 11, 17, 8 ], [ 6, 9 ], [ 7, 20 ], [ 10, 15, 19, 21 ], [ 12, 24 ] ]
The $G_1\cong D_8$ is marked by
$$
G=H_{65,3}\supset G_1=
$$
$$
[(\alpha_{1}\alpha_{13})(\alpha_{4}\alpha_{11})(\alpha_{5}\alpha_{22})(\alpha_{7}\alpha_{20})
(\alpha_{8}\alpha_{17})(\alpha_{10}\alpha_{15})(\alpha_{12}\alpha_{24})(\alpha_{14}\alpha_{23}),
$$
$$
(\alpha_{3}\alpha_{22})(\alpha_{4}\alpha_{17})(\alpha_{5}\alpha_{16})(\alpha_{6}\alpha_{9})
(\alpha_{7}\alpha_{20})(\alpha_{10}\alpha_{19})(\alpha_{13}\alpha_{14})(\alpha_{15}\alpha_{21})]
$$
with suborbits
$\{\alpha_{6},\alpha_{9}\}$, $\{\alpha_{7},\alpha_{20}\}$, $\{\alpha_{1},\alpha_{13},\alpha_{14},\alpha_{23}\}$,
$\{\alpha_{4},\alpha_{11},\alpha_{17},\alpha_{8}\}$ .
The orbit $\{\alpha_{6},\alpha_{9}\}$ has the type $(2\aaa_1)_I$, and the orbit
$\{\alpha_{7},\alpha_{20}\}$ has the type $(2\aaa_1)_{II}$ for $G_1$.
Both $G$ and $G_1$ are marked by $\ast$ in Tables 1---4.

\medskip

Case 32:
({\bf n}=10,\
$
\left(\begin{array}{cccc}
(2\aaa_1)_I & (6\aaa_1)_{II}   & (6\aaa_1)_I     &  (6\aaa_1)_I     \\
            & 4\aaa_1          & (8\aaa_1)_{I}   &  (8\aaa_1)_I     \\
       &                       & 4\aaa_1         &   (8\aaa_1)_{II} \\
       &                       &                 &  4\aaa_1
\end{array}\right)
\subset 14\aaa_1)$

$\Longleftarrow\ ({\bf n}=22,\ (2\aaa_1,4\aaa_1,8\aaa_1)\subset 14\aaa_1)$.

Similar to Case 16. By \cite{Nik11},
we should consider two cases, when \newline
$({\bf n}=22,(2\aaa_1,4\aaa_1,8\aaa_1)\subset 14\aaa_1)$ is marked
by $N_{23}$ and by $N_{21}$.

If $({\bf n}=22,(2\aaa_1,4\aaa_1,8\aaa_1)\subset 14\aaa_1)$ is marked by $N_{23}$, then
$G$ is marked by $G=H_{22,2}$ from Case 12 with orbits
%G:=Group([(3,4)(6,15)(8,11)(9,22)(12,23)(14,20)(16,19)(17,21),
%(2,22)(3,13)(4,18)(7,9)(10,24)(11,20)(15,17)(16,19),
%(2,13)(3,22)(4,9)(7,18)(8,14)(11,20)(12,23)(16,19)]);
%gap> Orbits(G);
%[ [ 2, 22, 13, 9, 3, 7, 4, 18 ], [ 6, 15, 17, 21 ], [ 8, 11, 14, 20 ], [ 10, 24 ], [ 12, 23 ], [ 16, 19 ] ]
%
%$
%\{\alpha_{2},\alpha_{22},\alpha_{13},\alpha_{7},\alpha_{9},
%\alpha_{3},\alpha_{18},\alpha_{4}\},\
%\{\alpha_{6},\alpha_{15},\alpha_{21},\alpha_{17}\},\
%\{\alpha_{8},
%\linebreak
%\alpha_{11},\alpha_{14},\alpha_{20}\},\
%\{\alpha_{10},\alpha_{24}\},\
%\{\alpha_{12},\alpha_{23}\},\  \{\alpha_{16},\alpha_{19}\}$;
$\{\alpha_{12},\alpha_{23}\}$,
$\{\alpha_{6},\alpha_{15},\alpha_{21},\alpha_{17}\}$,
$\{\alpha_{2},\alpha_{22},\alpha_{13},\alpha_{7},\alpha_{9},
\alpha_{3},\alpha_{18},\alpha_{4}\}$,
and  $G_1\cong D_8$ is marked by
%gap> H2:=Image(F[2]);
%Group([ (3,4)(6,15)(8,11)(9,22)(12,23)(14,20)(16,19)(17,21), (2,9)(3,18)(4,13)(6,21)(7,22)(8,14)(10,24)(16,19) ])
%gap> Orbits(H2);
%[ [ 2, 9, 22, 7 ], [ 3, 4, 18, 13 ], [ 6, 15, 21, 17 ], [ 8, 11, 14, 20 ], [ 10, 24 ], [ 12, 23 ], [ 16, 19 ] ]
$$
G=H_{22,2}\supset G_1
$$
$$
[(\alpha_{3}\alpha_{4})(\alpha_{6}\alpha_{15})(\alpha_{8}\alpha_{11})(\alpha_{9}\alpha_{22})
(\alpha_{12}\alpha_{23})(\alpha_{14}\alpha_{20})(\alpha_{16}\alpha_{19})(\alpha_{17}\alpha_{21}),
$$
$$
(\alpha_{2}\alpha_{9})(\alpha_{3}\alpha_{18})(\alpha_{4}\alpha_{13})(\alpha_{6}\alpha_{21})
(\alpha_{7}\alpha_{22})(\alpha_{8}\alpha_{14})(\alpha_{10}\alpha_{24})(\alpha_{16}\alpha_{19})]
$$
with suborbits
$\{\alpha_{12},\alpha_{23},\}$, $\{\alpha_{6},\alpha_{15},\alpha_{21},\alpha_{17}\}$,
$\{\alpha_{2},\alpha_{9},\alpha_{22},\alpha_{7}\}$, $\{\alpha_{3},\alpha_{4},\alpha_{18},\alpha_{13}\}$.

%N_{21}
%gap> H4:=Image(HH[4]);
%G:=Group([ (7,9)(10,15)(11,14)(12,13)(16,22)(17,23)(18,24)(19,21),(4,6)(7,18)(8,17)(9,16)(10,12)(19,22)(20,23)(21,24),
%(1,3)(4,6)(7,9)(10,12)(13,15)(16,18)(19,21)(22,24) ]);
%gap> Orbits(H4);
%[ [ 1, 3 ], [ 4, 6 ], [ 7, 9, 18, 16, 24, 22, 21, 19 ], [ 8, 17, 23, 20 ], [ 10, 15, 12, 13 ], [ 11, 14 ] ]

If $({\bf n}=22,(2\aaa_1,8\aaa_1)\subset 10\aaa_1)$ is marked by $N_{21}$, then
$G$ is marked by $G=H_{22,1}$ from Case 6 with orbits
%$
%\{\alpha_{1,1},\alpha_{3,1}\},\
%\{\alpha_{1,2},\alpha_{3,2}\},\
%\{\alpha_{1,3},\alpha_{3,3},\alpha_{3,6},\alpha_{1,6},\alpha_{3,8},
%\alpha_{1,8},\linebreak
%\alpha_{3,7},\alpha_{1,7}\},\
%\{\alpha_{2,3},\alpha_{2,6},\alpha_{2,8},\alpha_{2,7}\},\
%\{\alpha_{1,4},\alpha_{3,5},\alpha_{3,4},\alpha_{1,5}\},\
%\{\alpha_{2,4},\alpha_{2,5}\}
%$.
$\{\alpha_{1,1},\alpha_{3,1}\}$,
$\{\alpha_{1,4},\alpha_{3,5},\alpha_{3,4},\alpha_{1,5}$\},
$\{\alpha_{1,3},\alpha_{3,3},\alpha_{3,6},\alpha_{1,6},\alpha_{3,8},
\alpha_{1,8},\alpha_{3,7},\alpha_{1,7}\}$,
and  $G_1\cong D_8$ is marked by
%gap> H3:=Image(F[3]);
%Group([ (4,6)(7,18)(8,17)(9,16)(10,12)(19,22)(20,23)(21,24), (1,3)(4,6)(10,13)(11,14)(12,15)(16,24)(17,23)(18,22) ])
%gap> Orbits(H3);
%[ [ 1, 3 ], [ 4, 6 ], [ 7, 18, 22, 19 ], [ 8, 17, 23, 20 ], [ 9, 16, 24, 21 ], [ 10, 12, 13, 15 ], [ 11, 14 ] ]
$$
G=H_{22,1}\supset G_1=
$$
$$
[(\alpha_{1,2}\alpha_{3,2})(\alpha_{1,3}\alpha_{3,6})(\alpha_{2,3}\alpha_{2,6})(\alpha_{3,3}\alpha_{1,6})
(\alpha_{1,4}\alpha_{3,4})(\alpha_{1,7}\alpha_{1,8})(\alpha_{2,7}\alpha_{2,8})(\alpha_{3,7}\alpha_{3,8}),
$$
$$
(\alpha_{1,1}\alpha_{3,1})(\alpha_{1,2}\alpha_{3,2})(\alpha_{1,4}\alpha_{1,5})(\alpha_{2,4}\alpha_{2,5})
(\alpha_{3,4}\alpha_{3,5})(\alpha_{1,6}\alpha_{3,8})(\alpha_{2,6}\alpha_{2,8})(\alpha_{3,6}\alpha_{1,8})]
$$
with suborbits
$\{\alpha_{1,1}\}$,$\{\alpha_{3,1}\}$,$\{\alpha_{1,4},\alpha_{3,5}, \alpha_{3,4}, \alpha_{1,5}\}$,$\{\alpha_{1,3},\alpha_{3,6},\alpha_{1,8},\alpha_{1,7}\}$,$\{\alpha_{3,3},\alpha_{1,6},
\alpha_{3,8},\alpha_{3,7}\}$.

\medskip

Case 33:
({\bf n}=10,\
$
\left(\begin{array}{cccc}
(2\aaa_1)_{II} & 6\aaa_1     & 6\aaa_1          &  6\aaa_1        \\
               & 4\aaa_1     &(8\aaa_1)_I       &  (8\aaa_1)_I    \\
               &             & 4\aaa_1          &  (8\aaa_1)_{II} \\
       &                     &                  &  4\aaa_1
\end{array}\right)
\subset 14\aaa_1)
$

$\Longleftarrow\ ({\bf n}=22,\ (2\aaa_1,4\aaa_1,8\aaa_1)\subset 14\aaa_1)$.

Similar to Case 16.
By \cite{Nik11},
we should consider two cases, when \newline
$({\bf n}=22,(2\aaa_1,4\aaa_1,8\aaa_1)\subset 14\aaa_1)$ is marked
by $N_{23}$ and by $N_{21}$.

If $({\bf n}=22,(2\aaa_1,4\aaa_1,8\aaa_1)\subset 14\aaa_1)$ is marked by $N_{23}$, then
$G$ is marked by $G=H_{22,2}$ from Case 12 with orbits
%G:=Group([(3,4)(6,15)(8,11)(9,22)(12,23)(14,20)(16,19)(17,21),
%(2,22)(3,13)(4,18)(7,9)(10,24)(11,20)(15,17)(16,19),
%(2,13)(3,22)(4,9)(7,18)(8,14)(11,20)(12,23)(16,19)]);
%gap> Orbits(G);
%[ [ 2, 22, 13, 9, 3, 7, 4, 18 ], [ 6, 15, 17, 21 ], [ 8, 11, 14, 20 ], [ 10, 24 ], [ 12, 23 ], [ 16, 19 ] ]
%
%$
%\{\alpha_{2},\alpha_{22},\alpha_{13},\alpha_{7},\alpha_{9},
%\alpha_{3},\alpha_{18},\alpha_{4}\},\
%\{\alpha_{6},\alpha_{15},\alpha_{21},\alpha_{17}\},\
%\{\alpha_{8},
%\linebreak
%\alpha_{11},\alpha_{14},\alpha_{20}\},\
%\{\alpha_{10},\alpha_{24}\},\
%\{\alpha_{12},\alpha_{23}\},\  \{\alpha_{16},\alpha_{19}\}$;
$\{\alpha_{12},\alpha_{23}\}$,
$\{\alpha_{6},\alpha_{15},\alpha_{21},\alpha_{17}\}$,
$\{\alpha_{2},\alpha_{22},\alpha_{13},\alpha_{7},\alpha_{9},
\alpha_{3},\alpha_{18},\alpha_{4}\}$,
and $G_1\cong D_8$ is marked by
%gap> H1:=Image(F[1]);
%Group([ (3,4)(6,15)(8,11)(9,22)(12,23)(14,20)(16,19)(17,21), (2,3)(4,7)(8,14)(9,18)(10,24)(12,23)(13,22)(15,17) ])
%gap> Orbits(H1);
%[ [ 2, 3, 4, 7 ], [ 6, 15, 17, 21 ], [ 8, 11, 14, 20 ], [ 9, 22, 18, 13 ], [ 10, 24 ], [ 12, 23 ], [ 16, 19 ] ]
$$
G=H_{22,2}\supset G_1=
$$
$$
[(\alpha_{3}\alpha_{4})(\alpha_{6}\alpha_{15})(\alpha_{8}\alpha_{11})(\alpha_{9}\alpha_{22})
(\alpha_{12}\alpha_{23})(\alpha_{14}\alpha_{20})(\alpha_{16}\alpha_{19})(\alpha_{17}\alpha_{21}),
$$
$$
(\alpha_{2}\alpha_{3})(\alpha_{4}\alpha_{7})(\alpha_{8}\alpha_{14})(\alpha_{9}\alpha_{18})
(\alpha_{10}\alpha_{24})(\alpha_{12}\alpha_{23})(\alpha_{13}\alpha_{22})(\alpha_{15}\alpha_{17})]
$$
with suborbits
$\{\alpha_{12},\alpha_{23},\}$,$\{\alpha_{6},\alpha_{15},\alpha_{17},\alpha_{21}\}$,
$\{\alpha_{2},\alpha_{3},\alpha_{4},\alpha_{7}\}$, $\{\alpha_{9},\alpha_{22}, \alpha_{18},\alpha_{13}\}$.

%N_{21}
%gap> H4:=Image(HH[4]);
%G:=Group([ (7,9)(10,15)(11,14)(12,13)(16,22)(17,23)(18,24)(19,21),(4,6)(7,18)(8,17)(9,16)(10,12)(19,22)(20,23)(21,24),
%(1,3)(4,6)(7,9)(10,12)(13,15)(16,18)(19,21)(22,24) ]);
%gap> Orbits(H4);
%[ [ 1, 3 ], [ 4, 6 ], [ 7, 9, 18, 16, 24, 22, 21, 19 ], [ 8, 17, 23, 20 ], [ 10, 15, 12, 13 ], [ 11, 14 ] ]

If $({\bf n}=22,(2\aaa_1,8\aaa_1)\subset 10\aaa_1)$ is marked by $N_{21}$, then
$G$ is marked by $G=H_{22,1}$ from Case 6 with orbits
%$
%\{\alpha_{1,1},\alpha_{3,1}\},\
%\{\alpha_{1,2},\alpha_{3,2}\},\
%\{\alpha_{1,3},\alpha_{3,3},\alpha_{3,6},\alpha_{1,6},\alpha_{3,8},
%\alpha_{1,8},\linebreak
%\alpha_{3,7},\alpha_{1,7}\},\
%\{\alpha_{2,3},\alpha_{2,6},\alpha_{2,8},\alpha_{2,7}\},\
%\{\alpha_{1,4},\alpha_{3,5},\alpha_{3,4},\alpha_{1,5}\},\
%\{\alpha_{2,4},\alpha_{2,5}\}
%$.
$\{\alpha_{1,1},\alpha_{3,1}\}$,
$\{\alpha_{1,4},\alpha_{3,5},\alpha_{3,4},\alpha_{1,5}$\},
$\{\alpha_{1,3},\alpha_{3,3},\alpha_{3,6},\alpha_{1,6},\alpha_{3,8},
\alpha_{1,8},\alpha_{3,7},\alpha_{1,7}\}$,
and $G_1\cong D_8$ is marked by
%gap> H4:=Image(F[4]);
%Group([ (1,3)(7,16)(8,17)(9,18)(13,15)(19,24)(20,23)(21,22), (1,3)(4,6)(10,13)(11,14)(12,15)(16,24)(17,23)(18,22) ])
%gap> Orbits(H4);
%[ [ 1, 3 ], [ 4, 6 ], [ 7, 16, 24, 19 ], [ 8, 17, 23, 20 ], [ 9, 18, 22, 21 ], [ 10, 13, 15, 12 ], [ 11, 14 ] ]
$$
G=H_{22,1}\supset G_1=
$$
$$
[(\alpha_{1,1}\alpha_{3,1})(\alpha_{1,3}\alpha_{1,6})(\alpha_{2,3}\alpha_{2,6})(\alpha_{3,3}\alpha_{3,6})
(\alpha_{1,5}\alpha_{3,5})(\alpha_{1,7}\alpha_{3,8})(\alpha_{2,7}\alpha_{2,8})(\alpha_{3,7}\alpha_{1,8}),
$$
$$
(\alpha_{1,1}\alpha_{3,1})(\alpha_{1,2}\alpha_{3,2})(\alpha_{1,4}\alpha_{1,5})(\alpha_{2,4}\alpha_{2,5})
(\alpha_{3,4}\alpha_{3,5})(\alpha_{1,6}\alpha_{3,8})(\alpha_{2,6}\alpha_{2,8})(\alpha_{3,6}\alpha_{1,8})]
$$
with suborbits
$\{\alpha_{1,1}\}$, $\{\alpha_{3,1}\}$,$\{\alpha_{1,4},\alpha_{3,5}, \alpha_{3,4}, \alpha_{1,5}\}$,$\{\alpha_{1,3},\alpha_{1,6},\alpha_{3,8},\alpha_{1,7}\}$,$\{\alpha_{3,3},\alpha_{3,6},
\alpha_{1,8},\alpha_{3,7}\}$.

\medskip

Case 34:

({\bf n}=10,\
$
\left(\begin{array}{cccc}
(2\aaa_1)_{I} & (6\aaa_1)_I  & (6\aaa_1)_I   & (6\aaa_1)_I      \\
               & 4\aaa_1     & 4\aaa_2       &  (8\aaa_1)_{II}  \\
               &             & 4\aaa_1       &  4\aaa_2         \\
       &                     &                  &  4\aaa_1
\end{array}\right)
\subset 2\aaa_1\amalg 4\aaa_3)
$,\

({\bf n}=10,\
$
\left(\begin{array}{cccc}
(2\aaa_1)_{II} & 6\aaa_1     & 6\aaa_1       & 6\aaa_1      \\
               & 4\aaa_1     & 4\aaa_2       &  (8\aaa_1)_{II}  \\
               &             & 4\aaa_1       &  4\aaa_2         \\
       &                     &                  &  4\aaa_1
\end{array}\right)
\subset 2\aaa_1\amalg 4\aaa_3)
$

$\Longleftarrow\ ({\bf n}=22,\
\left(\begin{array}{rrr}
 2\aaa_1 & 6\aaa_1 & 10\aaa_1  \\
         & 4\aaa_1 & 4\aaa_3   \\
         &         & 8\aaa_1
\end{array}\right)\subset 2\aaa_1\amalg 4\aaa_3)$.

\medskip

Similar to Cases 1 and 9. By \cite{Nik11}, the $G$
is marked by $N_{21}$ and
$G=H_{22,1}$ from Case 9 with orbits
$\{\alpha_{1,1},\alpha_{3,1}\}$,
$\{\alpha_{2,3},\alpha_{2,6},\alpha_{2,8},\alpha_{2,7}\}$,
$\{\alpha_{1,3},\alpha_{3,3},\alpha_{3,6},\alpha_{1,6},
\alpha_{3,8}, \alpha_{1,8},\alpha_{3,7},\alpha_{1,7}\}$.
The $G_1\cong D_8$ is marked by
$$
G=H_{22,1}\supset G_1=
$$
%gap> H3:=Image(F[3]);
%Group([ (4,6)(7,18)(8,17)(9,16)(10,12)(19,22)(20,23)(21,24), (1,3)(7,18,19,22)(8,17,20,23)(9,16,21,24)(10,13,12,15)
%(11,14) ]);
$$
[(\alpha_{1,2}\alpha_{3,2})(\alpha_{1,3}\alpha_{3,6})(\alpha_{2,3}\alpha_{2,6})(\alpha_{3,3}\alpha_{1,6})
(\alpha_{1,4}\alpha_{3,4})(\alpha_{1,7}\alpha_{1,8})(\alpha_{2,7}\alpha_{2,8})(\alpha_{3,7}\alpha_{3,8}),
$$
$$
(\alpha_{1,1}\alpha_{3,1})
(\alpha_{1,3}\alpha_{3,6}\alpha_{1,7}\alpha_{1,8})(\alpha_{2,3}\alpha_{2,6}\alpha_{2,7}\alpha_{2,8})
(\alpha_{3,3}\alpha_{1,6}\alpha_{3,7}\alpha_{3,8})(\alpha_{1,4}\alpha_{1,5}\alpha_{3,4}\alpha_{3,5})
(\alpha_{2,4}\alpha_{2,5})]
$$
with suborbits
%gap> Orbits(H3);
%[ [ 1, 3 ], [ 8, 17, 20, 23 ],  [ 7, 18, 19, 22 ], [ 9, 16, 21, 24 ]]
$\{\alpha_{1,1},\alpha_{3,1}\}$, $\{\alpha_{2,3}\alpha_{2,6},\alpha_{2,7},\alpha_{2,8}\}$,
$\{\alpha_{1,3}\alpha_{3,6},\alpha_{1,7},\alpha_{1,8}\}$,
$\{\alpha_{3,3}\alpha_{1,6},\alpha_{3,7},\alpha_{3,8}\}$ for the first case, and by
$$
G=H_{22,1}\supset G_1=
$$
%gap> H4:=Image(F[4]);
%Group([ (1,3)(7,16)(8,17)(9,18)(13,15)(19,24)(20,23)(21,22), (4,6)(7,16,19,24)(8,17,20,23)(9,18,21,22)(10,15,12,13)
%(11,14) ]);
$$
[(\alpha_{1,1}\alpha_{3,1})(\alpha_{1,3}\alpha_{1,6})(\alpha_{2,3}\alpha_{2,6})(\alpha_{3,3}\alpha_{3,6})
(\alpha_{1,5}\alpha_{3,5})(\alpha_{1,7}\alpha_{3,8})(\alpha_{2,7}\alpha_{2,8})(\alpha_{3,7}\alpha_{1,8}),
$$
$$
(\alpha_{1,2}\alpha_{3,2})(\alpha_{1,3}\alpha_{1,6}\alpha_{1,7}\alpha_{3,8})
(\alpha_{2,3}\alpha_{2,6}\alpha_{2,7}\alpha_{2,8})
(\alpha_{3,3}\alpha_{3,6}\alpha_{3,7}\alpha_{1,8})
(\alpha_{1,4}\alpha_{3,5}\alpha_{3,4}\alpha_{1,5})(\alpha_{2,4}\alpha_{2,5})]
$$
with suborbits
%gap> Orbits(H4);
%[ [ 1, 3 ], [ 8, 17, 20, 23 ], [ 7, 16, 19, 24 ], [ 9, 18, 21, 22 ]]
$\{\alpha_{1,1},\alpha_{3,1}\}$,
$\{\alpha_{2,3},\alpha_{2,6},\alpha_{2,7},\alpha_{2,8}\}$,
$\{\alpha_{1,3},\alpha_{1,6},\alpha_{1,7},\alpha_{3,8}\}$,
$\{\alpha_{3,3},\alpha_{3,6},\alpha_{3,7},\alpha_{1,8}\}$
for the second case.
All these $G$ and $G_1$ are marked by $\ast$ in Tables 1---4.

\medskip

Case 35:

({\bf n}=10,\
$
\left(\begin{array}{cccc}
(2\aaa_1)_I & 2\aaa_3          & (6\aaa_1)_I     &  10\aaa_1  \\
            & 4\aaa_1          & (8\aaa_1)_{I}   &  12\aaa_1  \\
       &                       & 4\aaa_1         &  4\aaa_3   \\
       &                       &                 &  8\aaa_1
\end{array}\right)
\subset 6\aaa_3) \
\Longleftarrow ({\bf n}=34,\ ((6\aaa_1)_I,12\aaa_1)\subset 6\aaa_3)$.

Similar to Case 1. By \cite{Nik11}, the $G\cong \SSS_4$
is marked by $N_{21}$ and
%gap> Image(HH[2]);
%G:=Group([ (7,12,22)(8,11,23)(9,10,24)(13,18,21)(14,17,20)(15,16,19),
%(1,3)(7,21,9,19)(8,20)(10,18,15,24)(11,17,14,23)(12,16,13,22) ]);
%gap> Orbits(Image(HH[2]));
%[ [ 1, 3 ], [ 7, 12, 21, 22, 16, 13, 9, 19, 18, 10, 15, 24 ], [ 8, 11, 20, 23, 17, 14 ] ]
$$
G=H_{34,1}=
$$
$$
[(\alpha_{1,3}\alpha_{3,4}\alpha_{1,8})
(\alpha_{2,3}\alpha_{2,4}\alpha_{2,8})
(\alpha_{3,3}\alpha_{1,4}\alpha_{3,8})
(\alpha_{1,5}\alpha_{3,6}\alpha_{3,7})
(\alpha_{2,5}\alpha_{2,6}\alpha_{2,7})
(\alpha_{3,5}\alpha_{1,6}\alpha_{1,7}),
$$
$$
(\alpha_{1,1}\alpha_{3,1})
(\alpha_{1,3}\alpha_{3,7}\alpha_{3,3}\alpha_{1,7})
(\alpha_{2,3}\alpha_{2,7})
(\alpha_{1,4}\alpha_{3,6}\alpha_{3,5}\alpha_{3,8})
(\alpha_{2,4}\alpha_{2,6}\alpha_{2,5}\alpha_{2,8})
(\alpha_{3,4}\alpha_{1,6}\alpha_{1,5}\alpha_{1,8})]
$$
with orbits
%$
%\{\alpha_{1,1}, \alpha_{3,1} \},\  \{\alpha_{1,3}, \alpha_{3,4},
%\alpha_{3,7}, \alpha_{1,8}, \alpha_{1,6}, \alpha_{1,5},
%\alpha_{3,3}, \alpha_{1,7}, \linebreak
%\alpha_{3,6},
%\alpha_{1,4}, \alpha_{3,5}, \alpha_{3,8}\},\  \{\alpha_{2,3}, \alpha_{2,4},
%\alpha_{2,7}, \alpha_{2,8}, \alpha_{2,6}, \alpha_{2,5} \}
%$;
$\{\alpha_{2,3},\alpha_{2,4},\alpha_{2,7},\alpha_{2,8},\alpha_{2,6},\alpha_{2,5} \}$,
$\{\alpha_{1,3}, \alpha_{3,4},\alpha_{3,7}, \alpha_{1,8}, \alpha_{1,6}, \alpha_{1,5},
\alpha_{3,3}, \alpha_{1,7}, \alpha_{3,6}, \alpha_{1,4}, \linebreak \alpha_{3,5}, \alpha_{3,8}\}$.
%gap> H1:=Image(F[1]);
%Group([ (7,19)(8,20)(9,21)(10,12)(13,15)(16,24)(17,23)(18,22), (1,3)(10,16)(11,17)(12,18)(13,24)(14,23)(15,22)
%(19,21) ])
%gap> Orbits(H1);
%[ [ 1, 3 ], [ 7, 19, 21, 9 ], [ 8, 20 ], [ 10, 12, 16, 18, 24, 22, 13, 15 ], [ 11, 17, 23, 14 ] ]
The $G_1\cong D_8$ is marked by
$$
G=H_{34,1}\supset G_1=
$$
$$
[(\alpha_{1,3}\alpha_{1,7})(\alpha_{2,3}\alpha_{2,7})(\alpha_{3,3}\alpha_{3,7})(\alpha_{1,4}\alpha_{3,4})
(\alpha_{1,5}\alpha_{3,5})(\alpha_{1,6}\alpha_{3,8})(\alpha_{2,6}\alpha_{2,8})(\alpha_{3,6}\alpha_{1,8}),
$$
$$
(\alpha_{1,1}\alpha_{3,1})(\alpha_{1,4}\alpha_{1,6})(\alpha_{2,4}\alpha_{2,6})(\alpha_{3,4}\alpha_{3,6})
(\alpha_{1,5}\alpha_{3,8})(\alpha_{2,5}\alpha_{2,8})(\alpha_{3,5}\alpha_{1,8})(\alpha_{1,7}\alpha_{3,7})]
$$
with suborbits
$\{\alpha_{2,3},\alpha_{2,7}\}$,
$\{\alpha_{2,4},\alpha_{2,6},\alpha_{2,8},\alpha_{2,5}\}$,
$\{\alpha_{1,3},\alpha_{1,7},\alpha_{3,7},\alpha_{3,3}\}$,
$\{\alpha_{1,4},\alpha_{3,4},\alpha_{1,6},\alpha_{3,6}, \linebreak
\alpha_{3,8},\alpha_{1,8},\alpha_{1,5},\alpha_{3,5}\}$.
Both $G$ and $G_1$ are marked by $\ast$ in Tables 1---4.

\medskip

Case 36:

({\bf n}=10,\
$
\left(\begin{array}{cccc}
4\aaa_1 & (8\aaa_1)_{II}     & (8\aaa_1)_I        &  (8\aaa_1)_I  \\
              & 4\aaa_1            & (8\aaa_1)_I  &  (8\aaa_1)_I  \\
               &                   & 4\aaa_1      &  (8\aaa_1)_{II}\\
       &                     &                    &  4\aaa_1
\end{array}\right)
\subset 16\aaa_1) \Longleftarrow ({\bf n}=56,\ 16\aaa_1)$.

Similar to Cases 1 and 8. By \cite{Nik10}, the $G\cong \Gamma_{25}a_1$
is marked by $N_{23}$ and
$G=H_{56,2}$ from Case 8 with the orbit
$\{\alpha_{2}, \alpha_{3}, \alpha_{23},
\alpha_{24}, \alpha_{5},\alpha_{8}, \alpha_{17}, \alpha_{19}, \alpha_{6},
\alpha_{7}, \alpha_{11}, \alpha_{20}, \alpha_{18},
\alpha_{10},\alpha_{16},\alpha_{15}\}$.
The $G_1\cong D_8$ is marked by
$$
G=H_{56,2}\supset G_1=
$$
%gap> H5:=Image(F[5]);
%Group([ (2,17)(6,16)(7,8)(9,13)(10,18)(12,22)(19,20)(23,24), (1,14)(3,7)(5,8)(9,13)(10,24)(11,20)(15,19)(16,17) ])
%gap>
%gap> Orbits(H5);
%[ [ 1, 14 ], [ 2, 17, 16, 6 ], [ 3, 7, 8, 5 ], [ 9, 13 ], [ 10, 18, 24, 23 ], [ 11, 20, 19, 15 ], [ 12, 22 ] ]
$$
[(\alpha_{2}\alpha_{17})(\alpha_{6}\alpha_{16})(\alpha_{7}\alpha_{8})(\alpha_{9}\alpha_{13})
(\alpha_{10}\alpha_{18})(\alpha_{12}\alpha_{22})(\alpha_{19}\alpha_{20})(\alpha_{23}\alpha_{24}),
$$
$$
(\alpha_{1}\alpha_{14})(\alpha_{3}\alpha_{7})(\alpha_{5}\alpha_{8})(\alpha_{9}\alpha_{13})
(\alpha_{10}\alpha_{24})(\alpha_{11}\alpha_{20})(\alpha_{15}\alpha_{19})(\alpha_{16}\alpha_{17})]
$$
with suborbits
$\{\alpha_{2},\alpha_{17},\alpha_{16},\alpha_{6}\}$, 
$\{\alpha_{3},\alpha_{7},\alpha_{8},\alpha_{5}\}$,
$\{\alpha_{10},\alpha_{18},\alpha_{24},\alpha_{23}\}$, 
$\{\alpha_{11},\alpha_{20},\alpha_{19},\alpha_{15}\}$. \linebreak 
Both $G$ and $G_1$ are marked by $\ast$ in Tables 1---4.

\medskip

Case 37:
({\bf n}=10,\
$
\left(\begin{array}{cccc}
4\aaa_1 & (8\aaa_1)_I        & (8\aaa_1)_I      &  (8\aaa_1)_I   \\
               & 4\aaa_1     & 4\aaa_2          &  (8\aaa_1)_{II}\\
               &             & 4\aaa_1          &  4\aaa_2       \\
       &                     &                  &  4\aaa_1
\end{array}\right)
\subset 4\aaa_1\amalg 4\aaa_3)
$

$\Longleftarrow
({\bf n}=22,\
\left(\begin{array}{rrr}
 4\aaa_1 & (8\aaa_1)_{II} & 12\aaa_1  \\
         & 4\aaa_1 & 4\aaa_3   \\
         &         & 8\aaa_1
\end{array}\right)\subset 4\aaa_1\amalg 4\aaa_3)$.

%G:=Group([ (7,9)(10,15)(11,14)(12,13)(16,22)(17,23)(18,24)(19,21),(4,6)(7,18)(8,17)(9,16)(10,12)(19,22)(20,23)(21,24),
%(1,3)(4,6)(7,9)(10,12)(13,15)(16,18)(19,21)(22,24) ]);
%gap> Orbits(H4);
%[ [ 1, 3 ], [ 4, 6 ], [ 7, 9, 18, 16, 24, 22, 21, 19 ], [ 8, 17, 23, 20 ], [ 10, 15, 12, 13 ], [ 11, 14 ] ]

Similar to Cases 1 and 11. By \cite{Nik11},
the $G\cong C_2\times D_8$ is marked
by $N_{21}$ and $G=H_{22,1}$ from Case 11
with orbits
%[ 10, 15, 12, 13 ], [ 8, 17, 23, 20 ],[ 7, 9, 18, 16, 24, 22, 21, 19 ],
$\{\alpha_{1,4},\alpha_{3,5},\alpha_{3,4},\alpha_{1,5}\}$,
$\{\alpha_{2,3},\alpha_{2,6},\alpha_{2,8},\alpha_{2,7}\}$,
$\{\alpha_{1,3},\alpha_{3,3},\alpha_{3,6},\alpha_{1,6}, \linebreak \alpha_{3,8},
\alpha_{1,8},\alpha_{3,7},\alpha_{1,7}\}$.
%gap> H3:=Image(F[3]);
%Group([ (4,6)(7,18)(8,17)(9,16)(10,12)(19,22)(20,23)(21,24), (1,3)(7,18,19,22)(8,17,20,23)(9,16,21,24)(10,13,12,15)
%(11,14) ])
%gap> Orbits(H3);
%[ [ 1, 3 ], [ 4, 6 ], [ 7, 18, 19, 22 ], [ 8, 17, 20, 23 ], [ 9, 16, 21, 24 ], [ 10, 12, 13, 15 ], [ 11, 14 ] ]
The $G_1\cong D_8$ is marked by
$$
G=H_{22,1}\supset G_1=
$$
$$
[(\alpha_{1,2}\alpha_{3,2})(\alpha_{1,3}\alpha_{3,6})(\alpha_{2,3}\alpha_{2,6})(\alpha_{3,3}\alpha_{1,6})
(\alpha_{1,4}\alpha_{3,4})(\alpha_{1,7}\alpha_{1,8})(\alpha_{2,7}\alpha_{2,8})(\alpha_{3,7}\alpha_{3,8}),
$$
$$
(\alpha_{1,1}\alpha_{3,1})
(\alpha_{1,3}\alpha_{3,6}\alpha_{1,7}\alpha_{1,8})(\alpha_{2,3}\alpha_{2,6}\alpha_{2,7}\alpha_{2,8})
(\alpha_{3,3}\alpha_{1,6}\alpha_{3,7}\alpha_{3,8})(\alpha_{1,4}\alpha_{1,5}\alpha_{3,4}\alpha_{3,5})
(\alpha_{2,4}\alpha_{2,5})]
$$
with suborbits
$\{\alpha_{1,4},\alpha_{3,4},\alpha_{1,5},\alpha_{3,5}\}$,
$\{\alpha_{1,3},\alpha_{3,6},\alpha_{1,7},\alpha_{1,8}\}$, $\{\alpha_{2,3},\alpha_{2,6},\alpha_{2,7},\alpha_{2,8}\}$,
$\{\alpha_{3,3},\alpha_{1,6}, \linebreak \alpha_{3,7},\alpha_{3,8}\}$.
Both $G$ and $G_1$ are marked by $\ast$ in Tables 1---4.

\medskip

Case 38:
({\bf n}=10,\
$
\left(\begin{array}{cccc}
4\aaa_1 & (8\aaa_1)_{I}     & (8\aaa_1)_I    & 4\aaa_1\amalg 2\aaa_2   \\
               & 4\aaa_1     & 4\aaa_2       &  4\aaa_1\amalg 2\aaa_2  \\
               &             & 4\aaa_1       &  4\aaa_1\amalg 2\aaa_2  \\
       &                     &               &  2\aaa_2
\end{array}\right)
\subset 4\aaa_1\amalg 6\aaa_2)$\

$\Longleftarrow ({\bf n}=34,\ (4\aaa_1,6\aaa_2)\subset 4\aaa_1\amalg 6\aaa_2)$.

Similar to Cases 1 and 22. By \cite{Nik11}, the $G\cong \SSS_4$ is marked by $N_{22}$ and
$G=H_{34,1}$ of Case 22 with orbits
$\{\alpha_{1,1}, \alpha_{1,3}, \alpha_{2,4}, \alpha_{1,8}\}$, $\{\alpha_{1,2},\alpha_{2,5},
\alpha_{1,5},\alpha_{1,11},\alpha_{2,2},\alpha_{2,11},
\alpha_{1,12},\alpha_{2,12},\alpha_{2,7},\alpha_{1,7},\alpha_{1,9}, \linebreak
\alpha_{2,9}\}$. The $G_1\cong D_8$ is marked by
$G_1\subset G=H_{34,1}$ of Case 22
with suborbits
$\{\alpha_{1,1}, \alpha_{1,3}, \alpha_{2,4}, \linebreak \alpha_{1,8}\}$,
$\{\alpha_{1,2},\alpha_{1,11},\alpha_{2,7},\alpha_{2,12}\}$, $\{\alpha_{2,2},\alpha_{2,11},\alpha_{1,7},\alpha_{1,12}\}$,
$\{\alpha_{1,5},\alpha_{1,9},\alpha_{2,9},\alpha_{2,5}\}$.
Both $G$ and $G_1$ are marked by $\ast$ in Tables 1---4.

\medskip

Case 39:
$({\bf n}=12,\ (8\aaa_1,8\aaa_1)\subset 16\aaa_1)\Longleftarrow ({\bf n}=75,\ 16\aaa_1)$.

Similar to Cases 1 and 10. By \cite{Nik10}, the $G\cong 4^2\AAA_4$ is marked
by $N_{23}$ and $G=H_{75,1}$ from Case 10 with the orbit
$
\{\alpha_{1},\alpha_{9},\alpha_{6},\alpha_{7},\alpha_{4},
\alpha_{21},\alpha_{20},
\alpha_{22},\alpha_{5},\alpha_{14},\alpha_{10},\alpha_{23},
\alpha_{3},\alpha_{24},\alpha_{18},\alpha_{16}\}$.
The $G_1\cong Q_8$ is marked by
$$
G=H_{75,1}\supset G_1=
$$
%gap> H1:=Image(F[1]);
%Group([ (1,4,9,5)(3,20,24,14)(6,18,21,7)(10,16,22,23)(11,12)(13,17), (1,3,9,24)(4,14,5,20)(6,22,21,10)(7,23,18,16)
%(11,17)(12,13) ])
%gap> Orbits(H1);
%[ [ 1, 4, 3, 9, 14, 20, 5, 24 ], [ 6, 18, 22, 21, 16, 23, 7, 10 ], [ 11, 12, 17, 13 ] ]
$$
[(\alpha_{1}\alpha_{4}\alpha_{9}\alpha_{5})(\alpha_{3}\alpha_{20}\alpha_{24}\alpha_{14})
(\alpha_{6}\alpha_{18}\alpha_{21}\alpha_{7})(\alpha_{10}\alpha_{16}\alpha_{22}\alpha_{23})
(\alpha_{11}\alpha_{12})(\alpha_{13}\alpha_{17}),
$$
$$
(\alpha_{1}\alpha_{3}\alpha_{9}\alpha_{24})(\alpha_{4}\alpha_{14}\alpha_{5}\alpha_{20})
(\alpha_{6}\alpha_{22}\alpha_{21}\alpha_{10})(\alpha_{7}\alpha_{23}\alpha_{18}\alpha_{16})
(\alpha_{11}\alpha_{17})(\alpha_{12}\alpha_{13})]
$$
with suborbits
$\{\alpha_{1},\alpha_{4},\alpha_{3},\alpha_{9},\alpha_{14},\alpha_{20},\alpha_{5},\alpha_{24}\}$,
$\{\alpha_{6},\alpha_{18},\alpha_{22},\alpha_{21},\alpha_{16},\alpha_{23},\alpha_{7},\alpha_{10}\}$.
Both $G$ and $G_1$ are marked by $\ast$ in Tables 1---4.

\medskip

Case 40:
$({\bf n}=12,\ (\aaa_2,\aaa_2)\subset 2\aaa_2)\Longleftarrow ({\bf n}=26,\ 2\aaa_2)$.

By \cite{Nik10}, the $G\cong SD_{16}$ is marked by $N_{22}$ and
$$
G=H_{26,1}=
$$
%G:=Group([ (3,4)(7,8)(9,12,22,18)(10,11,21,17)(13,20,16,23)(14,19,15,24),
%(5,8)(6,7)(11,24)(12,23)(13,16)(14,15)(17,19)(18,20) ]);
$$
[(\alpha_{1,2}\alpha_{2,2})(\alpha_{1,4}\alpha_{2,4})
(\alpha_{1,5}\alpha_{2,6}\alpha_{2,11}\alpha_{2,9})
(\alpha_{2,5}\alpha_{1,6}\alpha_{1,11}\alpha_{1,9})
(\alpha_{1,7}\alpha_{2,10}\alpha_{2,8}\alpha_{1,12})
(\alpha_{2,7}\alpha_{1,10}\alpha_{1,8}\alpha_{2,12}),
$$
$$
(\alpha_{1,3}\alpha_{2,4})(\alpha_{2,3}\alpha_{1,4})(\alpha_{1,6}\alpha_{2,12})(\alpha_{2,6}\alpha_{1,12})
(\alpha_{1,7}\alpha_{2,8})(\alpha_{2,7}\alpha_{1,8})(\alpha_{1,9}\alpha_{1,10})(\alpha_{2,9}\alpha_{2,10})]
$$
with the orbit
%%[ [ 3, 4 ], [ 5, 8, 7, 6 ], [ 9, 12, 22, 23, 18, 13, 20, 16 ],
%%[ 10, 11, 21, 24, 17, 14, 19, 15 ] ]
%$
%\{\alpha_{1,2},\alpha_{2,2}\},\
%\{\alpha_{1,3},\alpha_{2,4},\alpha_{1,4},\alpha_{2,3}\},\
%\{\alpha_{1,5},\alpha_{2,6},\alpha_{2,11},\alpha_{1,12},\
%\linebreak
%\alpha_{2,9},\alpha_{1,7},\alpha_{2,10},\alpha_{2,8}\},\
%\{\alpha_{2,5},\alpha_{1,6},\alpha_{1,11},\alpha_{2,12},
%\alpha_{1,9},\alpha_{2,7},\alpha_{1,10},\alpha_{1,8}\}$.
$\{\alpha_{1,3},\alpha_{2,4},\alpha_{1,4},\alpha_{2,3}\}$.
%gap> H1:=Image(F[1]);
%Group([ (5,6)(7,8)(9,13,22,16)(10,14,21,15)(11,24,17,19)(12,23,18,20), (3,4)(7,8)(9,12,22,18)(10,11,21,17)
%(13,20,16,23)(14,19,15,24) ])
%gap>
%gap> Orbits(H1);
%[ [ 3, 4 ], [ 5, 6 ], [ 7, 8 ], [ 9, 13, 12, 22, 20, 23, 16, 18 ], [ 10, 14, 11, 21, 19, 24, 15, 17 ] ]
The $G_1\cong Q_8$ is marked by
$$
G=H_{26,1}\supset G_1
$$
$$
[(\alpha_{1,3}\alpha_{2,3})(\alpha_{1,4}\alpha_{2,4})(\alpha_{1,5}\alpha_{1,7}\alpha_{2,11}\alpha_{2,8})
(\alpha_{2,5}\alpha_{2,7}\alpha_{1,11}\alpha_{1,8})
(\alpha_{1,6}\alpha_{2,12}\alpha_{1,9}\alpha_{1,10})(\alpha_{2,6}\alpha_{1,12}\alpha_{2,9}\alpha_{2,10}),
$$
$$
(\alpha_{1,2}\alpha_{2,2})(\alpha_{1,4}\alpha_{2,4})(\alpha_{1,5}\alpha_{2,6}\alpha_{2,11}\alpha_{2,9})
(\alpha_{2,5}\alpha_{1,6}\alpha_{1,11}\alpha_{1,9})
(\alpha_{1,7}\alpha_{2,10}\alpha_{2,8}\alpha_{1,12})(\alpha_{2,7}\alpha_{1,10}\alpha_{1,8}\alpha_{2,12})]
$$
with suborbits
$\{\alpha_{1,3},\alpha_{2,3}\}$, $\{\alpha_{1,4},\alpha_{2,4}\}$.
Both $G$ and $G_1$ are marked by $\ast$ in Tables 1---4.

\medskip

Case 41:
$({\bf n}=16,\
(5\aaa_1,5\aaa_1,5\aaa_1)\subset 15\aaa_1)
\Longleftarrow ({\bf n}=55,\ 15\aaa_1)$.

Similar to Case 1. By \cite{Nik10}, the $G\cong \AAA_5$ is marked
by $N_{23}$ and
%%G:=Group([(2,8,7,17,11)(3,14,10,13,12)(4,19,16,15,9)(6,23,20,18,22),
%%(2,13)(3,22)(4,9)(7,18)(8,14)(11,20)(12,23)(16,19)]);
$$
G=H_{55,2}=
$$
$$[(\alpha_{2}\alpha_{8}\alpha_{7}\alpha_{17}\alpha_{11})
(\alpha_{3}\alpha_{14}\alpha_{10}\alpha_{13}\alpha_{12})
(\alpha_{4}\alpha_{19}\alpha_{16}\alpha_{15}\alpha_{9})
(\alpha_{6}\alpha_{23}\alpha_{20}\alpha_{18}\alpha_{22}),
$$
$$
(\alpha_{2}\alpha_{13})(\alpha_{3}\alpha_{22})
(\alpha_{4}\alpha_{9})(\alpha_{7}\alpha_{18})
(\alpha_{8}\alpha_{14})(\alpha_{11}\alpha_{20})
(\alpha_{12}\alpha_{23})(\alpha_{16}\alpha_{19})]
$$
with the orbit
%$
%\{\alpha_{2},\alpha_{8},\alpha_{13},\alpha_{7},\alpha_{14},
%\alpha_{12},\alpha_{17},\alpha_{18},\alpha_{10},\alpha_{3},
%\alpha_{23},\alpha_{11},\alpha_{22},\alpha_{20},
%\linebreak
%\alpha_{6}\},\
%\{\alpha_{4},\alpha_{19},\alpha_{9},\alpha_{16},\alpha_{15}\}$.
$\{\alpha_{2},\alpha_{8},\alpha_{13},\alpha_{7},\alpha_{14},
\alpha_{12},\alpha_{17},\alpha_{18},\alpha_{10},\alpha_{3},
\alpha_{23},\alpha_{11},\alpha_{22},\alpha_{20},\alpha_{6}\}$.
%gap> H1:=Image(F[1]);
%Group([ (2,6)(4,15)(7,23)(9,16)(10,20)(11,14)(12,17)(13,22), (3,7)(6,11)(8,22)(9,15)(10,23)(13,17)(14,18)(16,19) ])
%gap> Orbits(H1);
%[ [ 2, 6, 11, 14, 18 ], [ 3, 7, 23, 10, 20 ], [ 4, 15, 9, 16, 19 ], [ 8, 22, 13, 17, 12 ] ]
The $G_1\cong D_{10}$ is marked by
$$
G=H_{55,2}\supset G_1=
$$
$$
[(\alpha_{2}\alpha_{6})(\alpha_{4}\alpha_{15})(\alpha_{7}\alpha_{23})(\alpha_{9}\alpha_{16})
(\alpha_{10}\alpha_{20})(\alpha_{11}\alpha_{14})(\alpha_{12}\alpha_{17})(\alpha_{13}\alpha_{22}),
$$
$$
(\alpha_{3}\alpha_{7})(\alpha_{6}\alpha_{11})(\alpha_{8}\alpha_{22})(\alpha_{9}\alpha_{15})
(\alpha_{10}\alpha_{23})(\alpha_{13}\alpha_{17})(\alpha_{14}\alpha_{18})(\alpha_{16}\alpha_{19})]
$$
with suborbits
$\{\alpha_{2},\alpha_{6},\alpha_{11},\alpha_{14},\alpha_{18}\}$,
$\{\alpha_{3},\alpha_{7},\alpha_{23},\alpha_{10},\alpha_{20}\}$,
$\{\alpha_{8},\alpha_{22},\alpha_{13},\alpha_{17},\alpha_{12}\}$.
Both $G$ and $G_1$ are marked by $\ast$ in Tables 1---4.

\medskip

Case 42:
$({\bf n}=17,\ (\aaa_1,\aaa_1)
\subset 2\aaa_1)\Longleftarrow ({\bf n}=34,\ 2\aaa_1)$.

Similar to Cases 1 and 29. By \cite{Nik10}, the $G\cong \SSS_4$ is marked
by $N_{23}$ and
%%G:=Group([(1,16,19,15)(3,5,12,14)(4,9)(7,20,23,22)(8,11,17,18)(10,13),
%%(2,13)(3,22)(4,9)(7,18)(8,14)(11,20)(12,23)(16,19)]);
$G=H_{34,2}$ of Case 29 with the orbit
$\{\alpha_{4},\alpha_{9}\}$.
%gap> H1:=Image(F[1]);
%Group([ (2,10,13)(3,20,14)(5,18,7)(8,11,22)(12,23,17)(15,16,19), (1,15)(3,18)(5,17)(7,22)(8,14)(11,12)(16,19)
%(20,23) ])
%gap> Orbits(H1);
%[ [ 1, 15, 16, 19 ], [ 2, 10, 13 ], [ 3, 20, 18, 14, 23, 7, 8, 17, 5, 22, 11, 12 ] ]
The $G_1\cong \AAA_{4}$ is marked by
$$
G=H_{34,2}\supset G_1=
$$
$$
[(\alpha_{2}\alpha_{10}\alpha_{13})(\alpha_{3}\alpha_{20}\alpha_{14})(\alpha_{5}\alpha_{18}\alpha_{7})
(\alpha_{8}\alpha_{11}\alpha_{22})(\alpha_{12}\alpha_{23}\alpha_{17})(\alpha_{15}\alpha_{16}\alpha_{19}),
$$
$$
(\alpha_{1}\alpha_{15})(\alpha_{3}\alpha_{18})(\alpha_{5}\alpha_{17})(\alpha_{7}\alpha_{22})
(\alpha_{8}\alpha_{14})(\alpha_{11}\alpha_{12})(\alpha_{16}\alpha_{19})(\alpha_{20}\alpha_{23})]
$$
with suborbits $\{\alpha_{4}\}$, $\{\alpha_{9}\}$.
Both $G$ and $G_1$ are marked by $\ast$ in Tables 1---4.

\medskip

Case 43: $({\bf n}=17,\ (\aaa_1,3\aaa_1)
\subset 4\aaa_1)\Longleftarrow ({\bf n}=49,\ 4\aaa_1)$.

Similar to Cases 1 and 5. By \cite{Nik10},
the $G\cong 2^4C_3$ is marked
by $N_{23}$ and $G=H_{49,1}$ from Case 5
%%G1:=Group([(1,22,19)(3,16,17)(4,20,9)(7,10,8)(12,13,23)(14,18,21),
%(2,12)(3,8)(4,20)(7,16)(9,11)(13,23)(14,22)(18,19),
%(2,13)(3,22)(4,9)(7,18)(8,14)(11,20)(12,23)(16,19)]);
%gap> Orbits(G);
%%[ [ 1, 22, 19, 14, 3, 18, 16, 8, 21, 7, 17, 10 ], [ 2, 12, 13, 23 ], [ 4, 20, 9, 11 ] ]
with the orbit $\{\alpha_{2},\alpha_{12},\alpha_{13},\alpha_{23}\}$.
%gap> H4:=Image(F[4]);
%Group([ (1,3,18)(7,17,22)(8,19,21)(9,20,11)(10,14,16)(12,13,23), (1,10)(3,8)(4,11)(7,18)(9,20)(14,22)(16,19)(17,21) ])
%gap> Orbits(H4);
%[ [ 1, 3, 10, 18, 8, 14, 7, 19, 16, 22, 17, 21 ], [ 4, 11, 9, 20 ], [ 12, 13, 23 ] ]
The $G_1\cong \AAA_{4}$ is marked by
$$
G=H_{49,1}\supset G_1=
$$
$$
[(\alpha_{1}\alpha_{3}\alpha_{18})(\alpha_{7}\alpha_{17}\alpha_{22})(\alpha_{8}\alpha_{19}\alpha_{21})
(\alpha_{9}\alpha_{20}\alpha_{11})(\alpha_{10}\alpha_{14}\alpha_{16})(\alpha_{12}\alpha_{13}\alpha_{23}),
$$
$$
(\alpha_{1}\alpha_{10})(\alpha_{3}\alpha_{8})(\alpha_{4}\alpha_{11})(\alpha_{7}\alpha_{18})
(\alpha_{9}\alpha_{20})(\alpha_{14}\alpha_{22})(\alpha_{16}\alpha_{19})(\alpha_{17}\alpha_{21})]
$$
with suborbits
$\{\alpha_{2}\}$, $\{\alpha_{12},\alpha_{13},\alpha_{23}\}$.
Both $G$ and $G_1$ are marked by $\ast$ in Tables 1---4.

\medskip

Case 44:
$({\bf n}=17,\ (4\aaa_1,4\aaa_1)
\subset 8\aaa_1)\Longleftarrow ({\bf n}=34,\ 8\aaa_1)$.

Similar to Case 1. By \cite{Nik10}, the $G\cong \SSS_4$
is marked by $N_{23}$ and
%%G1:=Group([(1,18,15,22)(2,4,20,11)(3,7)(8,10,16,24)(9,12,23,13)(14,19),
%%(2,13)(3,22)(4,9)(7,18)(8,14)(11,20)(12,23)(16,19)]);
$$
G=H_{34,1}=
$$
$$
[(\alpha_{1}\alpha_{18}\alpha_{15}\alpha_{22})
(\alpha_{2}\alpha_{4}\alpha_{20}\alpha_{11})
(\alpha_{3}\alpha_{7})
(\alpha_{8}\alpha_{10}\alpha_{16}\alpha_{24})
(\alpha_{9}\alpha_{12}\alpha_{23}\alpha_{13})
(\alpha_{14}\alpha_{19}),
$$
$$
(\alpha_{2}\alpha_{13})(\alpha_{3}\alpha_{22})
(\alpha_{4}\alpha_{9})(\alpha_{7}\alpha_{18})
(\alpha_{8}\alpha_{14})(\alpha_{11}\alpha_{20})
(\alpha_{12}\alpha_{23})(\alpha_{16}\alpha_{19})]
$$
with the orbit
%$
%\{\alpha_{1},\alpha_{3},\alpha_{22},\alpha_{15},\alpha_{7},\alpha_{18}\},\
%\{\alpha_{2},\alpha_{13},\alpha_{20},\alpha_{9},
%\alpha_{12},\alpha_{4},\alpha_{23},
%\linebreak
%\alpha_{11}\},\
%\{\alpha_{8},\alpha_{14},\alpha_{16},\alpha_{10},
%\alpha_{19},\alpha_{24}\}$;
$\{\alpha_{2},\alpha_{13},\alpha_{20},\alpha_{9},\alpha_{12},\alpha_{4},\alpha_{23},\alpha_{11}\}$.
%gap> H1:=Image(F[1]);
%Group([ (1,3,22)(4,12,11)(7,18,15)(8,10,14)(9,20,23)(16,24,19), (2,9)(3,7)(4,13)(8,16)(11,12)(14,19)(18,22)(20,23) ])
%gap> Orbits(H1);
%[ [ 1, 3, 22, 7, 18, 15 ], [ 2, 9, 20, 23 ], [ 4, 12, 13, 11 ], [ 8, 10, 16, 14, 24, 19 ] ]
The $G_1\cong \AAA_{4}$ is marked by
$$
G=H_{34,1}\supset G_1=
$$
$$
[(\alpha_{1}\alpha_{3}\alpha_{22})(\alpha_{4}\alpha_{12}\alpha_{11})(\alpha_{7}\alpha_{18}\alpha_{15})
(\alpha_{8}\alpha_{10}\alpha_{14})(\alpha_{9}\alpha_{20}\alpha_{23})(\alpha_{16}\alpha_{24}\alpha_{19}),
$$
$$
(\alpha_{2}\alpha_{9})(\alpha_{3}\alpha_{7})(\alpha_{4}\alpha_{13})(\alpha_{8}\alpha_{16})
(\alpha_{11}\alpha_{12})(\alpha_{14}\alpha_{19})(\alpha_{18}\alpha_{22})(\alpha_{20}\alpha_{23})]
$$
with suborbits
$\{\alpha_{2},\alpha_{9},\alpha_{20},\alpha_{23}\}$, $\{\alpha_{4},\alpha_{12},\alpha_{13},\alpha_{11}\}$.
Both $G$ and $G_1$ are marked by $\ast$ in Tables 1---4.

\medskip

Case 45:
$({\bf n}=17,\ (4\aaa_1,12\aaa_1)
\subset 16\aaa_1)\Longleftarrow ({\bf n}=49,\ 16\aaa_1)$.

Similar to Case 1. By \cite{Nik10}, the $G\cong 2^4C_3$
is marked by $N_{23}$ and
$$
G=H_{49,2}=
$$
%%G2:=Group([(1,21,6)(2,14,4)(3,8,7)(9,22,23)(11,19,20)(12,18,13),
%%(2,12)(3,8)(4,20)(7,16)(9,11)(13,23)(14,22)(18,19),
%%(2,13)(3,22)(4,9)(7,18)(8,14)(11,20)(12,23)(16,19)]);
$$
[(\alpha_{1}\alpha_{21}\alpha_{6})
(\alpha_{2}\alpha_{14}\alpha_{4})
(\alpha_{3}\alpha_{8}\alpha_{7})
(\alpha_{9}\alpha_{22}\alpha_{23})
(\alpha_{11}\alpha_{19}\alpha_{20})
(\alpha_{12}\alpha_{18}\alpha_{13}),
$$
$$
(\alpha_{2}\alpha_{12})(\alpha_{3}\alpha_{8})
(\alpha_{4}\alpha_{20})(\alpha_{7}\alpha_{16})
(\alpha_{9}\alpha_{11})(\alpha_{13}\alpha_{23})
(\alpha_{14}\alpha_{22})(\alpha_{18}\alpha_{19}),
$$
$$
(\alpha_{2}\alpha_{13})(\alpha_{3}\alpha_{22})
(\alpha_{4}\alpha_{9})(\alpha_{7}\alpha_{18})
(\alpha_{8}\alpha_{14})(\alpha_{11}\alpha_{20})
(\alpha_{12}\alpha_{23})(\alpha_{16}\alpha_{19})]
$$
with the orbit
%$\{\alpha_{1},\alpha_{21},\alpha_{6}\},\
%\{\alpha_{2},\alpha_{14},\alpha_{12},\alpha_{11},\alpha_{23},
%\alpha_{19},\alpha_{4},\alpha_{22},\alpha_{18}, \alpha_{3},
%\alpha_{9},
%\alpha_{13},\alpha_{7},
%\alpha_{20},\alpha_{8},\alpha_{16}\}$.
$\{\alpha_{2},\alpha_{14},\alpha_{12},\alpha_{11},\alpha_{23},
\alpha_{19},\alpha_{4},\alpha_{22},\alpha_{18}, \alpha_{3},
\alpha_{9},
\alpha_{13},\alpha_{7},
\alpha_{20},\alpha_{8},\alpha_{16}\}$.
%gap> H1:=Image(F[1]);
%Group([ (1,6,21)(3,22,13)(4,16,14)(7,19,23)(8,18,20)(9,11,12), (2,3)(4,18)(7,9)(8,12)(11,16)(13,22)(14,23)(19,20) ])
%gap> Orbits(H1);
%[ [ 1, 6, 21 ], [ 2, 3, 22, 13 ], [ 4, 16, 18, 14, 11, 20, 23, 12, 8, 19, 7, 9 ] ]
The $G_1\cong \AAA_{4}$ is marked by
$$
G=H_{49,2}\supset G_1=
$$
$$
[(\alpha_{1}\alpha_{6}\alpha_{21})(\alpha_{3}\alpha_{22}\alpha_{13})(\alpha_{4}\alpha_{16}\alpha_{14})
(\alpha_{7}\alpha_{19}\alpha_{23})(\alpha_{8}\alpha_{18}\alpha_{20})(\alpha_{9}\alpha_{11}\alpha_{12}),
$$
$$
(\alpha_{2}\alpha_{3})(\alpha_{4}\alpha_{18})(\alpha_{7}\alpha_{9})(\alpha_{8}\alpha_{12})
(\alpha_{11}\alpha_{16})(\alpha_{13}\alpha_{22})(\alpha_{14}\alpha_{23})(\alpha_{19}\alpha_{20})]
$$
with suborbits
$\{\alpha_{2},\alpha_{3},\alpha_{22},\alpha_{13}\}$,
$\{\alpha_{4},\alpha_{16},\alpha_{18},\alpha_{14},\alpha_{11},\alpha_{20},\alpha_{23},
\alpha_{12},\alpha_{8},\alpha_{19},\alpha_{7},\alpha_{9}\}$.
Both $G$ and $G_1$ are marked by $\ast$ in Tables 1---4.

\medskip

Case 46:
$({\bf n}=17,\ (6\aaa_1,6\aaa_1)
\subset 12\aaa_1)\Longleftarrow ({\bf n}=49,\ 12\aaa_1)$.

Similar to Cases 1 and 5. By \cite{Nik10}, the $G\cong 2^4C_3$
is marked by $N_{23}$ and
$G=H_{49,1}$ from Case 5 with the orbit
$\{\alpha_{1},\alpha_{22},\alpha_{21},\alpha_{10},\alpha_{19},\alpha_{14},
\alpha_{8},\alpha_{17},\alpha_{18},\alpha_{7},\alpha_{3},\alpha_{16}\}$.
The $G_1\cong \AAA_4$ is marked by
$$
G=H_{49,1}\supset G_1=
$$
$$
[(\alpha_{1}\alpha_{3}\alpha_{18})(\alpha_{7}\alpha_{17}\alpha_{22})(\alpha_{8}\alpha_{19}\alpha_{21})
(\alpha_{9}\alpha_{20}\alpha_{11})(\alpha_{10}\alpha_{14}\alpha_{16})(\alpha_{12}\alpha_{13}\alpha_{23}),
$$
$$
(\alpha_{2}\alpha_{12})(\alpha_{3}\alpha_{8})(\alpha_{4}\alpha_{20})(\alpha_{7}\alpha_{16})
(\alpha_{9}\alpha_{11})(\alpha_{13}\alpha_{23})(\alpha_{14}\alpha_{22})(\alpha_{18}\alpha_{19})]
$$
with suborbits
$\{\alpha_{1},\alpha_{3},\alpha_{18},\alpha_{8},\alpha_{19},\alpha_{21}\}$,
$\{\alpha_{7},\alpha_{17},\alpha_{16},\alpha_{22},\alpha_{10},\alpha_{14}\}$.
%gap> H1:=Image(F[1]);
%Group([ (1,3,18)(7,17,22)(8,19,21)(9,20,11)(10,14,16)(12,13,23), (2,12)(3,8)(4,20)(7,16)(9,11)(13,23)(14,22)(18,19) ])
%gap>
%gap> Orbits(H1);
%[ [ 1, 3, 18, 8, 19, 21 ], [ 2, 12, 13, 23 ], [ 4, 20, 11, 9 ], [ 7, 17, 16, 22, 10, 14 ] ]
Both $G$ and $G_1$ are marked by $\ast$ in Tables 1---4.

\medskip

Case 47:
$({\bf n}=17,\ (6\aaa_1,6\aaa_1)
\subset 6\aaa_2)\Longleftarrow ({\bf n}=34,\ 6\aaa_2)$.

Similar to Cases 1 and 22. By \cite{Nik10},
the $G\cong \SSS_4$ is marked by $N_{22}$ and
$G=H_{34,1}$ of Case 22
%%gap> h34c1;
%%G1:=Group([ (3,10)(4,9)(7,16)(8,15)(13,18)(14,17)(21,23)(22,24),
%%(1,5,8)(2,6,7)(3,9,22)(4,10,21)(13,17,24)(14,18,23) ]);
with the orbit
%%Orbits(h34c1);
%%[ [ 1, 5, 8, 15 ], [ 2, 6, 7, 16 ],
%%[ 3, 10, 9, 21, 4, 22, 23, 24, 14, 13, 17, 18 ] ]
$\{\alpha_{1,2},\alpha_{2,5},\alpha_{1,5},\alpha_{1,11},\alpha_{2,2},
\alpha_{2,11},\alpha_{1,12},\alpha_{2,12},\alpha_{2,7},\alpha_{1,7},\alpha_{1,9},\alpha_{2,9}\}$.
%%gap> H1:=Image(F[1]);
%%Group([ (3,18,23)(4,17,24)(5,15,8)(6,16,7)(9,22,14)(10,21,13), (1,5)(2,6)(3,14)(4,13)(7,16)(8,15)(9,18)(10,17) ])
%%gap> Orbits(H1);
%%[ [ 1, 5, 15, 8 ], [ 2, 6, 16, 7 ], [ 3, 18, 14, 23, 9, 22 ], [ 4, 17, 13, 24, 10, 21 ] ]
The $G_1\cong \AAA_4$ is marked by
$$
G=H_{34,1}\supset G_1=
$$
$$
[(\alpha_{1,2}\alpha_{2,9}\alpha_{1,12})(\alpha_{2,2}\alpha_{1,9}\alpha_{2,12})(\alpha_{1,3}\alpha_{1,8}\alpha_{2,4})
(\alpha_{2,3}\alpha_{2,8}\alpha_{1,4})(\alpha_{1,5}\alpha_{2,11}\alpha_{2,7})(\alpha_{2,5}\alpha_{1,11}\alpha_{1,7}),
$$
$$
(\alpha_{1,1}\alpha_{1,3})(\alpha_{2,1}\alpha_{2,3})(\alpha_{1,2}\alpha_{2,7})(\alpha_{2,2}\alpha_{1,7})
(\alpha_{1,4}\alpha_{2,8})(\alpha_{2,4}\alpha_{1,8})(\alpha_{1,5}\alpha_{2,9})(\alpha_{2,5}\alpha_{1,9})]
$$
with suborbits
$\{\alpha_{1,2},\alpha_{2,9},\alpha_{2,7},\alpha_{1,12},\alpha_{1,5},\alpha_{2,11}\}$,
$\{\alpha_{2,2},\alpha_{1,9},\alpha_{1,7},\alpha_{2,12},\alpha_{2,5},\alpha_{1,11}\}$.
Both $G$ and $G_1$ are marked by $\ast$ in Tables 1---4.

\medskip

Case 48:
$({\bf n}=17,\ (\aaa_1,\aaa_1,\aaa_1)
\subset 3\aaa_1)\Longleftarrow ({\bf n}=61,\ 3\aaa_1)$.

Similar to Case 1. By \cite{Nik10}, the $G\cong \AAA_{4,3}$
is marked by $N_{23}$ and
$$
G=H_{61,1}=
$$
%G1:=Group([(1,10,23)(2,4,11)(3,21,20)(8,9,17)(13,22,14)(16,19,24),
%(2,14)(3,17)(4,8)(7,15)(9,13)(10,12)(11,21)(16,19),
%(2,13)(3,22)(4,9)(7,18)(8,14)(11,20)(12,23)(16,19)]);
$$
[(\alpha_{1}\alpha_{10}\alpha_{23})(\alpha_{2}\alpha_{4}\alpha_{11})
(\alpha_{3}\alpha_{21}\alpha_{20})(\alpha_{8}\alpha_{9}\alpha_{17})
(\alpha_{13}\alpha_{22}\alpha_{14})(\alpha_{16}\alpha_{19}\alpha_{24}),
$$
$$
(\alpha_{2}\alpha_{14})(\alpha_{3}\alpha_{17})(\alpha_{4}\alpha_{8})(\alpha_{7}\alpha_{15})
(\alpha_{9}\alpha_{13})(\alpha_{10}\alpha_{12})(\alpha_{11}\alpha_{21})(\alpha_{16}\alpha_{19}),
$$
$$
(\alpha_{2}\alpha_{13})(\alpha_{3}\alpha_{22})(\alpha_{4}\alpha_{9})(\alpha_{7}\alpha_{18})
(\alpha_{8}\alpha_{14})(\alpha_{11}\alpha_{20})(\alpha_{12}\alpha_{23})(\alpha_{16}\alpha_{19})]
$$
with the orbit
%$
%\{\alpha_{1},\alpha_{10},\alpha_{23},\alpha_{12}\},\
%\{\alpha_{2},\alpha_{4},\alpha_{21},\alpha_{13},\alpha_{11}\alpha_{8},
%\alpha_{20},\alpha_{17},\alpha_{22}\alpha_{14},\linebreak
%\alpha_{9},\alpha_{3}\},\
%\{\alpha_{7},\alpha_{18},\alpha_{15}\},\
%\{\alpha_{16},\alpha_{19},\alpha_{24}\}.$
$\{\alpha_{7},\alpha_{18},\alpha_{15}\}$.
%gap> H3:=Image(F[3]);
%Group([ (2,14,20)(3,22,17)(4,21,9)(8,13,11)(10,23,12)(16,24,19), (1,10)(2,13)(3,14)(4,20)(8,22)(9,11)(12,23)(17,21) ])
%gap> Orbits(H3);
%[ [ 1, 10, 23, 12 ], [ 2, 14, 13, 20, 3, 11, 4, 22, 8, 9, 21, 17 ], [ 16, 24, 19 ] ]
The $G_1\cong \AAA_4$ is marked by
$$
G=H_{61,1}\supset G_1=
$$
$$
[(\alpha_{2}\alpha_{14}\alpha_{20})(\alpha_{3}\alpha_{22}\alpha_{17})(\alpha_{4}\alpha_{21}\alpha_{9})
(\alpha_{8}\alpha_{13}\alpha_{11})(\alpha_{10}\alpha_{23}\alpha_{12})(\alpha_{16}\alpha_{24}\alpha_{19}),
$$
$$
(\alpha_{1}\alpha_{10})(\alpha_{2}\alpha_{13})(\alpha_{3}\alpha_{14})(\alpha_{4}\alpha_{20})
(\alpha_{8}\alpha_{22})(\alpha_{9}\alpha_{11})(\alpha_{12}\alpha_{23})(\alpha_{17}\alpha_{21})]
$$
with suborbits
$\{\alpha_{7}\}$, $\{\alpha_{18}\}$, $\{\alpha_{15}\}$.
Both $G$ and $G_1$ are marked by $\ast$ in Tables 1---4.

\medskip

Case 49:
$({\bf n}=17,\ (\aaa_1,\aaa_1,4\aaa_1)
\subset 6\aaa_1)\Longleftarrow ({\bf n}=34,\ (2\aaa_1,4\aaa_1)\subset 6\aaa_1)$.

Similar to Cases 1, 29 and 42. By \cite{Nik11}, the $G\cong \SSS_4$
is marked by $N_{23}$ and
%%G:=Group([(1,16,19,15)(3,5,12,14)(4,9)(7,20,23,22)(8,11,17,18)(10,13),
%%(2,13)(3,22)(4,9)(7,18)(8,14)(11,20)(12,23)(16,19)]);
%$
%\{\alpha_{1},\alpha_{16},\alpha_{19},\alpha_{15}\},\
%\{\alpha_{2},\alpha_{10},\alpha_{13}\},\
%\{\alpha_{3},\alpha_{22},\alpha_{12},\alpha_{18},\alpha_{17},
%\linebreak
%\alpha_{20},\alpha_{7},\alpha_{11},\alpha_{23},\alpha_{8},
%\alpha_{5},\alpha_{14}\},\
%\{\alpha_{4},\alpha_{9}\}$;
$G=H_{34,2}$ of Case 29 with orbits
$\{\alpha_{4},\alpha_{9}\}$, $\{\alpha_{1},\alpha_{16},\alpha_{19},\alpha_{15}\}$.
%gap> H1:=Image(F[1]);
%Group([ (2,10,13)(3,20,14)(5,18,7)(8,11,22)(12,23,17)(15,16,19), (1,15)(3,18)(5,17)(7,22)(8,14)(11,12)(16,19)
%(20,23) ])
%gap> Orbits(H1);
%[ [ 1, 15, 16, 19 ], [ 2, 10, 13 ], [ 3, 20, 18, 14, 23, 7, 8, 17, 5, 22, 11, 12 ] ]
The $G_1\cong \AAA_{4}$ is marked by
$G_1\subset G=H_{34,2}$ of Case 42 with suborbits
$\{\alpha_{4}\}$, $\{\alpha_{9}\}$, $\{\alpha_{1},\alpha_{16},\alpha_{19},\alpha_{15}\}$.
Both $G$ and $G_1$ are marked by $\ast$ in Tables 1---4.

\medskip

Case 50:
$({\bf n}=17,\ (\aaa_1,\aaa_1,6\aaa_1)
\subset 8\aaa_1)\Longleftarrow ({\bf n}=34,\ (2\aaa_1,(6\aaa_1)_{II})\subset 8\aaa_1)$.

Similar to Case 1. By \cite{Nik11}, the $G\cong \SSS_4$
is marked by $N_{23}$ and
%%G3:=Group([(1,6,18,7)(3,24,20,17)(4,9)(8,13,14,19)(10,15,12,23)(11,22),
%%(2,13)(3,22)(4,9)(7,18)(8,14)(11,20)(12,23)(16,19)]);
$$
G=H_{34,3}=
$$
$$
[(\alpha_{1}\alpha_{6}\alpha_{18}\alpha_{7})
(\alpha_{3}\alpha_{24}\alpha_{20}\alpha_{17})
(\alpha_{4}\alpha_{9})
(\alpha_{8}\alpha_{13}\alpha_{14}\alpha_{19})
(\alpha_{10}\alpha_{15}\alpha_{12}\alpha_{23})
(\alpha_{11}\alpha_{22}),
$$
$$
(\alpha_{2}\alpha_{13})(\alpha_{3}\alpha_{22})
(\alpha_{4}\alpha_{9})(\alpha_{7}\alpha_{18})
(\alpha_{8}\alpha_{14})(\alpha_{11}\alpha_{20})
(\alpha_{12}\alpha_{23})(\alpha_{16}\alpha_{19})]
$$
with orbits
%$
%\{\alpha_{1},\alpha_{18},\alpha_{6},\alpha_{7}\},\
%\{\alpha_{2},\alpha_{8},\alpha_{14},\alpha_{16},\alpha_{13},\alpha_{19}\},\
%\{\alpha_{3},\alpha_{22},\alpha_{20},
%\linebreak
%\alpha_{24},\alpha_{11},\alpha_{17}\},\
%\{\alpha_{4},\alpha_{9}\},\
%\{\alpha_{10},\alpha_{12},\alpha_{15},\alpha_{23}\}$;
$\{\alpha_{4},\alpha_{9}\}$,
$\{\alpha_{2},\alpha_{8},\alpha_{14},\alpha_{16},\alpha_{13},\alpha_{19}\}$.
%gap> H1:=Image(F[1]);
%Group([ (2,14,19)(3,22,24)(6,7,18)(8,13,16)(11,17,20)(12,15,23), (1,6)(2,16)(3,20)(7,18)(10,15)(11,22)(12,23)
%(13,19) ])
%gap> Orbits(H1);
%[ [ 1, 6, 7, 18 ], [ 2, 14, 16, 19, 8, 13 ], [ 3, 22, 20, 24, 11, 17 ], [ 10, 15, 23, 12 ] ]
The $G_1\cong \AAA_{4}$ is marked by
$$
G=H_{34,3}\supset G_1=
$$
$$
[(\alpha_{2}\alpha_{14}\alpha_{19})(\alpha_{3}\alpha_{22}\alpha_{24})(\alpha_{6}\alpha_{7}\alpha_{18})
(\alpha_{8}\alpha_{13}\alpha_{16})(\alpha_{11}\alpha_{17}\alpha_{20})(\alpha_{12}\alpha_{15}\alpha_{23}),
$$
$$
(\alpha_{1}\alpha_{6})(\alpha_{2}\alpha_{16})(\alpha_{3}\alpha_{20})(\alpha_{7}\alpha_{18})
(\alpha_{10}\alpha_{15})(\alpha_{11}\alpha_{22})(\alpha_{12}\alpha_{23})(\alpha_{13}\alpha_{19})]
$$
with suborbits
$\{\alpha_{4}\}$, $\{\alpha_{9}\}$,
$\{\alpha_{2},\alpha_{8},\alpha_{14},\alpha_{16},\alpha_{13},\alpha_{19}\}$.
Both $G$ and $G_1$ are marked by $\ast$ in Tables 1---4.

\medskip

Case 51:
$({\bf n}=17,\ (\aaa_1,\aaa_1,12\aaa_1)\subset 14\aaa_1)
\Longleftarrow ({\bf n}=34,\ (2\aaa_1,12\aaa_1)\subset 14\aaa_1)$.

Similar to Cases 1, 29 and 42. By \cite{Nik11}, the $G\cong \SSS_4$
is marked by $N_{23}$ and
%%G:=Group([(1,16,19,15)(3,5,12,14)(4,9)(7,20,23,22)(8,11,17,18)(10,13),
%%(2,13)(3,22)(4,9)(7,18)(8,14)(11,20)(12,23)(16,19)]);
%$
%\{\alpha_{1},\alpha_{16},\alpha_{19},\alpha_{15}\},\
%\{\alpha_{2},\alpha_{10},\alpha_{13}\},\
%\{\alpha_{3},\alpha_{22},\alpha_{12},\alpha_{18},\alpha_{17},
%\linebreak
%\alpha_{20},\alpha_{7},\alpha_{11},\alpha_{23},\alpha_{8},
%\alpha_{5},\alpha_{14}\},\
%\{\alpha_{4},\alpha_{9}\}$;
$G=H_{34,2}$ of Case 29 with orbits
$\{\alpha_{4},\alpha_{9}\}$,
$\{\alpha_{3},\alpha_{22},\alpha_{12},\alpha_{18},\alpha_{17},
\alpha_{20},\alpha_{7},\alpha_{11},\alpha_{23},\alpha_{8},\alpha_{5},\alpha_{14}\}$.
%gap> H1:=Image(F[1]);
%Group([ (2,10,13)(3,20,14)(5,18,7)(8,11,22)(12,23,17)(15,16,19), (1,15)(3,18)(5,17)(7,22)(8,14)(11,12)(16,19)
%(20,23) ])
%gap> Orbits(H1);
%[ [ 1, 15, 16, 19 ], [ 2, 10, 13 ], [ 3, 20, 18, 14, 23, 7, 8, 17, 5, 22, 11, 12 ] ]
The $G_1\cong \AAA_{4}$ is marked by
$G_1\subset G=H_{34,2}$ of Case 42 with suborbits
$\{\alpha_{4}\}$, $\{\alpha_{9}\}$,
$\{\alpha_{3},\alpha_{22},\alpha_{12},\alpha_{18},\alpha_{17}, \linebreak
\alpha_{20},\alpha_{7},
\alpha_{11},\alpha_{23},\alpha_{8},\alpha_{5},\alpha_{14}\}$.
Both $G$ and $G_1$ are marked by $\ast$ in Tables 1---4.

\medskip

Case 52:
$({\bf n}=17,\ (\aaa_1,3\aaa_1,4\aaa_1)
\subset 8\aaa_1)\Longleftarrow ({\bf n}=65,\ 8\aaa_1)$.

Similar to Case 1. By \cite{Nik10}, the $G\cong 2^4D_6$
is marked by $N_{23}$ and
%G:=Group([(1,21,23)(3,22,12)(4,20,13)(5,18,16)(8,11,10)(9,14,17),
%(1,11,23,8)(4,6,13,20)(5,16)(7,17,9,14)(10,21,19,15)(12,24)]);
$$
G=H_{65,4}=
$$
$$[(\alpha_{1}\alpha_{21}\alpha_{23})(\alpha_{3}\alpha_{22}\alpha_{12})
(\alpha_{4}\alpha_{20}\alpha_{13})
(\alpha_{5}\alpha_{18}\alpha_{16})
(\alpha_{8}\alpha_{11}\alpha_{10})(\alpha_{9}\alpha_{14}\alpha_{17}),
$$
$$
(\alpha_{1}\alpha_{11}\alpha_{23}\alpha_{8})
(\alpha_{4}\alpha_{6}\alpha_{13}\alpha_{20})
(\alpha_{5}\alpha_{16})
(\alpha_{7}\alpha_{17}\alpha_{9}\alpha_{14})
(\alpha_{10}\alpha_{21}\alpha_{19}\alpha_{15})
(\alpha_{12}\alpha_{24})]
$$
with the orbit
%$\{\alpha_{1},\alpha_{8},\alpha_{21},\alpha_{23},\alpha_{15},
%\alpha_{11},\alpha_{10},\alpha_{19}\},\
%\{\alpha_{3},\alpha_{22},
%\alpha_{24},\alpha_{12}\},\
%\{\alpha_{4},
%\linebreak
%\alpha_{13},
%\alpha_{20},\alpha_{6}\},\
%\{\alpha_{5},\alpha_{16},
%\alpha_{18}\},\ \{\alpha_{7},\alpha_{17},\alpha_{14},\alpha_{9},\}$.
$\{\alpha_{1},\alpha_{8},\alpha_{21},\alpha_{23},\alpha_{15},
\alpha_{11},\alpha_{10},\alpha_{19}\}$.
%gap> H1:=Image(F[1]);
%Group([ (4,13,20)(5,16,18)(7,17,9)(8,11,19)(12,24,22)(15,21,23), (3,12)(4,13)(6,20)(7,14)(8,10)(9,17)(11,19)(22,24) ])
%gap> Orbits(H1);
%[ [ 3, 12, 24, 22 ], [ 4, 13, 20, 6 ], [ 5, 16, 18 ], [ 7, 17, 14, 9 ], [ 8, 11, 10, 19 ], [ 15, 21, 23 ] ]
The $G_1\cong \AAA_{4}$ is marked by
$$
G=H_{65,4}\supset G_1=
$$
$$
[(\alpha_{4}\alpha_{13}\alpha_{20})(\alpha_{5}\alpha_{16}\alpha_{18})(\alpha_{7}\alpha_{17}\alpha_{9})
(\alpha_{8}\alpha_{11}\alpha_{19})(\alpha_{12}\alpha_{24}\alpha_{22})(\alpha_{15}\alpha_{21}\alpha_{23}),
$$
$$
(\alpha_{3}\alpha_{12})(\alpha_{4}\alpha_{13})(\alpha_{6}\alpha_{20})(\alpha_{7}\alpha_{14})
(\alpha_{8}\alpha_{10})(\alpha_{9}\alpha_{17})(\alpha_{11}\alpha_{19})(\alpha_{22}\alpha_{24})]
$$
with suborbits
$\{\alpha_{1}\}$, $\{\alpha_{15},\alpha_{21},\alpha_{23}\}$,
$\{\alpha_{8},\alpha_{11},\alpha_{10},\alpha_{19}\}$.
Both $G$ and $G_1$ are marked by $\ast$ in Tables 1---4.

\medskip

Case 53:
$({\bf n}=17,\ (\aaa_1,3\aaa_1,12\aaa_1)
\subset 16\aaa_1)\Longleftarrow ({\bf n}=75,\ 16\aaa_1)$.

Similar to Cases 1 and 10. By \cite{Nik10}, the $G\cong 4^2\AAA_4$
is marked by $N_{23}$ and
$G=H_{75,1}$ from Case 10 with the orbit
$\{\alpha_{1},\alpha_{9},\alpha_{6},\alpha_{7},\alpha_{4},
\alpha_{21},\alpha_{20},
\alpha_{22},\alpha_{5},\alpha_{14},\alpha_{10},\alpha_{23},
\alpha_{3},\alpha_{24},\alpha_{18},\alpha_{16}\}$.
The $G_1\cong \AAA_4$ is marked by
$$
G=H_{75,1}\supset G_1=
$$
%gap> H1:=Image(F[1]);
%Group([ (3,4,22)(5,20,23)(6,10,16)(7,9,18)(11,17,12)(14,24,21), (3,16)(4,5)(6,21)(10,20)(11,12)(13,17)(14,22)
%(23,24) ])
%gap> Orbits(H1);
%[ [ 3, 4, 16, 22, 5, 6, 14, 20, 10, 21, 24, 23 ], [ 7, 9, 18 ], [ 11, 17, 12, 13 ] ]
%gap> Orbits(G);
%[ [ 1, 9, 6, 7, 4, 21, 20, 22, 5, 14, 10, 23, 3, 24, 18, 16 ], [ 11, 12, 13, 17 ] ]
$$
[(\alpha_{3}\alpha_{4}\alpha_{22})(\alpha_{5}\alpha_{20}\alpha_{23})(\alpha_{6}\alpha_{10}\alpha_{16})
(\alpha_{7}\alpha_{9}\alpha_{18})(\alpha_{11}\alpha_{17}\alpha_{12})(\alpha_{14}\alpha_{24}\alpha_{21}),
$$
$$
(\alpha_{3}\alpha_{16})(\alpha_{4}\alpha_{5})(\alpha_{6}\alpha_{21})(\alpha_{10}\alpha_{20})
(\alpha_{11}\alpha_{12})(\alpha_{13}\alpha_{17})(\alpha_{14}\alpha_{22})(\alpha_{23}\alpha_{24})]
$$
with suborbits
$\{\alpha_{1}\}$, $\{\alpha_{7},\alpha_{9},\alpha_{18}\}$,
$\{\alpha_{3},\alpha_{4},\alpha_{16},\alpha_{22},\alpha_{5},\alpha_{6},
\alpha_{14},\alpha_{20},\alpha_{10},\alpha_{21},\alpha_{24},\alpha_{23}\}$.
Both $G$ and $G_1$ are marked by $\ast$ in Tables 1---4.

\medskip

Case 54:
$({\bf n}=17,\ (\aaa_1,4\aaa_1,4\aaa_1)
\subset 9\aaa_1)\Longleftarrow ({\bf n}=34,\ (\aaa_1,8\aaa_1)\subset 9\aaa_1)$.

This is similar to Cases 1 and 44. By \cite{Nik11}, the $G\cong \SSS_4$
is marked by $N_{23}$ and
%%G1:=Group([(1,18,15,22)(2,4,20,11)(3,7)(8,10,16,24)(9,12,23,13)(14,19),
%%(2,13)(3,22)(4,9)(7,18)(8,14)(11,20)(12,23)(16,19)]);
$G=H_{34,1}$ of Case 44
with orbits
%$
%\{\alpha_{1},\alpha_{3},\alpha_{22},\alpha_{15},\alpha_{7},\alpha_{18}\},\
%\{\alpha_{2},\alpha_{13},\alpha_{20},\alpha_{9},
%\alpha_{12},\alpha_{4},\alpha_{23},
%\linebreak
%\alpha_{11}\},\
%\{\alpha_{8},\alpha_{14},\alpha_{16},\alpha_{10},
%\alpha_{19},\alpha_{24}\}$;
$\{\alpha_5\}$,
$\{\alpha_{2},\alpha_{13},\alpha_{20},\alpha_{9},\alpha_{12},\alpha_{4},\alpha_{23},\alpha_{11}\}$.
%gap> H1:=Image(F[1]);
%Group([ (1,3,22)(4,12,11)(7,18,15)(8,10,14)(9,20,23)(16,24,19), (2,9)(3,7)(4,13)(8,16)(11,12)(14,19)(18,22)(20,23) ])
%gap> Orbits(H1);
%[ [ 1, 3, 22, 7, 18, 15 ], [ 2, 9, 20, 23 ], [ 4, 12, 13, 11 ], [ 8, 10, 16, 14, 24, 19 ] ]
The $G_1\cong \AAA_{4}$ is marked by
$G_1\subset G=H_{34,1}$ of Case 44 with suborbits $\{\alpha_5\}$,
$\{\alpha_{2},\alpha_{9},\alpha_{20},\alpha_{23}\}$, $\{\alpha_{4},\alpha_{12},\alpha_{13},\alpha_{11}\}$.
Both $G$ and $G_1$ are marked by $\ast$ in Tables 1---4.

\medskip

Case 55:
$({\bf n}=17,\ (\aaa_1,6\aaa_1,6\aaa_1)
\subset \aaa_1\amalg 6\aaa_2)
\Longleftarrow ({\bf n}=34,\ (\aaa_1,6\aaa_2)\subset \aaa_1\amalg 6\aaa_2)$.

Similar to Cases 1 and 22. By \cite{Nik11}, the $G\cong \SSS_4$
is marked by $N_{22}$ and $G=H_{34,1}$ of Case 22
with orbits
$\{\alpha_{1,6}\}$, $\{\alpha_{1,2},\alpha_{2,5},
\alpha_{1,5},\alpha_{1,11},\alpha_{2,2},\alpha_{2,11},
\alpha_{1,12},\alpha_{2,12},\alpha_{2,7},\alpha_{1,7},\alpha_{1,9},\alpha_{2,9}\}$.
The $G_1\cong \AAA_4$ is marked by
%gap> H1:=Image(F[1]);
%Group([ (3,18,23)(4,17,24)(5,15,8)(6,16,7)(9,22,14)(10,21,13), (1,5)(2,6)(3,14)(4,13)(7,16)(8,15)(9,18)(10,17) ])
%gap> Orbits(H1);
%[ [ 1, 5, 15, 8 ], [ 2, 6, 16, 7 ], [ 3, 18, 14, 23, 9, 22 ], [ 4, 17, 13, 24, 10, 21 ] ]
$$
G=H_{34,1}\supset G_1=
$$
$$
[(\alpha_{1,2}\alpha_{2,9}\alpha_{1,12})(\alpha_{4 }\alpha_{1,9}\alpha_{2,12})(\alpha_{1,3}\alpha_{1,8}\alpha_{2,4})
(\alpha_{2,3}\alpha_{2,8}\alpha_{1,4})(\alpha_{1,5}\alpha_{2,11}\alpha_{2,7})(\alpha_{2,5}\alpha_{1,11}\alpha_{1,7}),
$$
$$
(\alpha_{1,1}\alpha_{1,3})(\alpha_{2,1}\alpha_{2,3})(\alpha_{1,2}\alpha_{2,7})(\alpha_{2,2}\alpha_{1,7})
(\alpha_{1,4}\alpha_{2,8})(\alpha_{2,4}\alpha_{1,8})(\alpha_{1,5}\alpha_{2,6})(\alpha_{2,5}\alpha_{1,9})]
$$
with suborbits
$\{\alpha_{1,6}\}$, $\{\alpha_{1,2},\alpha_{2,9},\alpha_{2,7},\alpha_{1,12},\alpha_{1,5},\alpha_{2,11}\}$,
$\{\alpha_{2,2},\alpha_{1,9},\alpha_{1,7},\alpha_{2,12},\alpha_{2,5},\alpha_{1,11}\}$. \newline 
Both $G$ and $G_1$ are marked by $\ast$ in Tables 1---4.

\medskip

Case 56:
$({\bf n}=17,\  (3\aaa_1,4\aaa_1,4\aaa_1)\subset 11\aaa_1)
\Longleftarrow ({\bf n}=34,\ (3\aaa_1,8\aaa_1)\subset 11\aaa_1)$.

Similar to Case 1. By \cite{Nik11}, the $G\cong \SSS_4$
is marked by $N_{23}$ and
%%G4:=Group([(1,6,16,19)(4,14,23,9)(5,11,20,21)(7,8,12,18)(13,15)(17,22),
%%(2,13)(3,22)(4,9)(7,18)(8,14)(11,20)(12,23)(16,19)]);
$$
G=H_{34,4}=
$$
$$
[(\alpha_{1}\alpha_{6}\alpha_{16}\alpha_{19})
(\alpha_{4}\alpha_{14}\alpha_{23}\alpha_{9})
(\alpha_{5}\alpha_{11}\alpha_{20}\alpha_{21})
(\alpha_{7}\alpha_{8}\alpha_{12}\alpha_{18})
(\alpha_{13}\alpha_{15})(\alpha_{17}\alpha_{22}),
$$
$$
(\alpha_{2}\alpha_{13})(\alpha_{3}\alpha_{22})
(\alpha_{4}\alpha_{9})(\alpha_{7}\alpha_{18})
(\alpha_{8}\alpha_{14})(\alpha_{11}\alpha_{20})
(\alpha_{12}\alpha_{23})(\alpha_{16}\alpha_{19})]
$$
with orbits
%$
%\{\alpha_{1},\alpha_{16},\alpha_{6},\alpha_{19}\},\
%\{\alpha_{2},\alpha_{15},\alpha_{13}\},\
%\{\alpha_{3},\alpha_{17},\alpha_{22}\},\
%\{\alpha_{4},\alpha_{7},
%\linebreak
%\alpha_{18},\alpha_{23},
%\alpha_{9},\alpha_{12},\alpha_{8},\alpha_{14}\},\
%\{\alpha_{5},\alpha_{11},\alpha_{20},\alpha_{21}\}$.
$\{\alpha_{2},\alpha_{15},\alpha_{13}\}$,
$\{\alpha_{4},\alpha_{7},\alpha_{18},\alpha_{23},
\alpha_{9},\alpha_{12},\alpha_{8},\alpha_{14}\}$.
%gap> H1:=Image(F[1]);
%Group([ (2,13,15)(3,22,17)(4,23,18)(5,11,20)(6,16,19)(7,12,9), (1,6)(4,18)(5,21)(7,9)(8,23)(11,20)(12,14)(16,19) ])
%gap> Orbits(H1);
%[ [ 1, 6, 16, 19 ], [ 2, 13, 15 ], [ 3, 22, 17 ], [ 4, 23, 18, 8 ], [ 5, 11, 21, 20 ], [ 7, 12, 9, 14 ] ]
The $G_1\cong \AAA_{4}$ is marked by
$$
G=H_{34,4}\supset G_1=
$$
$$
[(\alpha_{2}\alpha_{13}\alpha_{15})(\alpha_{3}\alpha_{22}\alpha_{17})(\alpha_{4}\alpha_{23}\alpha_{18})
(\alpha_{5}\alpha_{11}\alpha_{20})(\alpha_{6}\alpha_{16}\alpha_{19})(\alpha_{7}\alpha_{12}\alpha_{9}),
$$
$$
(\alpha_{1}\alpha_{6})(\alpha_{4}\alpha_{18})(\alpha_{5}\alpha_{21})(\alpha_{7}\alpha_{9})
(\alpha_{8}\alpha_{23})(\alpha_{11}\alpha_{20})(\alpha_{12}\alpha_{14})(\alpha_{16}\alpha_{19})]
$$
with suborbits
$\{\alpha_{2},\alpha_{15},\alpha_{13}\}$,
$\{\alpha_{4},\alpha_{23},\alpha_{18},\alpha_{8}\}$,
$\{\alpha_{7},\alpha_{12},\alpha_{9},\alpha_{14}\}$.
Both $G$ and $G_1$ are marked by $\ast$ in Tables 1---4.

\medskip

Case 57:
$({\bf n}=17,\ (4\aaa_1,4\aaa_1,4\aaa_1)\subset 12\aaa_1)
\Longleftarrow ({\bf n}=61,\ 12\aaa_1)$.

Similar to Cases 1 and 48. By \cite{Nik10}, the $G$ is marked by by $N_{23}$ and
$G=H_{61,1}$ of Case 48 with the orbit
$\{\alpha_{2},\alpha_{4},\alpha_{21},\alpha_{13},\alpha_{11},\alpha_{8},
\alpha_{20},\alpha_{17},\alpha_{22},\alpha_{14},\alpha_{9},\alpha_{3}\}$.
%gap> H2:=Image(F[2]);
%Group([ (2,13,21)(4,20,8)(7,15,18)(9,14,11)(10,12,23)(16,24,19), (1,10)(2,13)(3,14)(4,20)(8,22)(9,11)(12,23)(17,21) ])
%gap> Orbits(H2);
%[ [ 1, 10, 12, 23 ], [ 2, 13, 21, 17 ], [ 3, 14, 11, 9 ], [ 4, 20, 8, 22 ], [ 7, 15, 18 ], [ 16, 24, 19 ] ]
The $G_1\cong \AAA_{4}$ is marked by
$$
G=H_{61,1}\supset G_1=
$$
$$
[(\alpha_{2}\alpha_{13}\alpha_{21})(\alpha_{4}\alpha_{20}\alpha_{8})(\alpha_{7}\alpha_{15}\alpha_{18})
(\alpha_{9}\alpha_{14}\alpha_{11})(\alpha_{10}\alpha_{12}\alpha_{23})(\alpha_{16}\alpha_{24}\alpha_{19}),
$$
$$
(\alpha_{1}\alpha_{10})(\alpha_{2}\alpha_{13})(\alpha_{3}\alpha_{14})(\alpha_{4}\alpha_{20})
(\alpha_{8}\alpha_{22})(\alpha_{9}\alpha_{11})(\alpha_{12}\alpha_{23})(\alpha_{17}\alpha_{21})]
$$
with suborbits
$\{\alpha_{2},\alpha_{13},\alpha_{21},\alpha_{17}\}$, $\{\alpha_{3},\alpha_{14},\alpha_{11},\alpha_{9}\}$,
$\{\alpha_{4},\alpha_{20},\alpha_{8},\alpha_{22}\}$.
Both $G$ and $G_1$ are marked by $\ast$ in Tables 1---4.

\medskip

Case 58:
({\bf n}=17,\
$
\left(\begin{array}{rrr}
4\aaa_1 & 4\aaa_2 & 8\aaa_1  \\
         & 4\aaa_1 & 4\aaa_2   \\
         &         & 4\aaa_1
\end{array}\right)
\subset 4\aaa_3)\
\Longleftarrow ({\bf n}=34,\ (4\aaa_1,8\aaa_1)\subset 4\aaa_3)$.

Similar to Case 1. By \cite{Nik11},
the $G\cong \SSS_4$ is marked
by $N_{21}$ and
%gap> Image(HH[4]);
%G:=Group([ (1,7,19)(2,8,20)(3,9,21)(10,22,16)(11,23,17)(12,24,18),
%(1,3)(7,21,9,19)(8,20)(10,18,15,24)(11,17,14,23)(12,16,13,22)]);
%gap> Orbits(Image(HH[4]));
%[ [ 1, 7, 3, 19, 21, 9 ], [ 2, 8, 20 ], [ 10, 22, 18, 16, 12, 15, 13, 24 ], [ 11, 23, 17, 14 ] ]
$$
G=H_{34,3}=
$$
$$
[(\alpha_{1,1}\alpha_{1,3}\alpha_{1,7})
(\alpha_{2,1}\alpha_{2,3}\alpha_{2,7})
(\alpha_{3,1}\alpha_{3,3}\alpha_{3,7})
(\alpha_{1,4}\alpha_{1,8}\alpha_{1,6})
$$
$$
(\alpha_{2,4}\alpha_{2,8}\alpha_{2,6})
(\alpha_{3,4}\alpha_{3,8}\alpha_{3,6}),
$$
$$
(\alpha_{1,1}\alpha_{3,1})
(\alpha_{1,3}\alpha_{3,7}\alpha_{3,3}\alpha_{1,7})
(\alpha_{2,3}\alpha_{2,7})
(\alpha_{1,4}\alpha_{3,6}\alpha_{3,5}\alpha_{3,8})
$$
$$
(\alpha_{2,4}\alpha_{2,6}\alpha_{2,5}\alpha_{2,8})
(\alpha_{3,4}\alpha_{1,6}\alpha_{1,5}\alpha_{1,8}) ]
$$
with orbits
%$\{\alpha_{1,1},\alpha_{1,3},\alpha_{3,1},\alpha_{1,7},\alpha_{3,7},\alpha_{3,3}\},\
%\{\alpha_{2,1},\alpha_{2,3},\alpha_{2,7}\},\
%\{\alpha_{1,4}, \linebreak
%\alpha_{1,8},\alpha_{3,6},\alpha_{1,6},
%\alpha_{3,4},\alpha_{3,5},\alpha_{1,5},\alpha_{3,8}\},\
%\{\alpha_{2,4},\alpha_{2,8},\alpha_{2,6},\alpha_{2,5}\}$.
$\{\alpha_{2,4},\alpha_{2,8},\alpha_{2,6},\alpha_{2,5}\}$,
$\{\alpha_{1,4},
\alpha_{1,8},\alpha_{3,6},\alpha_{1,6},
\alpha_{3,4},\alpha_{3,5},\alpha_{1,5},\alpha_{3,8}\}$.
%gap> H1:=Image(F[1]);
%Group([ (1,7,19)(2,8,20)(3,9,21)(10,22,16)(11,23,17)(12,24,18), (7,9)(10,15)(11,14)(12,13)(16,22)(17,23)(18,24)
%(19,21) ])
%gap> Orbits(H1);
%[ [ 1, 7, 19, 9, 21, 3 ], [ 2, 8, 20 ], [ 10, 22, 15, 16 ], [ 11, 23, 14, 17 ], [ 12, 24, 13, 18 ] ]
The $G_1\cong \AAA_{4}$ is marked by
$$
G=H_{34,3}\supset G_1=
$$
$$
[(\alpha_{1,1}\alpha_{1,3}\alpha_{1,7})(\alpha_{2,1}\alpha_{2,3}\alpha_{2,7})(\alpha_{3,1}\alpha_{3,3}\alpha_{3,7})
(\alpha_{1,4}\alpha_{1,8}\alpha_{1,6})(\alpha_{2,4}\alpha_{2,8}\alpha_{2,6})(\alpha_{3,4}\alpha_{3,8}\alpha_{3,6}),
$$
$$
(\alpha_{1,3}\alpha_{3,3})(\alpha_{1,4}\alpha_{3,5})(\alpha_{2,4}\alpha_{2,5})(\alpha_{3,4}\alpha_{1,5})
(\alpha_{1,6}\alpha_{1,8})(\alpha_{2,6}\alpha_{2,8})(\alpha_{3,6}\alpha_{3,8})(\alpha_{1,7}\alpha_{3,7})]
$$
with suborbits
$\{\alpha_{1,4},\alpha_{1,8},\alpha_{3,5},\alpha_{1,6}\}$,
$\{\alpha_{2,4},\alpha_{2,8},\alpha_{2,5},\alpha_{2,6}\}$,
$\{\alpha_{3,4},\alpha_{3,8},\alpha_{1,5},\alpha_{3,6}\}$.
Both $G$ and $G_1$ are marked by $\ast$ in Tables 1---4.

\medskip

Case 59:
$({\bf n}=17,\ (4\aaa_1,4\aaa_1,6\aaa_1)\subset 14\aaa_1)\
\Longleftarrow ({\bf n}=34,\ ((6\aaa_1)_I,8\aaa_1)\subset 14\aaa_1)$.

Similar to Cases 1 and 44. By \cite{Nik10}, the $G\cong \SSS_4$
is marked by $N_{23}$ and
%%G1:=Group([(1,18,15,22)(2,4,20,11)(3,7)(8,10,16,24)(9,12,23,13)(14,19),
%%(2,13)(3,22)(4,9)(7,18)(8,14)(11,20)(12,23)(16,19)]);
$G=H_{34,1}$ of Case 44
with orbits
%$
%\{\alpha_{1},\alpha_{3},\alpha_{22},\alpha_{15},\alpha_{7},\alpha_{18}\},\
%\{\alpha_{2},\alpha_{13},\alpha_{20},\alpha_{9},
%\alpha_{12},\alpha_{4},\alpha_{23},
%\linebreak
%\alpha_{11}\},\
%\{\alpha_{8},\alpha_{14},\alpha_{16},\alpha_{10},
%\alpha_{19},\alpha_{24}\}$;
$\{\alpha_{1},\alpha_{3},\alpha_{22},\alpha_{15},\alpha_{7},\alpha_{18}\}$,
$\{\alpha_{2},\alpha_{13},\alpha_{20},\alpha_{9},\alpha_{12},\alpha_{4},\alpha_{23},\alpha_{11}\}$.
%gap> H1:=Image(F[1]);
%Group([ (1,3,22)(4,12,11)(7,18,15)(8,10,14)(9,20,23)(16,24,19), (2,9)(3,7)(4,13)(8,16)(11,12)(14,19)(18,22)(20,23) ])
%gap> Orbits(H1);
%[ [ 1, 3, 22, 7, 18, 15 ], [ 2, 9, 20, 23 ], [ 4, 12, 13, 11 ], [ 8, 10, 16, 14, 24, 19 ] ]
The $G_1\cong \AAA_{4}$ is marked by
$G_1\subset G=H_{34,1}$ of Case 44 with suborbits
$\{\alpha_{2},\alpha_{9},\alpha_{20},\alpha_{23}\}$, $\{\alpha_{4},\alpha_{12},\alpha_{13},\alpha_{11}\}$,
$\{\alpha_{1},\alpha_{3},\alpha_{22},\alpha_{15},\alpha_{7},\alpha_{18}\}$.
Both $G$ and $G_1$ are marked by $\ast$ in Tables 1---4.

\medskip

Case 60:
$({\bf n}=17,\ (4\aaa_1,6\aaa_1,6\aaa_1)
\subset 16\aaa_1)\Longleftarrow ({\bf n}=75,\ 16\aaa_1)$.

Similar to Cases 1 and 10. By \cite{Nik10},
the $G\cong 4^2\AAA_4$ is marked by $N_{23}$ and
$G=H_{75,1}$ from Case 10 with the orbit
$\{\alpha_{1},\alpha_{9},\alpha_{6},\alpha_{7},\alpha_{4},
\alpha_{21},\alpha_{20},
\alpha_{22},\alpha_{5},\alpha_{14},\alpha_{10},\alpha_{23},
\alpha_{3},\alpha_{24},\alpha_{18},\alpha_{16}\}$.
The $G_1\cong \AAA_4$ is marked by
$$
G=H_{75,1}\supset G_1=
$$
%gap> H2:=Image(F[2]);
%Group([ (1,3,10)(4,6,21)(7,24,14)(9,16,20)(11,13,17)(18,23,22), (3,16)(4,5)(6,21)(10,20)(11,12)(13,17)(14,22)
%(23,24) ])
%gap> Orbits(H2);
%[ [ 1, 3, 10, 16, 20, 9 ], [ 4, 6, 5, 21 ], [ 7, 24, 14, 23, 22, 18 ], [ 11, 13, 12, 17 ] ]
$$
[(\alpha_{1}\alpha_{3}\alpha_{10})(\alpha_{4}\alpha_{6}\alpha_{21})(\alpha_{7}\alpha_{24}\alpha_{14})
(\alpha_{9}\alpha_{16}\alpha_{20})(\alpha_{11}\alpha_{13}\alpha_{17})(\alpha_{18}\alpha_{23}\alpha_{22}),
$$
$$
(\alpha_{3}\alpha_{16})(\alpha_{4}\alpha_{5})(\alpha_{6}\alpha_{21})(\alpha_{10}\alpha_{20})
(\alpha_{11}\alpha_{12})(\alpha_{13}\alpha_{17})(\alpha_{14}\alpha_{22})(\alpha_{23}\alpha_{24})]
$$
with suborbits
$\{\alpha_{4},\alpha_{6},\alpha_{5},\alpha_{21}\}$,
$\{\alpha_{1},\alpha_{3},\alpha_{10},\alpha_{16},\alpha_{20},\alpha_{9}\}$,
$\{\alpha_{7},\alpha_{24},\alpha_{14},\alpha_{23},\alpha_{22},\alpha_{18}\}$.
Both $G$ and $G_1$ are marked by $\ast$ in Tables 1---4.

\medskip

Case 61:
$({\bf n}=17,\ (4\aaa_1,6\aaa_1,6\aaa_1)\subset 4\aaa_1\amalg 6\aaa_2)
\Longleftarrow ({\bf n}=34,\ (4\aaa_1,6\aaa_2)\subset 4\aaa_1\amalg 6\aaa_2)$.

Similar to Cases 1, 22 and 55. By \cite{Nik11}, the $G\cong \SSS_4$ is marked
by  $N_{22}$ and $G=H_{34,1}$ of Case 22
with orbits
$\{\alpha_{1,1}, \alpha_{1,3}, \alpha_{2,4}, \alpha_{1,8}\}$, $\{\alpha_{1,2},\alpha_{2,5},
\alpha_{1,5},\alpha_{1,11},\alpha_{2,2},\alpha_{2,11},
\alpha_{1,12},\alpha_{2,12},\alpha_{2,7},\linebreak \alpha_{1,7}, 
\alpha_{1,9},\alpha_{2,9}\}$.
The $G_1\cong \AAA_4$ is marked by
$G_1\subset G=H_{34,1}$ of Case 55
with suborbits
$\{\alpha_{1,1}, \alpha_{1,3}, \alpha_{2,4}, \alpha_{1,8}\}$,
$\{\alpha_{1,2},\alpha_{2,9},\alpha_{2,7},\alpha_{1,12},\alpha_{1,5},\alpha_{2,11}\}$,
$\{\alpha_{2,2},\alpha_{1,9},\alpha_{1,7},\alpha_{2,12},\alpha_{2,5},\alpha_{1,11}\}$.
Both $G$ and $G_1$ are marked by $\ast$ in Tables 1---4.

\medskip

Case 62:
$({\bf n}=21,\ (4\aaa_1, 4\aaa_1,4\aaa_1)
\subset 12\aaa_1)\Longleftarrow ({\bf n}=49,\ 12\aaa_1)$.

Similar to Cases 1 and 5. By \cite{Nik10}, the
$G\cong 2^4C_3$ is marked by $N_{23}$ and
$G=H_{49,1}$ of Case 5 with the orbit
$\{\alpha_{1},\alpha_{22},\alpha_{21},\alpha_{10},\alpha_{19},\alpha_{14},
\alpha_{8},\alpha_{17},\alpha_{18},\alpha_{7},\alpha_{3},\alpha_{16}\}$.
The $G_1\cong (C_2)^4$ is marked by
$$
G=H_{49,1}\supset G_1=
$$
%gap> GeneratorsOfGroup(H1);
%[ (1,21)(2,23)(3,8)(4,9)(10,17)(11,20)(12,13)(14,22), (2,12)(3,8)(4,20)(7,16)(9,11)(13,23)(14,22)(18,19),
%  (1,17)(2,12)(3,22)(4,11)(8,14)(9,20)(10,21)(13,23), (2,13)(3,22)(4,9)(7,18)(8,14)(11,20)(12,23)(16,19) ]
$$
[(\alpha_{1}\alpha_{21})(\alpha_{2}\alpha_{23})(\alpha_{3}\alpha_{8})(\alpha_{4}\alpha_{9})
(\alpha_{10}\alpha_{17})(\alpha_{11}\alpha_{20})(\alpha_{12}\alpha_{13})(\alpha_{14}\alpha_{22}),
$$
$$
(\alpha_{2}\alpha_{12})(\alpha_{3}\alpha_{8})(\alpha_{4}\alpha_{20})(\alpha_{7}\alpha_{16})
(\alpha_{9}\alpha_{11})(\alpha_{13}\alpha_{23})(\alpha_{14}\alpha_{22})(\alpha_{18}\alpha_{19}),
$$
$$
(\alpha_{1}\alpha_{17})(\alpha_{2}\alpha_{12})(\alpha_{3}\alpha_{22})(\alpha_{4}\alpha_{11})
(\alpha_{8}\alpha_{14})(\alpha_{9}\alpha_{20})(\alpha_{10}\alpha_{21})(\alpha_{13}\alpha_{23}),
$$
$$
(\alpha_{2}\alpha_{13})(\alpha_{3}\alpha_{22})(\alpha_{4}\alpha_{9})(\alpha_{7}\alpha_{18})
(\alpha_{8}\alpha_{14})(\alpha_{11}\alpha_{20})(\alpha_{12}\alpha_{23})(\alpha_{16}\alpha_{19})]
$$
with suborbits
%gap> Orbits(H1);
%[ [ 1, 21, 17, 10 ], [ 3, 8, 22, 14 ], [ 4, 9, 20, 11 ], [ 7, 16, 18, 19 ] ]
$\{\alpha_{1},\alpha_{21},\alpha_{17},\alpha_{10} \}$, $\{\alpha_{3},\alpha_{8},\alpha_{22},\alpha_{14} \}$,
$\{\alpha_{7},\alpha_{16},\alpha_{18},\alpha_{19}\}$.
Both $G$ and $G_1$ are marked by $\ast$ in Tables 1---4.

\medskip

Case 63:
$({\bf n}=21,\ (4\aaa_1, 4\aaa_1,4\aaa_1,4\aaa_1)
\subset 16\aaa_1)\Longleftarrow ({\bf n}=75,\ 16\aaa_1)$.

Similar to Cases 1 and 10. By \cite{Nik10}, the $G\cong 4^2\AAA_4$
is marked by $N_{23}$ and
$G=H_{75,1}$ of Case 10 with the orbit
$\{\alpha_{1},\alpha_{9},\alpha_{6},\alpha_{7},\alpha_{4},
\alpha_{21},\alpha_{20},
\alpha_{22},\alpha_{5},\alpha_{14},\alpha_{10},\alpha_{23},
\alpha_{3},\alpha_{24},\alpha_{18},\alpha_{16}\}$.
The $G_1\cong (C_2)^4$ is marked by
$$
G=H_{75,1}\supset G_1=
$$
%gap> GeneratorsOfGroup(H1);
%[ (1,7)(3,16)(9,18)(10,22)(11,17)(12,13)(14,20)(23,24), (1,9)(3,23)(7,18)(10,14)(11,12)(13,17)(16,24)(20,22),
%  (3,24)(4,6)(5,21)(10,14)(11,17)(12,13)(16,23)(20,22), (3,16)(4,5)(6,21)(10,20)(11,12)(13,17)(14,22)(23,24) ]
%gap> Orbits(H1);
%[ [ 1, 7, 9, 18 ], [ 3, 16, 23, 24 ], [ 4, 6, 5, 21 ], [ 10, 22, 14, 20 ], [ 11, 17, 12, 13 ] ]
$$
[(\alpha_{1}\alpha_{7})(\alpha_{3}\alpha_{16})(\alpha_{9}\alpha_{18})(\alpha_{10}\alpha_{22})
(\alpha_{11}\alpha_{17})(\alpha_{12}\alpha_{13})(\alpha_{14}\alpha_{20})(\alpha_{23}\alpha_{24}),
$$
$$
(\alpha_{1}\alpha_{9})(\alpha_{3}\alpha_{23})(\alpha_{7}\alpha_{18})(\alpha_{10}\alpha_{14})
(\alpha_{11}\alpha_{12})(\alpha_{13}\alpha_{17})(\alpha_{16}\alpha_{24})(\alpha_{20}\alpha_{22}),
$$
$$
(\alpha_{3}\alpha_{24})(\alpha_{4}\alpha_{6})(\alpha_{5}\alpha_{21})(\alpha_{10}\alpha_{14})
(\alpha_{11}\alpha_{17})(\alpha_{12}\alpha_{13})(\alpha_{16}\alpha_{23})(\alpha_{20}\alpha_{22}),
$$
$$
(\alpha_{3}\alpha_{16})(\alpha_{4}\alpha_{5})(\alpha_{6}\alpha_{21})(\alpha_{10}\alpha_{20})
(\alpha_{11}\alpha_{12})(\alpha_{13}\alpha_{17})(\alpha_{14}\alpha_{22})(\alpha_{23}\alpha_{24})]
$$
with suborbits
$\{\alpha_{1},\alpha_{7},\alpha_{9},\alpha_{18}\}$,
$\{\alpha_{3},\alpha_{16},\alpha_{23},\alpha_{24}\}$,
$\{\alpha_{4},\alpha_{6},\alpha_{5},\alpha_{21}\}$,
$\{\alpha_{10},\alpha_{22},\alpha_{14},\alpha_{20}\}$. 
\newline 
Both $G$ and $G_1$ are marked by $\ast$ in Tables 1---4.

\medskip

Case 64:
$({\bf n}=22,\ (2\aaa_1,2\aaa_1)
\subset 4\aaa_1)\Longleftarrow ({\bf n}=39,\ 4\aaa_1)$.

Similar to Cases 1 and 14. By \cite{Nik10}, the
$G\cong 2^4C_2$ is marked by $N_{23}$ and
$G=H_{39,2}$ of Case 14 with the orbit
$\{\alpha_{2},\alpha_{12},\alpha_{13},\alpha_{23}\}$.
%gap> H1:=Image(F[1]);
%Group([ (3,14,8,22)(4,20,9,11)(5,7,15,19)(6,18,24,16)(10,17)(12,13), (5,18)(6,7)(10,17)(11,20)(12,13)(14,22)(15,16)
%(19,24), (2,23)(4,9)(5,24)(6,15)(7,16)(11,20)(12,13)(18,19) ])
%gap> Orbits(H1);
%[ [ 2, 23 ], [ 3, 14, 8, 22 ], [ 4, 20, 9, 11 ], [ 5, 7, 18, 24, 15, 6, 16, 19 ], [ 10, 17 ], [ 12, 13 ] ]
The $G_1\cong C_2\times D_8$ is marked by
$$
G=H_{39,2}\supset G_1=
$$
$$
[(\alpha_{3}\alpha_{14}\alpha_{8}\alpha_{22})(\alpha_{4}\alpha_{20}\alpha_{9}\alpha_{11})(\alpha_{5}\alpha_{7}\alpha_{15}\alpha_{19})
(\alpha_{6}\alpha_{18}\alpha_{24}\alpha_{16})(\alpha_{10}\alpha_{17})(\alpha_{12}\alpha_{13}),
$$
$$
(\alpha_{5}\alpha_{18})(\alpha_{6}\alpha_{7})(\alpha_{10}\alpha_{17})(\alpha_{11}\alpha_{20})
(\alpha_{12}\alpha_{13})(\alpha_{14}\alpha_{22})(\alpha_{15}\alpha_{16})(\alpha_{19}\alpha_{24}),
$$
$$
(\alpha_{2}\alpha_{23})(\alpha_{4}\alpha_{9})(\alpha_{5}\alpha_{24})(\alpha_{6}\alpha_{15})
(\alpha_{7}\alpha_{16})(\alpha_{11}\alpha_{20})(\alpha_{12}\alpha_{13})(\alpha_{18}\alpha_{19})]
$$
with suborbits
$\{\alpha_{2},\alpha_{23}\}$, $\{\alpha_{12},\alpha_{13}\}$.
Both $G$ and $G_1$ are marked by $\ast$ in Tables 1---4.

\medskip

Case 65:
$({\bf n}=22,\ ((4\aaa_1,4\aaa_1)
\subset 8\aaa_1)_I)\Longleftarrow ({\bf n}=40,\ 8\aaa_1)$.

Similar to cases 1 and 4. By \cite{Nik10}, the
$G\cong Q_8\ast Q_8$ is marked by $N_{23}$ and
$G=H_{40,1}$ of Case 4 with the orbit
$\{\alpha_{1},\alpha_{16},\alpha_{14},\alpha_{19},\alpha_{18},\alpha_{20},
\alpha_{24},\alpha_{2}\}$.
The $G_1\cong C_2\times D_8$ is
marked by
$$
G=H_{40,1}\supset G_1=
$$
$$
[(\alpha_{3}\alpha_{5})(\alpha_{6}\alpha_{12})(\alpha_{8}\alpha_{13})(\alpha_{10}\alpha_{22})
(\alpha_{15}\alpha_{21})(\alpha_{16}\alpha_{19})(\alpha_{17}\alpha_{23})(\alpha_{18}\alpha_{20}),
$$
$$
(\alpha_{1}\alpha_{16})(\alpha_{2}\alpha_{18})(\alpha_{4}\alpha_{9})(\alpha_{5}\alpha_{23})
(\alpha_{10}\alpha_{12})(\alpha_{14}\alpha_{20})(\alpha_{15}\alpha_{21})(\alpha_{19}\alpha_{24}),
$$
$$
(\alpha_{2}\alpha_{14})(\alpha_{3}\alpha_{6})(\alpha_{4}\alpha_{9})(\alpha_{5}\alpha_{12})
(\alpha_{10}\alpha_{23})(\alpha_{15}\alpha_{21})(\alpha_{17}\alpha_{22})(\alpha_{18}\alpha_{20})]
$$
%[ (3,5)(6,12)(8,13)(10,22)(15,21)(16,19)(17,23)(18,20), (1,16)(2,18)(4,9)(5,23)(10,12)(14,20)(15,21)(19,24),
%(2,14)(3,6)(4,9)(5,12)(10,23)(15,21)(17,22)(18,20) ]
with suborbits
$\{\alpha_{1},\alpha_{16},\alpha_{19},\alpha_{24}\}$,
$\{\alpha_{2},\alpha_{18},\alpha_{14},\alpha_{20}\}$.
%[ [ 1, 16, 19, 24 ], [ 2, 18, 14, 20 ],
Both $G$ and $G_1$ are marked by $\ast$ in Tables 1---4.

\medskip

Case 66:
$({\bf n}=22,\ (8\aaa_1,8\aaa_1)
\subset 16\aaa_1)\ \Longleftarrow ({\bf n}=39,\ 16\aaa_1)$.

Similar to Cases 1 and 7. By \cite{Nik10}, the
$G\cong 2^4C_2$ is marked by $N_{23}$ and
$G=H_{39,1}$ of Case 7 with the orbit
$\{\alpha_{2},\alpha_{12},\alpha_{13},\alpha_{3},\alpha_{18},\alpha_{8},
\alpha_{23},\alpha_{19},\alpha_{7},\alpha_{22},\alpha_{4},
\alpha_{14},\alpha_{20},\alpha_{11},\alpha_{16},\alpha_{9}\}$.
The $G_1\cong C_2\times D_8$ is marked by
$$
G=H_{39,1}\supset G_1=
$$
%gap> H1:=Image(F[1]);
%Group([ (5,10)(7,13)(8,12)(9,22)(11,23)(14,16)(17,21)(19,20), (2,7)(3,9)(4,22)(8,11)(12,16)(13,18)(14,20)
%(19,23), (2,3)(4,18)(7,9)(8,12)(11,16)(13,22)(14,23)(19,20) ])
$$
[(\alpha_{5}\alpha_{10})(\alpha_{7}\alpha_{13})(\alpha_{8}\alpha_{12})(\alpha_{9}\alpha_{22})
(\alpha_{11}\alpha_{23})(\alpha_{14}\alpha_{16})(\alpha_{17}\alpha_{21})(\alpha_{19}\alpha_{20}),
$$
$$
(\alpha_{2}\alpha_{7})(\alpha_{3}\alpha_{9})(\alpha_{4}\alpha_{22})(\alpha_{8}\alpha_{11})
(\alpha_{12}\alpha_{16})(\alpha_{13}\alpha_{18})(\alpha_{14}\alpha_{20})(\alpha_{19}\alpha_{23}),
$$
$$
(\alpha_{2}\alpha_{3})(\alpha_{4}\alpha_{18})(\alpha_{7}\alpha_{9})(\alpha_{8}\alpha_{12})
(\alpha_{11}\alpha_{16})(\alpha_{13}\alpha_{22})(\alpha_{14}\alpha_{23})(\alpha_{19}\alpha_{20})]
$$
with suborbits
%gap> Orbits(H1);
%[ [ 2, 7, 3, 13, 9, 18, 22, 4 ], [ 5, 10 ], [ 8, 12, 11, 16, 23, 14, 19, 20 ], [ 17, 21 ] ]
$\{\alpha_{2},\alpha_{7},\alpha_{3},\alpha_{13},\alpha_{9},\alpha_{18},\alpha_{22},\alpha_{4}\}$,
$\{\alpha_{8},\alpha_{12},\alpha_{11},\alpha_{16},\alpha_{23},\alpha_{14},\alpha_{19},\alpha_{20}\}$.
Both $G$ and $G_1$ are marked by $\ast$ in Tables 1---4.

\medskip

Case 67:
$({\bf n}=22,\ (2\aaa_1, 2\aaa_1,4\aaa_1)
\subset 8\aaa_1)\Longleftarrow ({\bf n}=56,\ 8\aaa_1)$.

Similar to Cases 1 and 23. By \cite{Nik10}, the $G\cong \Gamma_{15}a_1$
is marked by $N_{23}$ and $G=H_{56,1}$ of Case 23 with the orbit
$\{\alpha_{1}, \alpha_{14}, \alpha_{24}, \alpha_{17}, \alpha_{2},\alpha_{23}, \alpha_{5},
\alpha_{3} \}$.
%gap> H4:=Image(F[4]);
%Group([ (1,17)(2,3)(6,9)(7,10)(8,13)(15,20)(16,18)(19,21), (2,17)(6,16)(7,8)(9,13)(10,18)(12,22)(19,20)(23,24), (5,14)
%(6,9)(7,10)(8,18)(13,16)(15,21)(19,20)(23,24) ])
%gap> Orbits(H4);
%[ [ 1, 17, 2, 3 ], [ 5, 14 ], [ 6, 9, 16, 13, 18, 8, 10, 7 ], [ 12, 22 ], [ 15, 20, 21, 19 ], [ 23, 24 ] ]
The $G_1\cong C_2\times D_8$ is marked by
$$
G=H_{56,1}\supset G_1=
$$
$$
[(\alpha_{1}\alpha_{17})(\alpha_{2}\alpha_{3})(\alpha_{6}\alpha_{9})(\alpha_{7}\alpha_{10})
(\alpha_{8}\alpha_{13})(\alpha_{15}\alpha_{20})(\alpha_{16}\alpha_{18})(\alpha_{19}\alpha_{21}),
$$
$$
(\alpha_{2}\alpha_{17})(\alpha_{6}\alpha_{16})(\alpha_{7}\alpha_{8})(\alpha_{9}\alpha_{13})
(\alpha_{10}\alpha_{18})(\alpha_{12}\alpha_{22})(\alpha_{19}\alpha_{20})(\alpha_{23}\alpha_{24}),
$$
$$
(\alpha_{5}\alpha_{14})(\alpha_{6}\alpha_{9})(\alpha_{7}\alpha_{10})(\alpha_{8}\alpha_{18})
(\alpha_{13}\alpha_{16})(\alpha_{15}\alpha_{21})(\alpha_{19}\alpha_{20})(\alpha_{23}\alpha_{24})]
$$
with suborbits
$\{\alpha_{5},\alpha_{14}\}$, $\{\alpha_{23},\alpha_{24}\}$,
$\{\alpha_{1},\alpha_{17},\alpha_{2},\alpha_{3}\}$.
Both $G$ and $G_1$ are marked by $\ast$ in Tables 1---4.

\medskip

Case 68:
$({\bf n}=22,\ (2\aaa_1, 2\aaa_1,8\aaa_1)
\subset 12\aaa_1)\Longleftarrow ({\bf n}=65,\ 12\aaa_1)$.

Similar to Cases 1 and 24. By \cite{Nik10}, the $G\cong 2^4D_6$
is marked by $N_{23}$ and
$G=H_{65,3}$ of Case 24 with the orbit
$\{\alpha_{1},\alpha_{7},\alpha_{23},\alpha_{8},\alpha_{11},\alpha_{17},\alpha_{20},
\alpha_{9}, \alpha_{6},\alpha_{4},\alpha_{14},\alpha_{13}\}$.
%gap> H3:=Image(F[3]);
%Group([ (1,17)(4,23)(5,22)(6,9)(8,13)(11,14)(12,24)(19,21), (1,8)(3,22)(4,13)(5,16)(10,21)(11,23)(14,17)(15,19), (1,8)
%(4,14)(6,9)(7,20)(10,15)(11,23)(13,17)(19,21) ])
%gap> Orbits(H3);
%[ [ 1, 17, 8, 14, 13, 11, 4, 23 ], [ 3, 22, 5, 16 ], [ 6, 9 ], [ 7, 20 ], [ 10, 21, 15, 19 ], [ 12, 24 ] ]
The $G_1\cong C_2\times D_8$ is marked by
$$
G=H_{65,3}\supset G_1=
$$
$$
[(\alpha_{1}\alpha_{17})(\alpha_{4}\alpha_{23})(\alpha_{5}\alpha_{22})(\alpha_{6}\alpha_{9})
(\alpha_{8}\alpha_{13})(\alpha_{11}\alpha_{14})(\alpha_{12}\alpha_{24})(\alpha_{19}\alpha_{21}),
$$
$$
(\alpha_{1}\alpha_{8})(\alpha_{3}\alpha_{22})(\alpha_{4}\alpha_{13})(\alpha_{5}\alpha_{16})
(\alpha_{10}\alpha_{21})(\alpha_{11}\alpha_{23})(\alpha_{14}\alpha_{17})(\alpha_{15}\alpha_{19}),
$$
$$
(\alpha_{1}\alpha_{8})(\alpha_{4}\alpha_{14})(\alpha_{6}\alpha_{9})(\alpha_{7}\alpha_{20})
(\alpha_{10}\alpha_{15})(\alpha_{11}\alpha_{23})(\alpha_{13}\alpha_{17})(\alpha_{19}\alpha_{21})]
$$
with suborbits
$\{\alpha_{6},\alpha_{9}\}$, $\{\alpha_{7},\alpha_{20}\}$,
$\{\alpha_{1},\alpha_{17},\alpha_{8},\alpha_{14},\alpha_{13},\alpha_{11},\alpha_{4},\alpha_{23}\}$.
Both $G$ and $G_1$ are marked by $\ast$ in Tables 1---4.

\medskip

Case 69:
$({\bf n}=22,\ ((4\aaa_1,4\aaa_1)_I,8\aaa_1)
\subset 16\aaa_1)\Longleftarrow ({\bf n}=56,\ 16\aaa_1)$.

Similar to cases 1 and 8. By \cite{Nik10}, the
$G\cong \Gamma_{25}a_1$
is marked by $N_{23}$ and
$G=H_{56,2}$ of Case 8 with the orbit
$\{\alpha_{2}, \alpha_{3}, \alpha_{23},
\alpha_{24}, \alpha_{5},\alpha_{8}, \alpha_{17}, \alpha_{19}, \alpha_{6},
\alpha_{7}, \alpha_{11}, \alpha_{20}, \alpha_{18},
\alpha_{10},\alpha_{16},\alpha_{15}\}$.
The $G_1\cong C_2\times D_8$ is marked by
$$
G=H_{56,2}\supset G_1=
$$
%gap> H3:=Image(F[3]);
%Group([ (3,11)(5,15)(7,20)(8,19)(9,13)(10,24)(12,22)(18,23), (2,10)(3,20)(5,19)(6,24)(7,15)(8,11)(16,18)
%(17,23), (1,14)(3,20)(5,19)(7,11)(8,15)(12,22)(16,17)(18,23) ])
%gap> Orbits(H3);
%[ [ 1, 14 ], [ 2, 10, 24, 6 ], [ 3, 11, 20, 8, 7, 19, 15, 5 ], [ 9, 13 ], [ 12, 22 ], [ 16, 18, 17, 23 ] ]
$$
(\alpha_{3}\alpha_{11})(\alpha_{5}\alpha_{15})(\alpha_{7}\alpha_{20})(\alpha_{8}\alpha_{19})
(\alpha_{9}\alpha_{13})(\alpha_{10}\alpha_{24})(\alpha_{12}\alpha_{22})(\alpha_{18}\alpha_{23}),
$$
$$
(\alpha_{2}\alpha_{10})(\alpha_{3}\alpha_{20})(\alpha_{5}\alpha_{19})(\alpha_{6}\alpha_{24})
(\alpha_{7}\alpha_{15})(\alpha_{8}\alpha_{11})(\alpha_{16}\alpha_{18})(\alpha_{17}\alpha_{23}),
$$
$$
(\alpha_{1}\alpha_{14})(\alpha_{3}\alpha_{20})(\alpha_{5}\alpha_{19})(\alpha_{7}\alpha_{11})
(\alpha_{8}\alpha_{15})(\alpha_{12}\alpha_{22})(\alpha_{16}\alpha_{17})(\alpha_{18}\alpha_{23})]
$$
with suborbits
$\{\alpha_{2},\alpha_{10},\alpha_{24},\alpha_{6}\}$, $\{\alpha_{16},\alpha_{18},\alpha_{17},\alpha_{23}\}$,
$\{\alpha_{3},\alpha_{11},\alpha_{20},\alpha_{8},\alpha_{7},\alpha_{19},\alpha_{15},\alpha_{5}\}$.
Both $G$ and $G_1$ are marked by $\ast$ in Tables 1---4.

\medskip

Case 70:
$({\bf n}=34,\ (\aaa_1,\aaa_1)\subset 2\aaa_1)\Longleftarrow ({\bf n}=51,\ 2\aaa_1)$.

Similar to Case 1. By \cite{Nik10}, the $G\cong C_2\times \SSS_4$
is marked by $N_{23}$ and
%G:=Group([(1,15,18,20)(3,6,23,17)(4,19)(7,21,24,11)(8,9,14,16)(12,22),
%(2,13)(3,22)(4,9)(7,18)(8,14)(11,20)(12,23)(16,19)]);
$$
G=H_{51,3}=
$$
$$
[(\alpha_{1}\alpha_{15}\alpha_{18}\alpha_{20})
(\alpha_{3}\alpha_{6}\alpha_{23}\alpha_{17})
(\alpha_{4}\alpha_{19})
(\alpha_{7}\alpha_{21}\alpha_{24}\alpha_{11})
(\alpha_{8}\alpha_{9}\alpha_{14}\alpha_{16})
(\alpha_{12}\alpha_{22}),
$$
$$
(\alpha_{2}\alpha_{13})(\alpha_{3}\alpha_{22})
(\alpha_{4}\alpha_{9})(\alpha_{7}\alpha_{18})
(\alpha_{8}\alpha_{14})(\alpha_{11}\alpha_{20})
(\alpha_{12}\alpha_{23})(\alpha_{16}\alpha_{19})]
$$
with the orbit
%$\{\alpha_{1},\alpha_{15},\alpha_{21},\alpha_{24},
%\alpha_{18},\alpha_{7},\alpha_{20},\alpha_{11}\},\
%\{\alpha_{2},\alpha_{13}\},\
%\{\alpha_{3},\alpha_{12},\alpha_{6},
%\linebreak
%\alpha_{23},\alpha_{22},\alpha_{17}\},\
%\{\alpha_{4},\alpha_{9},
%\alpha_{19},\alpha_{16},\alpha_{14},\alpha_{8}\}$;
$\{\alpha_{2},\alpha_{13}\}$.
%Group([ (3,6,22)(4,16,8)(7,24,18)(9,14,19)(11,20,15)(12,23,17), (1,11,24,20)(3,22,23,12)(4,9,19,16)(6,17)(7,21,18,15)
%(8,14) ])
%gap> Orbits(H2);
%[ [ 1, 11, 20, 24, 15, 18, 7, 21 ], [ 3, 6, 22, 17, 23, 12 ], [ 4, 16, 9, 8, 14, 19 ] ]
The $G_1\cong \SSS_4$ is marked by
$$
G=H_{51,3}\supset G_1=
$$
$$
[(\alpha_{3}\alpha_{6}\alpha_{22})(\alpha_{4}\alpha_{16}\alpha_{8})(\alpha_{7}\alpha_{24}\alpha_{18})
(\alpha_{9}\alpha_{14}\alpha_{19})(\alpha_{11}\alpha_{20}\alpha_{15})(\alpha_{12}\alpha_{23}\alpha_{17}),
$$
$$
(\alpha_{1}\alpha_{11}\alpha_{24}\alpha_{20})(\alpha_{3}\alpha_{22}\alpha_{23}\alpha_{12})
(\alpha_{4}\alpha_{9}\alpha_{19}\alpha_{16})(\alpha_{6}\alpha_{17})
(\alpha_{7}\alpha_{21}\alpha_{18}\alpha_{15})(\alpha_{8}\alpha_{14})]
$$
with suborbits
$\{\alpha_{2}\}$, $\{\alpha_{13}\}$.
Both $G$ and $G_1$ are marked by $\ast$ in Tables 1---4.

\medskip

Case 71:
$({\bf n}=34,\ (\aaa_1,2\aaa_1)\subset 3\aaa_1)\Longleftarrow ({\bf n}=61,\ 3\aaa_1)$.

Similar to Cases 1 and 48. By \cite{Nik10}, $G\cong \AAA_{4,3}$
is marked by by $N_{23}$ and $G=H_{61,1}$ of Case 48
with the orbit $\{\alpha_{7},\alpha_{15},\alpha_{18}\}$.
%gap> H3:=Image(F[3]);
%Group([ (2,20,14)(3,17,22)(4,9,21)(8,11,13)(10,12,23)(16,19,24), (1,10,23,12)(2,3,21,9)(4,22,8,20)(7,15)(11,17,14,13)
%(16,19) ])
%gap> Orbits(H3);
%[ [ 1, 10, 12, 23 ], [ 2, 20, 3, 14, 4, 17, 21, 13, 9, 22, 8, 11 ], [ 7, 15 ], [ 16, 19, 24 ] ]
The $G_1\cong \SSS_4$ is marked by
$$
G=H_{61,1}\supset G_1=
$$
$$
[(\alpha_{2}\alpha_{20}\alpha_{14})(\alpha_{3}\alpha_{17}\alpha_{22})(\alpha_{4}\alpha_{9}\alpha_{21})
(\alpha_{8}\alpha_{11}\alpha_{13})(\alpha_{10}\alpha_{12}\alpha_{23})(\alpha_{16}\alpha_{19}\alpha_{24}),
$$
$$
(\alpha_{1}\alpha_{10}\alpha_{23}\alpha_{12})(\alpha_{2}\alpha_{3}\alpha_{21}\alpha_{9})(\alpha_{4}\alpha_{22}\alpha_{8}\alpha_{20})
(\alpha_{7}\alpha_{15})(\alpha_{11}\alpha_{17}\alpha_{14}\alpha_{13})(\alpha_{16}\alpha_{19})]
$$
with suborbits
$\{\alpha_{18}\}$, $\{\alpha_{7},\alpha_{15}\}$.
Both $G$ and $G_1$ are marked by $\ast$ in Tables 1---4.

\medskip

Case 72:
$({\bf n}=34,\ (\aaa_1,3\aaa_1)\subset 4\aaa_1)\Longleftarrow ({\bf n}=65,\ 4\aaa_1)$.

Similar to Cases 1 and 52. By \cite{Nik10}, the $G\cong 2^4D_6$
is marked by $N_{23}$ and
$G=H_{65,4}$ of Case 52 with the orbit
$\{\alpha_{3},\alpha_{22},\alpha_{24},\alpha_{12}\}$.
%gap> H1:=Image(F[1]);
%Group([ (4,20,13)(5,18,16)(7,9,17)(8,19,11)(12,22,24)(15,23,21), (1,8,21,19)(4,13,20,6)(5,18)(7,17,14,9)(10,15,11,23)
%(12,22) ])
%gap> Orbits(H1);
%[ [ 1, 8, 19, 21, 11, 15, 23, 10 ], [ 4, 20, 13, 6 ], [ 5, 18, 16 ], [ 7, 9, 17, 14 ], [ 12, 22, 24 ] ]
The $G_1\cong \SSS_4$ is marked by
$$
G=H_{65,4}\supset G_1=
$$
$$
[(\alpha_{4}\alpha_{20}\alpha_{13})(\alpha_{5}\alpha_{18}\alpha_{16})(\alpha_{7}\alpha_{9}\alpha_{17})
(\alpha_{8}\alpha_{19}\alpha_{11})(\alpha_{12}\alpha_{22}\alpha_{24})(\alpha_{15}\alpha_{23}\alpha_{21}),
$$
$$
(\alpha_{1}\alpha_{8}\alpha_{21}\alpha_{19})(\alpha_{4}\alpha_{13}\alpha_{20}\alpha_{6})(\alpha_{5}\alpha_{18})
(\alpha_{7}\alpha_{17}\alpha_{14}\alpha_{9})(\alpha_{10}\alpha_{15}\alpha_{11}\alpha_{23})(\alpha_{12}\alpha_{22})]
$$
with suborbits
$\{\alpha_{3}\}$, $\{\alpha_{22},\alpha_{24},\alpha_{12}\}$.
Both $G$ and $G_1$ are marked by $\ast$ in Tables 1---4.

\medskip

Case 73:
$({\bf n}=34,\ (4\aaa_1,4\aaa_1)\subset 8\aaa_1)\Longleftarrow ({\bf n}=51,\ 8\aaa_1)$.

Similar to Cases 1 and 70. By \cite{Nik10}, the $G\cong C_2\times \SSS_4$
is marked by $N_{23}$ and $G=H_{51,3}$ of Case 70 with the orbit
$\{\alpha_{1},\alpha_{15},\alpha_{21},\alpha_{24},
\alpha_{18},\alpha_{7},\alpha_{20},\alpha_{11}\}$.
%gap> H1:=Image(F[1]);
%Group([ (3,22,6)(4,8,16)(7,18,24)(9,19,14)(11,15,20)(12,17,23), (1,7,24,18)(2,13)(3,12,23,22)(4,9,19,16)(6,17)
%(11,21,20,15) ])
%gap> Orbits(H1);
%[ [ 1, 7, 18, 24 ], [ 2, 13 ], [ 3, 22, 12, 6, 17, 23 ], [ 4, 8, 9, 16, 19, 14 ], [ 11, 15, 21, 20 ] ]
The $G_1\cong \SSS_4$ is marked by
$$
G=H_{51,3}\supset G_1=
$$
$$
[(\alpha_{3}\alpha_{22}\alpha_{6})(\alpha_{4}\alpha_{8}\alpha_{16})(\alpha_{7}\alpha_{18}\alpha_{24})
(\alpha_{9}\alpha_{19}\alpha_{14})(\alpha_{11}\alpha_{15}\alpha_{20})(\alpha_{12}\alpha_{17}\alpha_{23}),
$$
$$
(\alpha_{1}\alpha_{7}\alpha_{24}\alpha_{18})(\alpha_{2}\alpha_{13})
(\alpha_{3}\alpha_{12}\alpha_{23}\alpha_{22})(\alpha_{4}\alpha_{9}\alpha_{19}\alpha_{16})
(\alpha_{6}\alpha_{17})(\alpha_{11}\alpha_{21}\alpha_{20}\alpha_{15})]
$$
with suborbits
$\{\alpha_1,\alpha_7,\alpha_{18},\alpha_{24}\}$,
$\{\alpha_{11},\alpha_{15},\alpha_{21},\alpha_{20}\}$.
Both $G$ and $G_1$ are marked by $\ast$ in Tables 1---4.

\medskip

Case 74:
$({\bf n}=34,\ (4\aaa_1,8\aaa_1)\subset 12\aaa_1)\Longleftarrow ({\bf n}=61,\ 12\aaa_1)$.

Similar to Cases 1 and 48. By \cite{Nik10}, the $G\cong \AAA_{4,3}$
is marked by by $N_{23}$ and $G=H_{61,1}$ of Case 48 with the orbit
$\{\alpha_{2},\alpha_{4},\alpha_{21},\alpha_{13},\alpha_{11},\alpha_{8},
\alpha_{20},\alpha_{17},\alpha_{22},\alpha_{14},\alpha_{9},\alpha_{3}\}$.
%gap> H2:=Image(F[2]);
%Group([ (2,13,21)(4,20,8)(7,15,18)(9,14,11)(10,12,23)(16,24,19), (1,10,23,12)(2,3,21,9)(4,22,8,20)(7,15)(11,17,14,13)
%(16,19) ])
%gap> Orbits(H2);
%[ [ 1, 10, 12, 23 ], [ 2, 13, 3, 21, 11, 9, 17, 14 ], [ 4, 20, 22, 8 ], [ 7, 15, 18 ], [ 16, 24, 19 ] ]
The $G_1\cong \SSS_4$ is marked by
$$
G=H_{61,1}\supset G_1=
$$
$$
[(\alpha_{2}\alpha_{13}\alpha_{21})(\alpha_{4}\alpha_{20}\alpha_{8})(\alpha_{7}\alpha_{15}\alpha_{18})
(\alpha_{9}\alpha_{14}\alpha_{11})(\alpha_{10}\alpha_{12}\alpha_{23})(\alpha_{16}\alpha_{24}\alpha_{19}),
$$
$$
(\alpha_{1}\alpha_{10}\alpha_{23}\alpha_{12})(\alpha_{2}\alpha_{3}\alpha_{21}\alpha_{9})
(\alpha_{4}\alpha_{22}\alpha_{8}\alpha_{20})(\alpha_{7}\alpha_{15})
(\alpha_{11}\alpha_{17}\alpha_{14}\alpha_{13})(\alpha_{16}\alpha_{19})]
$$
with suborbits
$\{\alpha_{4},\alpha_{20},\alpha_{22},\alpha_{8}\}$,
$\{\alpha_{2},\alpha_{13},\alpha_{3},\alpha_{21},\alpha_{11},\alpha_{9},\alpha_{17},\alpha_{14}\}$.
Both $G$ and $G_1$ are marked by $\ast$ in Tables 1---4.

\medskip

Case 75:
$({\bf n}=34,\ (4\aaa_1,12\aaa_1)\subset 16\aaa_1)\Longleftarrow ({\bf n}=65,\ 16\aaa_1)$.

Similar to Case 1. By \cite{Nik10}, the $G\cong 2^4D_6$
is marked by $N_{23}$ and
$$
G=H_{65,1}=
$$
%G:=Group([(1,6,21)(4,17,14)(7,15,8)(9,10,23)(11,20,19)(12,18,24),
%(1,11,23,8)(4,6,13,20)(5,16)(7,17,9,14)(10,21,19,15)(12,24)]);
$$
[(\alpha_{1}\alpha_{6}\alpha_{21})(\alpha_{4}\alpha_{17}\alpha_{14})
(\alpha_{7}\alpha_{15}\alpha_{8})(\alpha_{9}\alpha_{10}\alpha_{23})
(\alpha_{11}\alpha_{20}\alpha_{19})(\alpha_{12}\alpha_{18}\alpha_{24}),
$$
$$
(\alpha_{1}\alpha_{11}\alpha_{23}\alpha_{8})
(\alpha_{4}\alpha_{6}\alpha_{13}\alpha_{20})
(\alpha_{5}\alpha_{16})
(\alpha_{7}\alpha_{17}\alpha_{9}\alpha_{14})
(\alpha_{10}\alpha_{21}\alpha_{19}\alpha_{15})
(\alpha_{12}\alpha_{24})]
$$
with the orbit
%$\{\alpha_{1},\alpha_{19},\alpha_{6},\alpha_{7},\alpha_{4},
%\alpha_{17},\alpha_{20},\alpha_{11},\alpha_{14},
%\alpha_{21},\alpha_{15},\alpha_{23},
%\alpha_{8},\alpha_{10},\linebreak
%\alpha_{9},\alpha_{13}\},\
%\{\alpha_{5},\alpha_{16}\},\
%\{\alpha_{12},\alpha_{24},\alpha_{18}\}$;
$\{\alpha_{1},\alpha_{19},\alpha_{6},\alpha_{7},\alpha_{4},
\alpha_{17},\alpha_{20},\alpha_{11},\alpha_{14},
\alpha_{21},\alpha_{15},\alpha_{23},\alpha_{8},\alpha_{10},
\alpha_{9},\alpha_{13}\}$.
%gap> H1:=Image(F[1]);
%Group([ (4,19,7)(6,13,21)(8,11,23)(9,14,15)(10,20,17)(12,18,24), (1,4,7,19)(5,16)(6,10,11,13)(8,14,20,21)(9,15,23,17)
%(12,18) ])
%gap> Orbits(H1);
%[ [ 1, 4, 19, 7 ], [ 5, 16 ], [ 6, 13, 10, 21, 20, 11, 8, 17, 23, 14, 9, 15 ], [ 12, 18, 24 ] ]
The $G_1\cong \SSS_4$ is marked by
$$
G=H_{65,1}\supset G_1=
$$
$$
[(\alpha_{4}\alpha_{19}\alpha_{7})(\alpha_{6}\alpha_{13}\alpha_{21})(\alpha_{8}\alpha_{11}\alpha_{23})
(\alpha_{9}\alpha_{14}\alpha_{15})(\alpha_{10}\alpha_{20}\alpha_{17})(\alpha_{12}\alpha_{18}\alpha_{24}),
$$
$$
(\alpha_{1}\alpha_{4}\alpha_{7}\alpha_{19})(\alpha_{5}\alpha_{16})
(\alpha_{6}\alpha_{10}\alpha_{11}\alpha_{13})(\alpha_{8}\alpha_{14}\alpha_{20}\alpha_{21})
(\alpha_{9}\alpha_{15}\alpha_{23}\alpha_{17})(\alpha_{12}\alpha_{28})]
$$
with suborbits
$\{\alpha_{1},\alpha_{4},\alpha_{19},\alpha_{7}\}$,
$\{\alpha_{6},\alpha_{13},\alpha_{10},\alpha_{21},\alpha_{20},\alpha_{11},\alpha_{8},\alpha_{17},
\alpha_{23},\alpha_{14},\alpha_{9},\alpha_{15}\}$.
Both $G$ and $G_1$ are marked by $\ast$ in Tables 1---4.

\medskip

Case 76:
$({\bf n}=34,\ ((6\aaa_1)_{I},(6\aaa_1)_{II})\subset 12\aaa_1)\Longleftarrow ({\bf n}=65,\ 12\aaa_1)$.

Similar to Cases 1 and 24. By \cite{Nik10}, the $G\cong 2^4D_6$ is marked by
$N_{23}$ and $G=H_{65,3}$ of Case 24 with the orbit
$\{\alpha_{1},\alpha_{7},\alpha_{23},\alpha_{8},\alpha_{11},\alpha_{17},\alpha_{20},
\alpha_{9}, \alpha_{6},\alpha_{4},\alpha_{14},\alpha_{13}\}$.
%gap> H1:=Image(F[1]);
%Group([ (1,4,6)(5,16,22)(7,8,13)(9,23,17)(10,21,19)(11,14,20), (3,5,22,16)(4,6,17,9)(7,13,20,14)(8,11)(10,15,19,21)
%(12,24) ])
%gap> Orbits(H1);
%[ [ 1, 4, 6, 17, 9, 23 ], [ 3, 5, 16, 22 ], [ 7, 8, 13, 11, 20, 14 ], [ 10, 21, 15, 19 ], [ 12, 24 ] ]
The $G_1\cong \SSS_4$ is marked by
$$
G=H_{65,3}\supset G_1=
$$
$$
[(\alpha_{1}\alpha_{4}\alpha_{6})(\alpha_{5}\alpha_{16}\alpha_{22})(\alpha_{7}\alpha_{8}\alpha_{13})
(\alpha_{9}\alpha_{23}\alpha_{17})(\alpha_{10}\alpha_{21}\alpha_{19})(\alpha_{11}\alpha_{14}\alpha_{20}),
$$
$$
(\alpha_{3}\alpha_{5}\alpha_{22}\alpha_{16})(\alpha_{4}\alpha_{6}\alpha_{17}\alpha_{9})
(\alpha_{7}\alpha_{13}\alpha_{20}\alpha_{14})(\alpha_{8}\alpha_{11})
(\alpha_{10}\alpha_{15}\alpha_{19}\alpha_{21})(\alpha_{12}\alpha_{24})]
$$
with suborbits
$\{\alpha_{7},\alpha_{8},\alpha_{13},\alpha_{11},\alpha_{20},\alpha_{14}\}$ of
the type $I$, and
$\{\alpha_{1},\alpha_{4},\alpha_{6},\alpha_{17},\alpha_{9},\alpha_{23}\}$
of the type $II$.
Both $G$ and $G_1$ are marked by $\ast$ in Tables 1---4.

\medskip

Case 77:
$({\bf n}=39,\ (4\aaa_1,4\aaa_1)\subset 8\aaa_1)\Longleftarrow ({\bf n}=56,\ 8\aaa_1)$.

Similar to Cases 1 and 23. By \cite{Nik10}, the $G\cong \Gamma_{15}a_1$
is marked by $N_{23}$ and $G=H_{56,1}$ of Case 23 with the orbit
$\{\alpha_{1}, \alpha_{14}, \alpha_{24}, \alpha_{17}, \alpha_{2},\alpha_{23}, \alpha_{5},
\alpha_{3}\}$.
%gap> H1:=Image(F[1]);
%Group([ (2,17)(6,16)(7,8)(9,13)(10,18)(12,22)(19,20)(23,24), (1,17)(2,3)(6,9)(7,10)(8,13)(15,20)(16,18)(19,21), (5,23)
%(6,10)(7,9)(8,16)(13,18)(14,24)(15,20)(19,21) ])
%gap> Orbits(H1);
%[ [ 1, 17, 2, 3 ], [ 5, 23, 24, 14 ], [ 6, 16, 9, 10, 18, 8, 13, 7 ], [ 12, 22 ], [ 15, 20, 19, 21 ] ]
The $G_1\cong 2^4C_2$ is marked by
$$
G=H_{56,1}\supset G_1=
$$
$$
[(\alpha_{2}\alpha_{17})(\alpha_{6}\alpha_{16})(\alpha_{7}\alpha_{8})(\alpha_{9}\alpha_{13})
(\alpha_{10}\alpha_{18})(\alpha_{12}\alpha_{22})(\alpha_{19}\alpha_{20})(\alpha_{23}\alpha_{24}),
$$
$$
(\alpha_{1}\alpha_{17})(\alpha_{2}\alpha_{3})(\alpha_{6}\alpha_{9})(\alpha_{7}\alpha_{10})
(\alpha_{8}\alpha_{13})(\alpha_{15}\alpha_{20})(\alpha_{16}\alpha_{18})(\alpha_{19}\alpha_{21}),
$$
$$
(\alpha_{5}\alpha_{23})(\alpha_{6}\alpha_{10})(\alpha_{7}\alpha_{9})(\alpha_{8}\alpha_{16})
(\alpha_{13}\alpha_{18})(\alpha_{14}\alpha_{24})(\alpha_{15}\alpha_{20})(\alpha_{19}\alpha_{21})]
$$
with suborbits
$\{\alpha_{1},\alpha_{17},\alpha_{2},\alpha_{3}\}$, $\{\alpha_{5},\alpha_{23},\alpha_{24},\alpha_{14}\}$.
Both $G$ and $G_1$ are marked by $\ast$ in Tables 1---4.

\medskip

Case 78:
$({\bf n}=39,\ (4\aaa_1,8\aaa_1)\subset 12\aaa_1)\Longleftarrow ({\bf n}=65,\ 12\aaa_1)$.

Similar to Cases 1 and 24. By \cite{Nik10},
the $G\cong 2^4D_6$ is marked by by $N_{23}$ and
$G=H_{65,3}$ of Case 24 with the orbit
$\{\alpha_{1},\alpha_{7},\alpha_{23},\alpha_{8},\alpha_{11},\alpha_{17},\alpha_{20},
\alpha_{9}, \alpha_{6},\alpha_{4},\alpha_{14},\alpha_{13}\}$.
%gap> H1:=Image(F[1]);
%Group([ (4,7)(5,16)(6,14)(8,11)(9,13)(12,24)(15,21)(17,20), (1,8)(3,16)(5,22)(6,7)(9,20)(10,19)(11,23)(15,21), (3,5)
%(4,13)(6,7)(9,20)(10,15)(14,17)(16,22)(19,21) ])
%gap> Orbits(H1);
%[ [ 1, 8, 11, 23 ], [ 3, 16, 5, 22 ], [ 4, 7, 13, 6, 9, 14, 20, 17 ], [ 10, 19, 15, 21 ], [ 12, 24 ] ]
The $G_1\cong 2^4C_2$ is marked by
$$
G=H_{65,3}\supset G_1=
$$
$$
[(\alpha_{4}\alpha_{7})(\alpha_{5}\alpha_{16})(\alpha_{6}\alpha_{14})(\alpha_{8}\alpha_{11})
(\alpha_{9}\alpha_{13})(\alpha_{12}\alpha_{24})(\alpha_{15}\alpha_{21})(\alpha_{17}\alpha_{20}),
$$
$$
(\alpha_{1}\alpha_{8})(\alpha_{3}\alpha_{16})(\alpha_{5}\alpha_{22})(\alpha_{6}\alpha_{7})
(\alpha_{9}\alpha_{20})(\alpha_{10}\alpha_{19})(\alpha_{11}\alpha_{23})(\alpha_{15}\alpha_{21}),
$$
$$
(\alpha_{3}\alpha_{5})(\alpha_{4}\alpha_{13})(\alpha_{6}\alpha_{7})(\alpha_{9}\alpha_{20})
(\alpha_{10}\alpha_{15})(\alpha_{14}\alpha_{17})(\alpha_{16}\alpha_{22})(\alpha_{19}\alpha_{21})]
$$
with suborbits
$\{\alpha_{1},\alpha_{8},\alpha_{11},\alpha_{23}\}$,
$\{\alpha_{4},\alpha_{7},\alpha_{13},\alpha_{6},\alpha_{9},\alpha_{14},\alpha_{20},\alpha_{17}\}$.
Both $G$ and $G_1$ are mar\-ked by $\ast$ in Tables 1---4.

\medskip

Case 79:
$({\bf n}=39,\ (8\aaa_1,8\aaa_1)\subset 16\aaa_1)\Longleftarrow ({\bf n}=75,\ 16\aaa_1)$.

Similar to Cases 1 and 10. By \cite{Nik10}, the $G\cong 4^2\AAA_4$
is marked by $N_{23}$ and $G=H_{75,1}$ of Case 10 with the orbit
$\{\alpha_{1},\alpha_{9},\alpha_{6},\alpha_{7},\alpha_{4},
\alpha_{21},\alpha_{20},
\alpha_{22},\alpha_{5},\alpha_{14},\alpha_{10},\alpha_{23},
\alpha_{3},\alpha_{24},\alpha_{18},\alpha_{16}\}$.
The $G_1\cong 2^4C_2$ is marked by
$$
G=H_{75,1}\supset G_1=
$$
%gap> H2:=Image(F[2]);
%Group([ (1,3)(4,22)(5,10)(6,14)(7,23)(9,24)(16,18)(20,21), (3,24)(4,6)(5,21)(10,14)(11,17)(12,13)(16,23)
%(20,22), (3,16)(4,5)(6,21)(10,20)(11,12)(13,17)(14,22)(23,24) ])
%gap> Orbits(H2);
%[ [ 1, 3, 24, 16, 9, 23, 18, 7 ], [ 4, 22, 6, 5, 20, 14, 21, 10 ], [ 11, 17, 12, 13 ] ]
$$
[(\alpha_{1}\alpha_{3})(\alpha_{4}\alpha_{22})(\alpha_{5}\alpha_{10})(\alpha_{6}\alpha_{14})
(\alpha_{7}\alpha_{23})(\alpha_{9}\alpha_{24})(\alpha_{16}\alpha_{18})(\alpha_{20}\alpha_{21}),
$$
$$
(\alpha_{3}\alpha_{24})(\alpha_{4}\alpha_{6})(\alpha_{5}\alpha_{21})(\alpha_{10}\alpha_{14})
(\alpha_{11}\alpha_{17})(\alpha_{12}\alpha_{13})(\alpha_{16}\alpha_{23})(\alpha_{20}\alpha_{22}),
$$
$$
(\alpha_{3}\alpha_{16})(\alpha_{4}\alpha_{5})(\alpha_{6}\alpha_{21})(\alpha_{10}\alpha_{20})
(\alpha_{11}\alpha_{12})(\alpha_{13}\alpha_{17})(\alpha_{14}\alpha_{22})(\alpha_{23}\alpha_{24})]
$$
with suborbits
$\{\alpha_{1},\alpha_{3},\alpha_{24},\alpha_{16},\alpha_{9},\alpha_{23},\alpha_{18},\alpha_{7}\}$,
$\{\alpha_{4},\alpha_{22},\alpha_{6},\alpha_{5},\alpha_{20},\alpha_{14},\alpha_{21},\alpha_{10}\}$.
Both $G$ and $G_1$ are marked by $\ast$ in Tables 1---4.

\medskip

Case 80:
$({\bf n}=40,\ (8\aaa_1,8\aaa_1)\subset 16\aaa_1)\Longleftarrow ({\bf n}=56,\ 16\aaa_1)$.

Similar to Cases 1 and 8. By \cite{Nik10}, the $G\cong \Gamma_{25}a_1$
is marked by $N_{23}$ and
$G=H_{56,2}$ of Case 8 with the orbit
$\{\alpha_{2}, \alpha_{3}, \alpha_{23},
\alpha_{24}, \alpha_{5},\alpha_{8}, \alpha_{17}, \alpha_{19}, \alpha_{6},
\alpha_{7}, \alpha_{11}, \alpha_{20}, \alpha_{18},
\alpha_{10},\alpha_{16},\alpha_{15}\}$.
The $G_1\cong Q_8\ast Q_8$ is marked by
$$
G=H_{56,2}\supset G_1=
$$
%gap> GeneratorsOfGroup(H1);
%[ (2,16)(3,5)(6,17)(9,13)(10,23)(11,15)(12,22)(18,24), (2,10)(3,20)(5,19)(6,24)(7,15)(8,11)(16,18)(17,23),
%  (3,11)(5,15)(7,20)(8,19)(9,13)(10,24)(12,22)(18,23), (1,14)(3,20)(5,19)(7,11)(8,15)(12,22)(16,17)(18,23) ]
$$
[(\alpha_{2}\alpha_{16})(\alpha_{3}\alpha_{5})(\alpha_{6}\alpha_{17})(\alpha_{9}\alpha_{13})
(\alpha_{10}\alpha_{23})(\alpha_{11}\alpha_{15})(\alpha_{12}\alpha_{22})(\alpha_{18}\alpha_{24}),
$$
$$
(\alpha_{2}\alpha_{10})(\alpha_{3}\alpha_{20})(\alpha_{5}\alpha_{19})(\alpha_{6}\alpha_{24})
(\alpha_{7}\alpha_{15})(\alpha_{8}\alpha_{11})(\alpha_{16}\alpha_{18})(\alpha_{17}\alpha_{23}),
$$
$$
(\alpha_{3}\alpha_{11})(\alpha_{5}\alpha_{15})(\alpha_{7}\alpha_{20})(\alpha_{8}\alpha_{19})
(\alpha_{9}\alpha_{13})(\alpha_{10}\alpha_{24})(\alpha_{12}\alpha_{22})(\alpha_{18}\alpha_{23}),
$$
$$
(\alpha_{1}\alpha_{14})(\alpha_{3}\alpha_{20})(\alpha_{5}\alpha_{19})(\alpha_{7}\alpha_{11})
(\alpha_{8}\alpha_{15})(\alpha_{12}\alpha_{22})(\alpha_{16}\alpha_{17})(\alpha_{18}\alpha_{23})]
$$
with suborbits
%gap> Orbits(H1);
%[ [ 1, 14 ], [ 2, 16, 10, 18, 17, 23, 24, 6 ], [ 3, 5, 20, 11, 19, 15, 7, 8 ], [ 9, 13 ], [ 12, 22 ] ]
$\{\alpha_{2},\alpha_{16},\alpha_{10}\alpha_{18},\alpha_{17},\alpha_{23},\alpha_{24},\alpha_{6}\}$,
$\{\alpha_{3},\alpha_{5},\alpha_{20}\alpha_{11},\alpha_{19},\alpha_{15},\alpha_{7},\alpha_{8}\}$.
Both $G$ and $G_1$ are marked by $\ast$ in Tables 1---4.

\medskip

Case 81:
$({\bf n}=49,\ (4\aaa_1,4\aaa_1)\subset 8\aaa_1)\Longleftarrow ({\bf n}=65,\ 8\aaa_1)$.

Similar to Cases 1 and 52. By \cite{Nik10}, the $G\cong 2^4D_6$
is marked by $N_{23}$ and $G=H_{65,4}$ from Case 52 with the orbit
%G:=Group([(1,21,23)(3,22,12)(4,20,13)(5,18,16)(8,11,10)(9,14,17),
%(1,11,23,8)(4,6,13,20)(5,16)(7,17,9,14)(10,21,19,15)(12,24)]);
$\{\alpha_{1},\alpha_{8},\alpha_{21},\alpha_{23},\alpha_{15},
\alpha_{11},\alpha_{10},\alpha_{19}\}$.
%gap> H1:=Image(F[1]);
%Group([ (4,13,20)(5,16,18)(7,17,9)(8,11,19)(12,24,22)(15,21,23), (1,15)(3,12)(4,20)(6,13)(7,9)(14,17)(21,23)
%(22,24), (3,12)(4,13)(6,20)(7,14)(8,10)(9,17)(11,19)(22,24) ])
%gap> Orbits(H1);
%[ [ 1, 15, 21, 23 ], [ 3, 12, 24, 22 ], [ 4, 13, 20, 6 ], [ 5, 16, 18 ], [ 7, 17, 9, 14 ], [ 8, 11, 10, 19 ] ]
The $G_1\cong 2^4C_3$ is marked by
$$
G=H_{65,4}\supset G_1=
$$
$$
[(\alpha_{4}\alpha_{13}\alpha_{20})(\alpha_{5}\alpha_{16}\alpha_{18})(\alpha_{7}\alpha_{17}\alpha_{9})
(\alpha_{8}\alpha_{11}\alpha_{19})(\alpha_{12}\alpha_{24}\alpha_{22})(\alpha_{15}\alpha_{21}\alpha_{23}),
$$
$$
(\alpha_{1}\alpha_{15})(\alpha_{3}\alpha_{12})(\alpha_{4}\alpha_{20})(\alpha_{6}\alpha_{13})
(\alpha_{7}\alpha_{9})(\alpha_{14}\alpha_{17})(\alpha_{21}\alpha_{23})(\alpha_{22}\alpha_{24}),
$$
$$
(\alpha_{3}\alpha_{12})(\alpha_{4}\alpha_{13})(\alpha_{6}\alpha_{20})(\alpha_{7}\alpha_{14})
(\alpha_{8}\alpha_{10})(\alpha_{9}\alpha_{17})(\alpha_{11}\alpha_{19})(\alpha_{22}\alpha_{24})]
$$
with suborbits
$\{\alpha_{1},\alpha_{15},\alpha_{21},\alpha_{23}\}$, $\{\alpha_{8},\alpha_{11},\alpha_{10},\alpha_{19}\}$.
Both $G$ and $G_1$ are marked by $\ast$ in Tables 1---4.

\medskip

Case 82:
$({\bf n}=49,\ (4\aaa_1,12\aaa_1)\subset 16\aaa_1)\Longleftarrow ({\bf n}=75,\ 16\aaa_1)$.

Similar to Cases 1 and 10. By \cite{Nik10}, the $G\cong 4^2\AAA_4$
is marked by $N_{23}$ and $G=H_{75,1}$ of Case 10 with the orbit
$\{\alpha_{1},\alpha_{9},\alpha_{6},\alpha_{7},\alpha_{4},
\alpha_{21},\alpha_{20},
\alpha_{22},\alpha_{5},\alpha_{14},\alpha_{10},\alpha_{23},
\alpha_{3},\alpha_{24},\alpha_{18},\alpha_{16}\}$.
The $G_1\cong 2^4C_3$ is marked by
$$
G=H_{75,1}\supset G_1=
$$
%gap> H1:=Image(F[1]);
%Group([ (1,3,10)(4,6,21)(7,24,14)(9,16,20)(11,13,17)(18,23,22), (3,24)(4,6)(5,21)(10,14)(11,17)(12,13)(16,23)
%(20,22), (3,16)(4,5)(6,21)(10,20)(11,12)(13,17)(14,22)(23,24) ])
%gap> Orbits(H1);
%[ [ 1, 3, 10, 24, 16, 14, 20, 23, 7, 22, 9, 18 ], [ 4, 6, 5, 21 ], [ 11, 13, 17, 12 ] ]
$$
[(\alpha_{1}\alpha_{3}\alpha_{10})(\alpha_{4}\alpha_{6}\alpha_{21})(\alpha_{7}\alpha_{24}\alpha_{14})
(\alpha_{9}\alpha_{16}\alpha_{20})(\alpha_{11}\alpha_{13}\alpha_{17})(\alpha_{18}\alpha_{23}\alpha_{22}),
$$
$$
(\alpha_{3}\alpha_{24})(\alpha_{4}\alpha_{6})(\alpha_{5}\alpha_{21})(\alpha_{10}\alpha_{14})
(\alpha_{11}\alpha_{17})(\alpha_{12}\alpha_{13})(\alpha_{16}\alpha_{23})(\alpha_{20}\alpha_{22}),
$$
$$
(\alpha_{3}\alpha_{16})(\alpha_{4}\alpha_{5})(\alpha_{6}\alpha_{21})(\alpha_{10}\alpha_{20})
(\alpha_{11}\alpha_{12})(\alpha_{13}\alpha_{17})(\alpha_{14}\alpha_{22})(\alpha_{23}\alpha_{24})]
$$
with suborbits
$\{\alpha_{4},\alpha_{6},\alpha_{5},\alpha_{21}\}$,
$\{\alpha_{1},\alpha_{3},\alpha_{10},\alpha_{24},\alpha_{16},\alpha_{14},\alpha_{20},\alpha_{23},
\alpha_{7},\alpha_{22},\alpha_{9},\alpha_{18}\}$.
Both $G$ and $G_1$ are marked by $\ast$ in Tables 1---4.

\medskip

%%%%%%%%%%%%%%%%%%%%%%%%%%%%%%%%%%%%%%%%%%%%
%%%%%%%%%%%%%%%%%%%%%%%%%%%%%%%%%%%%%%%%%%%%
%%%%%%%%%%%%%%%%%%%%%%%%%%%%%%%%%%%%%%%%%%%%

This completes the proof of the Theorem.

%%%%%%%%%%%%%%%%%%%%%%%%%%%%%%%%%%%%%%%%%%%%%%%
\end{proof}

\medskip
%%%%%%%%%%%%%%%%%%%%%%%%%%%%%%%%%%%%%
%%%%%%%%%%%%%%%%%%%%%%%%%%%%%%%%%%%%%

Thus, we have finally obtained

\vskip1cm

\begin{theorem} Classification of Picard lattices $S=MS_X$ of K3 surfaces $X$ with
finite symplectic automorphic groups which are big enough
(larger than $D_6$, $C_4$, $(C_2)^2$, $C_3$, $C_2$ and $C_1$) and with
at least one $-2$ curve is given in the Tables 1 --- 4 of Section \ref{sec4:tables}
below in lines which are not marked by $o$.
\end{theorem}

We hope to consider remaining symplectic groups $D_6$, $C_4$, $(C_2)^2$, $C_3$, $C_2$, $C_1$
later as well. Now, we obtain for these groups

\begin{theorem} If the Picard lattice $S=MS_X$ of a K3 surface $X$ with
at least one $-2$ curve is different from all lattices of lines of Tables 1 --- 4
of Section \ref{sec4:tables} which are not marked by $o$ (for example, if the genus is different),
then the symplectic automorphism group of $X$
is small, it is one of groups: $D_6$, $C_4$, $(C_2)^2$, $C_3$, $C_2$ or $C_1$.
\end{theorem}

%%%%%%%%%%%%%%%%%%%%%%%%%%%%%%%%%%%%%
%%%%%%%%%%%%%%%%%%%%%%%%%%%%%%%%%%%%%

\section{Tables.}
\label{sec4:tables}

Here we present tables of lattices $S$ of degenerations of K\"ahlerian K3 surfaces
with finite symplectic automorphism groups and classification of Picard lattices
of K3 surfaces which were discussed in previous sections.

\begin{table}
\label{table1}
\caption{Types and lattices $S$ of degenerations of
codimension $1$ of K\"ahlerian K3 surfaces
with finite symplectic automorphism groups $G=Clos(G)$. All lines
are marked by $\ast$, by definition.}

%%\index{Table 1}

%%%\addtocontents{toc}{\contentsline {section}{\tocsection {}{T.1}{Table 1}}
%%%{\pageref{table1}}}

\begin{tabular}{|c||c|c|c|c|c|c|c|c|}
\hline
 {\bf n}& $|G|$& $i$&  $G$  & $\rk S_G$ &$q_{S_G}$  &$Deg$& $\rk S$ &$q_S$ \\
\hline
\hline
 $1$  &$2$ & $1$& $C_2$ &   $8$     &$2_{II}^{+8}$ & $\aaa_1$   &  $9$  &
$2_7^{+9}$\\
\hline
      &    &    &    &              &              & $2\aaa_1$  &  $9$  &
$2_{II}^{-6},4_3^{-1}$\\
\hline
\hline
 $2$  & $3$& $1$&$C_3$  &    $12$   & $3^{+6}$       & $\aaa_1$   &  $13$ &
$2_3^{-1},3^{+6}$\\
\hline
      &    &    &       &           &                & $3\aaa_1$  & $13$  &
$2_1^{+1},3^{-5}$\\
\hline
\hline
 $3$  &$4$ & $2$&$C_2^2$&    $12$   & $2_{II}^{-6},4_{II}^{-2}$& $\aaa_1$   & $13$  &
$2_3^{+7},4_{II}^{+2}$\\
\hline
      &    &    &       &           &                            & $2\aaa_1$& $13$  &
$2_{II}^{-4},4_7^{-3}$\\
\hline
      &    &    &       &           &                & $4\aaa_1$  & $13$  &
$2_{II}^{-6},8_3^{-1}$\\
\hline
\end{tabular}
\end{table}

\begin{table}

\begin{tabular}{|c||c|c|c|c|c|c|c|c|}
\hline
 {\bf n}& $|G|$& $i$&  $G$   & $\rk S_G$ &$q_{S_G}$ &$Deg$& $\rk S$ &$q_S$ \\
\hline
\hline
 $17$&$12$ & $3$&$\AAA_4$&  $16$& $2_{II}^{-2},4_{II}^{-2},3^{+2}$ & $\aaa_1$ &$17$  &
$2_7^{-3},4_{II}^{+2},3^{+2}$\\
\hline
     &     &    &        &           &     &$3\aaa_1$ &  $17$  &
$2_1^{-3},4_{II}^{+2},3^{-1}$\\
\hline
     &     &    &        &           &            &$4\aaa_1$ &  $17$  &
$2_{II}^{-2},8_3^{-1},3^{+2}$\\
\hline
     &     &    &        &           &             &$6\aaa_1$ &  $17$  &
$4_1^{-3},3^{+1}$\\
\hline
     &     &    &        &           &             &$12\aaa_1$&  $17$  &
$2_{II}^{-2},8_1^{+1},3^{-1}$\\
\hline
\hline
 $18$&$12$ &$4$ &$D_{12}$&  $16$&$2_{II}^{+4},3^{+4}$ &  $\aaa_1$&  $17$&
$2_7^{+5},3^{+4}$\\
\hline
     &     &    &        &           &   &  $2\aaa_1$    &  $17$  &
$2_{II}^{+2},4_7^{+1},3^{+4}$\\
\hline

\end{tabular}

\end{table}

\begin{table}

\begin{tabular}{|c||c|c|c|c|c|c|c|c|}
\hline
 {\bf n}& $|G|$& $i$&  $G$   & $\rk S_G$ &$q_{S_G}$       &$Deg$& $\rk S$ &$q_S$    \\
\hline
\hline
 $39$&$32$ &$27$&$2^4C_2$&  $17$     &$2_{II}^{+2},4_0^{+2},8_7^{+1}$ & $4\aaa_1$&$18$&
$4_6^{+4}$\\
\hline
     &     &    &        &           &                                & $8\aaa_1$&  $18$&
$2_{II}^{+2},4_7^{+1},8_7^{+1}$\\
\hline
     &     &    &        &           &                                & $16\aaa_1$&$18$&
$2_{II}^{+2},4_6^{+2}$\\
\hline
\hline
$40$ &$32$ &$49$&$Q_8*Q_8$& $17$     & $4_7^{+5}$ & $8\aaa_1$  & $18$ &
$4_6^{+4}$\\
\hline
\hline
$46$ &$36$ &$9$ &$3^2C_4$ & $18$     & $2_6^{-2},3^{+2},9^{-1}$& $6\aaa_1$   &  $19$  &
$4_7^{+1},3^{+1},9^{-1}$               \\
\hline
     &     &    &        &           &                           & $9\aaa_1$ &  $19$  &
$2_5^{-3},3^{+2}$ \\
\hline
     &     &    &        &           &                           & $9\aaa_2$ &  $19$  &
$2_5^{-1},3^{+2}$\\
\hline

\end{tabular}

\end{table}

%%%%%%%%%%%%%%%%%%%%%%%%%%%%%%%%%%%%%%%%%%
%%%%%%%%%%%%%%%%%%%%%%%%%%%%%%%%%%%%%%%%%%
%%%%%%%%%%%%%%%%%%%%%%%%%%%%%%%%%%%%%%%%%%
%%%%%%%%%%%%%%%%%%%%%%%%%%%%%%%%%%%%%%%%%%

\begin{table}
\label{table2}
\caption{Types and lattices $S$ of degenerations of
codimension $\ge 2$ of K\"ahlerian K3 surfaces
with finite symplectic automorphism groups
$G=Clos(G)$ for ${\bf n}\ge 12$.}

\begin{tabular}{|c||c|c|c|c|c|c|c|c|}
\hline
 {\bf n}& $|G|$& $i$&  $G$   & $\rk S_G$ &$Deg$& $\rk S$ &$q_S$            \\
\hline
\hline
  $12$&$8$& $4$ & $Q_8$ &  $17$     &
  $(8\aaa_1,8\aaa_1)\subset 16\aaa_1$ $o$ & $19$  &
$2_{II}^{-2},8_5^{-1}$ $\ast$\\
\hline
      &   &     &       &           &
$(8\aaa_1,8\aaa_1)\subset 8\aaa_2$ & $19$   &
$2_7^{-3},3^{-1}$ $\ast$\\
\hline
      &   &     &       &           &
$(8\aaa_1,\aaa_2)\subset 8\aaa_1\amalg \aaa_2$ & $19$   &
$2_2^{+2},16_3^{-1}$ $\ast$\\
\hline
      &   &     &       &           &
$(\aaa_2,\aaa_2)\subset 2\aaa_2$ $o$ & $19$   &
$2_5^{-1},8_{II}^{-2}$ $\ast$\\
\hline
\hline
 $16$ &$10$&$1$ &$D_{10}$&  $16$   & $(\aaa_1,\aaa_1)\subset 2\aaa_1$  & $18$  &
$2_6^{+2},5^{+4}$\\
\hline
      &   &     &       &            &  $(\aaa_1,5\aaa_1)\subset 6\aaa_1$  &$18$  &
$2_6^{+2},5^{-3}$\\
\hline
      &   &     &       &            &  $(5\aaa_1,5\aaa_1)\subset 10\aaa_1$  &$18$  &
$2_6^{+2},5^{+2}$ $\ast$\\
\hline
\end{tabular}
\end{table}

\begin{table}

\begin{tabular}{|c||c|c|c|c|c|c|c|c|}
\hline
 &   &   &  & &  $(6\aaa_1,6\aaa_1)\subset 6\aaa_2$\ $o$ &  $18$  &
$4_{II}^{-2},3^{+1}$  $\ast$ \\
\hline
 &   &   &  & &  $(6\aaa_1,12\aaa_1)\subset 6\aaa_3$ &  $18$  &
$4_6^{+2}$  $\ast$ \\
\hline
&   &   &  & &  $(\aaa_1,\aaa_1,\aaa_1)\subset 3\aaa_1$ $o$ &  $19$  &
$2_{1}^{+1},4_{II}^{+2},3^{+2}$  $\ast$ \\
\hline
&   &   &  & &  $(\aaa_1,\aaa_1,4\aaa_1)\subset 6\aaa_1$\ $o$ &  $19$  &
$2_{2}^{+2},8_{7}^{+1},3^{+2}$  $\ast$ \\
\hline
&   &   &  & &  $(\aaa_1,\aaa_1,6\aaa_1)\subset 8\aaa_1$\ $o$ &  $19$  &
$4_{7}^{-3},3^{+1}$  $\ast$ \\
\hline
&   &   &  & &  $(\aaa_1,\aaa_1,12\aaa_1)\subset 14\aaa_1$ \ $o$ &  $19$  &
$2_{2}^{+2},8_{5}^{-1},3^{-1}$  $\ast$ \\
\hline
&   &   &  & &  $(\aaa_1,3\aaa_1,4\aaa_1)\subset 8\aaa_1$ $o$ &  $19$ &
$2_{II}^{-2},8_{7}^{+1},3^{-1}$  $\ast$ \\
\hline
&   &   &  & &  $(\aaa_1,3\aaa_1,12\aaa_1)\subset 16\aaa_1$ $o$ &  $19$  &
$2_{II}^{-2},8_{5}^{-1}$  $\ast$ \\
\hline

\end{tabular}
\end{table}

\begin{table}

\begin{tabular}{|c||c|c|c|c|c|c|}

\hline
 {\bf n}& $|G|$& $i$&  $G$   &$Deg$& $\rk S$ & $q_S$            \\
\hline
\hline
 $18$&$12$ &$4$ &$D_{12}$&   $(\aaa_1,2\aaa_1)\subset 3\aaa_1$ &  $18$  &
$2_7^{+3},4_7^{+1},3^{+4}$\\
\hline
 &   &   &  &   $(\aaa_1,3\aaa_1)\subset 4\aaa_1$ &  $18$  &
$2_0^{+4},3^{-3}$ $\ast$\\
\hline
&   &    &  &    $(\aaa_1,6\aaa_1)\subset 7\aaa_1$ &  $18$  &
$2_1^{-3},4_7^{+1},3^{+3}$ $\ast$\\
\hline
&   &    &  &    $(2\aaa_1,3\aaa_1)\subset 5\aaa_1$ &  $18$  &
$2_1^{+3},4_7^{+1},3^{-3}$\\
\hline
&   &    &  &    $\left((2\aaa_1,6\aaa_1)\subset 8\aaa_1\right)_I$ &  $18$  &
$4_4^{-2},3^{+3}$ $\ast$ \\
\hline
&   &    &  &    $\left((2\aaa_1,6\aaa_1)\subset 8\aaa_1\right)_{II}$ &  $18$  &
$2_{II}^{-4},3^{+3}$ \\
\hline
&   &    &  &    $(2\aaa_1,6\aaa_1)\subset 2\ddd_4$ &  $18$  &
$2_{II}^{-2},3^{+3}$ $\ast$\\
\hline
\end{tabular}
\end{table}

\begin{table}
\begin{tabular}{|c||c|c|c|c|c|c|c|c|}

\hline
&   &    &  &
$
\left(\begin{array}{rrr}
2\aaa_1 & 5\aaa_1 & (8\aaa_1)_{II}  \\
         & 3\aaa_1 & 3\aaa_3   \\
         &         & 6\aaa_1
\end{array}\right)
\subset 2\aaa_1\amalg 3\aaa_3$
&  $19$  & $2_{II}^{+2},4_1^{+1},3^{+2}$\\
\hline
&   &    &  &
$
\left(\begin{array}{rrr}
2\aaa_1 & 5\aaa_1 & 2\ddd_4  \\
         & 3\aaa_1 & 9\aaa_1   \\
         &         & 6\aaa_1
\end{array}\right)
\subset 3\aaa_1\amalg 2\ddd_4$
&  $19$  & $2_5^{+3},3^{-2}$ $\ast$\\
\hline
&   &    &  &
$
\left(\begin{array}{rrr}
2\aaa_1 & (8\aaa_1)_I & (8\aaa_1)_{II}  \\
         & 6\aaa_1 & 12\aaa_1   \\
         &         & 6\aaa_1
\end{array}\right)
\subset 14\aaa_1$
&  $19$  & $2_{II}^{+2},4_1^{+1},3^{+2}$\\
\hline
&   &    &  &
$
\left(\begin{array}{rrr}
2\aaa_1 & (8\aaa_1)_I & (8\aaa_1)_I  \\
         & 6\aaa_1 & 6\aaa_2   \\
         &         & 6\aaa_1
\end{array}\right)
\subset 2\aaa_1\amalg 6\aaa_2$
&  $19$  & $4_3^{-1},3^{-1},9^{-1}$ $\ast$\\
\hline
&   &    &  &
$
\left(\begin{array}{rrr}
2\aaa_1 & 2\ddd_4 & (8\aaa_1)_I  \\
         & 6\aaa_1 & 12\aaa_1   \\
         &         & 6\aaa_1
\end{array}\right)
\subset 2\ddd_4\amalg 6\aaa_1$
&  $19$  & $4_1^{+1},3^{+2}$ $\ast$\\
\hline
&   &    &  &
$(3\aaa_1,6\aaa_1,6\aaa_1)\subset 15 \aaa_1$
&  $19$  & $2_3^{-1},4_4^{-2},3^{-1}$  $\ast$\\
\hline
&   &    &  &
$(3\aaa_1,6\aaa_1,6\aaa_1)\subset 3\aaa_1\amalg 6\aaa_2 $
&  $19$  & $2_1^{-3},9^{-1}$  $\ast$\\
\hline
&   &    &  &
$
\left(\begin{array}{rrr}
3\aaa_1 & 3\aaa_3 & 9\aaa_1  \\
         & 6\aaa_1 & 12\aaa_1   \\
         &         & 6\aaa_1
\end{array}\right)
\subset 3\aaa_3\amalg 6\aaa_1$
&  $19$  & $2_6^{-2},4_5^{-1},3^{+1}$ $\ast$\\
\hline

\end{tabular}
\end{table}

\begin{table}

\begin{tabular}{|c||c|c|c|c|c|c|c|c|}
\hline
 {\bf n}& $|G|$& $i$&  $G$   & $Deg$& $\rk S$ &$q_S$            \\
\hline
\hline
 $22$& $16$&$11$&$C_2\times D_8$&   $(2\aaa_1,2\aaa_1)\subset 4\aaa_1$ \ $o$
 &  $18$  & $4_6^{+4}$ $\ast$\\
\hline
     &     &    &        &              $(2\aaa_1,4\aaa_1)\subset 6\aaa_1$
     &  $18$  &
$4_7^{+3},8_7^{+1}$ $\ast$ \\
\hline
     &     &    &        &    $(2\aaa_1,8\aaa_1)\subset 10\aaa_1$    &  $18$  &
$4_6^{+4}$\\
\hline
     &     &    &        &    $((4\aaa_1,4\aaa_1)\subset 8\aaa_1)_I$\ $o$   &  $18$  &
$4_6^{+4}$ $\ast$\\
\hline
     &     &    &        &   $((4\aaa_1,4\aaa_1)\subset 8\aaa_1)_{II}$    &  $18$  &
$2_{II}^{-2},8_6^{-2}$ $\ast$\\
\hline
     &     &    &        &      $(4\aaa_1,8\aaa_1)\subset 12\aaa_1$    &  $18$  &
$2_{II}^{+2},4_7^{+1},8_7^{+1}$ $\ast$\\
\hline
     &     &    &        &           $(4\aaa_1,8\aaa_1)\subset 4\aaa_3$    &  $18$  &
$2_{II}^{-2},4_2^{-2}$ $\ast$\\
\hline

\end{tabular}

\end{table}

\begin{table}

\begin{tabular}{|c||c|c|c|c|c|c|c|c|}
\hline
 {\bf n}& $|G|$& $i$&  $G$ & $\rk S_G$  & $Deg$& $\rk S$ &$q_S$   \\
\hline
\hline
 $34$& $24$&$12$&$\SSS_4$    & $17$ &
$(\aaa_1,\aaa_1)\subset 2\aaa_1$ $o$ &  $19$  &
$4_1^{+3},3^{+2}$ $\ast$\\
\hline
& & & &  &  $(\aaa_1,2\aaa_1)\subset 3\aaa_1$ $o$ &  $19$  &
$2_1^{+1},4_{II}^{+2},3^{+2}$ $\ast$\\
\hline
& & & & & $(\aaa_1,3\aaa_1)\subset 4\aaa_1$ $o$ &  $19$  &
$4_{3}^{-3},3^{-1}$ $\ast$\\
\hline
& & & & & $(\aaa_1,4\aaa_1)\subset 5\aaa_1$&  $19$  &
$2_1^{+1},4_1^{+1},8_7^{+1},3^{+2}$ $\ast$\\
\hline
& & & & & $(\aaa_1,(6\aaa_1)_{I})\subset 7\aaa_1$&  $19$  &
$2_1^{+1},4_{6}^{-2},3^{+1}$ $\ast$\\
\hline
& & & & & $(\aaa_1,8\aaa_1)\subset 9\aaa_1$&  $19$  &
$2_1^{+1},4_{0}^{+2},3^{+2}$ $\ast$\\
\hline
& & & & & $(\aaa_1,12\aaa_1)\subset 13\aaa_1$&  $19$  &
$2_7^{+1},4_{5}^{-1},8_7^{+1},3^{-1}$\\
\hline
\end{tabular}

\end{table}

\vfill

%%%%%%%%%%%%%%%%%%%%%%%%%%%%%%%%%%%%%%%%%%%
%%%%%%%%%%%%%%%%%%%%%%%%%%%%%%%%%%%%%%%%%%%
%%%%%%%%%%%%%%%%%%%%%%%%%%%%%%%%%%%%%%%%%%%%%
%%%%%%%%%%%%%%%%%%%%%%%%%%%%%%%%%%%%%%%%%%%%
%%%%%%%%%%%%%%%%%%%%%%%%%%%%%%%%%%%%%%%%%%%%%%

\begin{table}
\label{table3}
\caption{Types and lattices $S$ of degenerations of
codimension $\ge 2$ of K\"ahlerian K3 surfaces
with symplectic automorphism group $D_8$.}

\begin{tabular}{|c||c|c|c|c|c|c|c|c|}
\hline
 {\bf n}& $|G|$& $i$&  $G$   & $\rk S_G$ &$Deg$& $\rk S$ &$q_S$            \\
\hline
\hline
  $10$&$8$& $3$ & $D_8$ &  $15$     &
  $(\aaa_1,\aaa_1)\subset 2\aaa_1$\ $o$ & $17$  &
$4_7^{+5}$ $\ast$\\
\hline
      &   &     &       &           &
$(\aaa_1,(2\aaa_1)_I)\subset 3\aaa_1$ & $17$   &
$2_7^{+1},4_0^{+4}$ \\
\hline
      &   &     &       &           &
$(\aaa_1,(2\aaa_1)_{II})\subset \aaa_3$ & $17$   &
$4_7^{+5}$ $\ast$\\
\hline
      &   &     &       &           &
$(\aaa_1,4\aaa_1)\subset 5\aaa_1$ & $17$   &
$2_7^{+1},4_1^{+3},8_7^{+1}$\\
\hline
      &   &     &       &           &
$(\aaa_1,8\aaa_1)\subset 9\aaa_1$ & $17$   &
$2_1^{+1},4_6^{+4}$\\
\hline
      &   &     &       &           &
$(\aaa_1,2\aaa_2)\subset \aaa_1\amalg 2\aaa_2$ & $17$   &
$2_7^{+1},4_{II}^{+4}$ $\ast$\\
\hline
      &   &     &       &           &
$((2\aaa_1)_I,(2\aaa_1)_{II})\subset 4\aaa_1$ & $17$   &
$4_7^{+5}$ $\ast$\\
\hline
\end{tabular}
\end{table}

\begin{table}
\begin{tabular}{|c||c|c|c|c|c|c|c|c|}
\hline
 {\bf n}& $G$   & $\rk S_G$ &$Deg$& $\rk S$ &$q_S$            \\
\hline
\hline
  $10$ & $D_8$ &  $15$     & $(\aaa_1,4\aaa_1,8\aaa_1)\subset 13\aaa_1$ & $18$  &
$2_7^{+1},4_0^{+2},8_7^{+1}$  \\
\hline
  &    &   & $(\aaa_1,4\aaa_1,8\aaa_1)\subset \aaa_1\amalg 4\aaa_3$ & $18$  &
$2_7^{+1},4_7^{+3}$  \\
\hline
  &    &   & $(\aaa_1,4\aaa_1,2\aaa_2)\subset 5\aaa_1\amalg 2\aaa_2$ & $18$  &
$2_3^{-1},4_{II}^{+2},8_3^{-1}$  \\
\hline
&   &     &
$
\left(\begin{array}{ccc}
(2\aaa_1)_{I} & 4\aaa_1 & (6\aaa_1)_I  \\
         & (2\aaa_1)_{II} & 6\aaa_1   \\
         &         & 4\aaa_1
\end{array}\right)
\subset 8\aaa_1
$ & $18$   & $4_5^{+3},8_1^{+1}$ \\
\hline
&   &     &
$
\left(\begin{array}{ccc}
(2\aaa_1)_{I} & 4\aaa_1 & (6\aaa_1)_{II}  \\
         & (2\aaa_1)_{II} & 6\aaa_1   \\
         &         & 4\aaa_1
\end{array}\right)
\subset 8\aaa_1
$ & $18$   & $4_7^{+3},8_7^{+1}$ $\ast$  \\
\hline
&   &     &
$
\left(\begin{array}{ccc}
(2\aaa_1)_{II} & 4\aaa_1 & 6\aaa_1  \\
         & (2\aaa_1)_{I} & 2\aaa_3   \\
         &         & 4\aaa_1
\end{array}\right)
\subset 2\aaa_1\amalg 2\aaa_3
$ & $18$ & $4_6^{+4}$  $\ast$ \\
\hline
&   &     &
$
\left(\begin{array}{ccc}
(2\aaa_1)_{I} & (6\aaa_1)_I & (6\aaa_1)_{II}  \\
         & 4\aaa_1 & (8\aaa_1)_I   \\
         &         & 4\aaa_1
\end{array}\right)
\subset 10\aaa_1
$ & $18$   & $2_6^{+2},8_0^{+2}$ \\
\hline
\end{tabular}
\end{table}

\begin{table}
\begin{tabular}{|c||c|c|c|c|c|c|c|c|}
\hline
 {\bf n}& $G$   & $\rk S_G$ &$Deg$& $\rk S$ &$q_S$            \\
\hline
\hline
  $10$ & $D_8$ &  $15$     &
$
\left(\begin{array}{ccc}
(2\aaa_1)_{I} & (6\aaa_1)_I & 10\aaa_1  \\
         & 4\aaa_1 & 12\aaa_1   \\
         &         & 8\aaa_1
\end{array}\right)
\subset 14\aaa_1
$ & $18$   & $2_2^{-2},4_1^{+1},8_3^{-1}$ $\ast$ \\
\hline
&   &     &
$
\left(\begin{array}{ccc}
(2\aaa_1)_{I} & (6\aaa_1)_I & 10\aaa_1  \\
         & 4\aaa_1 & 4\aaa_3   \\
         &         & 8\aaa_1
\end{array}\right)
\subset 2\aaa_1\amalg 4\aaa_3
$ & $18$   & $2_2^{-2},4_0^{+2}$  \\
\hline
&   &     &
$
\left(\begin{array}{ccc}
(2\aaa_1)_{I} & (6\aaa_1)_{II} & 10\aaa_1  \\
         & 4\aaa_1 & 4\aaa_3   \\
         &         & 8\aaa_1
\end{array}\right)
\subset 2\aaa_1\amalg 4\aaa_3
$ & $18$ & $2_2^{-2},4_{II}^{-2}$  $\ast$ \\
\hline
&   &     &
$((2\aaa_1)_{I},4\aaa_1,8\aaa_1)\subset 2\aaa_3\amalg 8\aaa_1$
& $18$   & $2_4^{-2},4_6^{-2}$ $\ast$ \\
\hline
&   &     &
$
\left(\begin{array}{ccc}
4\aaa_1  & (8\aaa_1)_I & (8\aaa_1)_{II}  \\
         & 4\aaa_1   & (8\aaa_1)_{I} \\
         &         & 4\aaa_1
\end{array}\right)
\subset 12\aaa_1$\ $o$
& $18$   & $2_{II}^{-2},4_1^{+1},8_5^{-1}$ $\ast$ \\
\hline
&   &     &
$
\left(\begin{array}{ccc}
4\aaa_1  & (8\aaa_1)_I & (8\aaa_1)_{I}  \\
         & 4\aaa_1       & 4\aaa_2   \\
         &               & 4\aaa_1
\end{array}\right)
\subset 4\aaa_1\amalg 4\aaa_2
$ & $18$   & $4_{7}^{+1},8_5^{-1},3^{+1}$ $\ast$ \\
\hline
&   &     &
$
\left(\begin{array}{ccc}
4\aaa_1  & 4\aaa_2 & (8\aaa_1)_{II}   \\
         & 4\aaa_1 & 4\aaa_2   \\
         &         & 4\aaa_1
\end{array}\right)
\subset 4\aaa_3
$\ $o$ & $18$   & $2_{II}^{-2},4_2^{-2}$ $\ast$ \\
\hline

\end{tabular}
\end{table}

%%%%%%%%%%%%%%%%%%%%%%%%%%%%%%%%%%%%%%%%%%%%%%%%%%%%%%%%%%%%%%%%%%%
%%%%%%%%%%%%%%%%%%%%%%%%%%%%%%%%%%%%%%%%%%%%%%%%%%%%%%%%%%%%%%%%%%
%%%%%%%%%%%%%%%%%%%%%%%%%%%%%%%%%%%%%%%%%%%%%%%%%%%%%%%%%%%%%%%%%%%%%%
%%%%%%%%%%%%%%%%%%% codim=4
%%%%%%%%%%%%%%%%%%%%%%%%%%%%%%%%%%%%%%%%%%%%%%%%%%%%%%%%%%%%%%%%%%%%%

\begin{table}
\begin{tabular}{|c||c|c|c|c|c|c|c|c|}
\hline
 {\bf n}& $G$  &$Deg$& $\rk S$ &$q_S$            \\
\hline
\hline
  $10$ & $D_8$ &
$
\left(\begin{array}{cccc}
\aaa_1 & 2\aaa_1 & 3\aaa_1      &  5\aaa_1       \\
       & \aaa_1  & 3\aaa_1      &  5\aaa_1       \\
       &         & (2\aaa_1)_I  & (6\aaa_1)_I   \\
       &         &              &  4\aaa_1
\end{array}\right)
\subset 8\aaa_1
$ $o$ & $19$   & $4_4^{-2},8_5^{-1}$ $\ast$ \\
\hline
  &   &
$ (\aaa_1,\,\aaa_1,\,(2\aaa_1)_I,\, 8\aaa_1)\subset 12\aaa_1$ $o$
& $19$   & $4_5^{+3}$ $\ast$ \\
\hline
  &   &
$
\left(\begin{array}{cccc}
\aaa_1 & 2\aaa_1 & 5\aaa_1      &  5\aaa_1       \\
       & \aaa_1  & 5\aaa_1      &  5\aaa_1       \\
       &         & 4\aaa_1      & (8\aaa_1)_I    \\
       &         &              &  4\aaa_1
\end{array}\right)
\subset 10\aaa_1
$\ $o$ & $19$   & $4_5^{-1},8_4^{-2}$ $\ast$ \\
\hline
  &   &
$ (\aaa_1,\,\aaa_1,\,4\aaa_1,\,8\aaa_1)\subset 14\aaa_1$ $o$
& $19$   & $4_4^{-2},8_5^{-1}$ \\
\hline
  &   &
$ (\aaa_1,\,\aaa_1,\,4\aaa_1,\, 8\aaa_1)\subset 2\aaa_1\amalg 4\aaa_3$
\ $o$ & $19$   & $4_5^{+3}$ $\ast$ \\
\hline
  &   &
$
\left(\begin{array}{cccc}
\aaa_1 & \aaa_3 & 3\aaa_1      &  5\aaa_1       \\
       & (2\aaa_1)_{II}  & 4\aaa_1      &  8\aaa_1       \\
       &         & (2\aaa_1)_I      & (6\aaa_1)_I    \\
       &         &              &  4\aaa_1
\end{array}\right)
\subset \aaa_3\amalg 6\aaa_1
$ & $19$   & $4_4^{-2},8_5^{-1}$ \\
\hline
  &   &
$
\left(\begin{array}{cccc}
\aaa_1 & \aaa_3              & 3\aaa_1          &  5\aaa_1     \\
       & (2\aaa_1)_{II}      &  4\aaa_1         &  6\aaa_1     \\
       &                     & (2\aaa_1)_I      &  2\aaa_3     \\
       &                     &                  &  4\aaa_1
\end{array}\right)
\subset 3\aaa_3
$\ $o$ & $19$   & $4_5^{+3}$ $\ast$ \\
\hline
\end{tabular}
\end{table}

\begin{table}
\begin{tabular}{|c||c|c|c|c|c|c|c|c|}
\hline
 {\bf n}& $G$  &$Deg$& $\rk S$ &$q_S$            \\
\hline
\hline
  $10$ & $D_8$ &
$ (\aaa_1,\,(2\aaa_1)_I,\,4\aaa_1,\, 8\aaa_1)
\subset \aaa_1\amalg 2\aaa_3,\amalg 8\aaa_1$
& $19$   & $2_3^{-1},4_6^{-2}$ $\ast$ \\
\hline
&   &
$
\left(\begin{array}{cccc}
\aaa_1 & 3\aaa_1      & 5\aaa_1            &  9\aaa_1       \\
       & (2\aaa_1)_I  & (6\aaa_1)_{I}      & 10\aaa_1       \\
       &              & 4\aaa_1            & 4\aaa_3        \\
       &         &                         & 8\aaa_1
\end{array}\right)
\subset 3\aaa_1\amalg 4\aaa_3
$
& $19$   & $2_5^{-1},4_4^{-2}$  \\
\hline
  &   &
$
\left(\begin{array}{cccc}
\aaa_1 & 5\aaa_1   & 5\aaa_1      &  5\aaa_1       \\
       & 4\aaa_1  & (8\aaa_1)_I  &  (8\aaa_1)_I       \\
       &          & 4\aaa_1      &  4\aaa_2    \\
       &         &               &  4\aaa_1
\end{array}\right)
\subset 5\aaa_1\amalg 4\aaa_2
$ & $19$   & $2_1^{+1},4_7^{+1},8_3^{-1},3^{+1}$ \\
\hline
  &   &
$
\left(\begin{array}{cccc}
\aaa_1 & 5\aaa_1   & 5\aaa_1            &  9\aaa_1     \\
       & 4\aaa_1   & (8\aaa_1)_{I}      & 12\aaa_1     \\
       &                     & 4\aaa_1  &  4\aaa_3     \\
       &                     &          &  8\aaa_1
\end{array}\right)
\subset 5\aaa_1\amalg 4\aaa_3
$ & $19$   & $2_7^{+1},4_5^{-1},8_5^{-1}$ $\ast$ \\
\hline
  &   &
$
\left(\begin{array}{cccc}
\aaa_1 & 5\aaa_1             & 5\aaa_1          &  9\aaa_1     \\
       & 4\aaa_1             & 4\aaa_2          &  12\aaa_1    \\
       &                     & 4\aaa_1          &  4\aaa_3     \\
       &                     &                  &  8\aaa_1
\end{array}\right)
\subset \aaa_1\amalg 4\ddd_4
$ & $19$   & $2_7^{+1},4_6^{+2}$   \\
\hline
\end{tabular}
\end{table}

\begin{table}
\begin{tabular}{|c||c|c|c|c|c|c|c|c|}
\hline
 {\bf n}& $G$  &$Deg$& $\rk S$ &$q_S$            \\
\hline
\hline
  $10$ & $D_8$ &
$
\left(\begin{array}{cccc}
(2\aaa_1)_{II}   & 6\aaa_1   & 4\aaa_1        &  6\aaa_1       \\
                 & 4\aaa_1   & (6\aaa_1)_{I}  &  (8\aaa_1)_I   \\
       &                     & (2\aaa_1)_I    &  2\aaa_3        \\
       &                     &                &  4\aaa_1
\end{array}\right)
\subset 6\aaa_1\amalg 2\aaa_3
$
& $19$   & $4_6^{+2},8_7^{+1}$  \\
\hline
  &   &
$
\left(\begin{array}{cccc}
(2\aaa_1)_{II}    & 4\aaa_1        & 6\aaa_1      & 6\aaa_1     \\
                  & (2\aaa_1)_{I}  & 2\aaa_3  &  (6\aaa_1)_{II} \\
       &          & 4\aaa_1        &  4\aaa_2    \\
       &          &                &  4\aaa_1
\end{array}\right)
\subset 2\aaa_1\amalg 2\aaa_5
$ & $19$   & $4_5^{+3}$ $\ast$ \\
\hline
  &   &
$
\left(\begin{array}{cccc}
(2\aaa_1)_I & (6\aaa_1)_{II}   & (6\aaa_1)_I     &  (6\aaa_1)_I     \\
            & 4\aaa_1          & (8\aaa_1)_{I}   &  (8\aaa_1)_I     \\
       &                       & 4\aaa_1         &   (8\aaa_1)_{II} \\
       &                       &                 &  4\aaa_1
\end{array}\right)
\subset 14\aaa_1
$\ $o$ & $19$   & $4_6^{+2},8_7^{+1}$ \\
\hline
  &   &
$
\left(\begin{array}{cccc}
(2\aaa_1)_{II} & 6\aaa_1     & 6\aaa_1          &  6\aaa_1        \\
               & 4\aaa_1     &(8\aaa_1)_I       &  (8\aaa_1)_I    \\
               &             & 4\aaa_1          &  (8\aaa_1)_{II} \\
       &                     &                  &  4\aaa_1
\end{array}\right)
\subset 14\aaa_1
$\ $o$ & $19$   & $4_4^{-2},8_5^{-1}$   \\
\hline
\end{tabular}
\end{table}

\begin{table}
\begin{tabular}{|c||c|c|c|c|c|c|c|c|}
\hline
 {\bf n}& $G$  &$Deg$& $\rk S$ &$q_S$            \\
\hline
\hline
  $10$ & $D_8$ &
$
\left(\begin{array}{cccc}
4\aaa_1       & (6\aaa_1)_I   & (8\aaa_1)_I  &  (8\aaa_1)_I \\
              & (2\aaa_1)_I   & 2\aaa_3      &  (6\aaa_1)_{II}   \\
       &                      &   4\aaa_1    &  4\aaa_2        \\
       &                      &              &  4\aaa_1
\end{array}\right)
\subset 4\aaa_1\amalg 2\aaa_5
$
& $19$   & $2_6^{+2},8_7^{+1}$  \\
\hline
  &   &
$
\left(\begin{array}{cccc}
(2\aaa_1)_{I}    & (6\aaa_1)_I     & (6\aaa_1)_{II}   & 10\aaa_1  \\
                 & 4\aaa_1         & (8\aaa_1)_I      &  12\aaa_1 \\
       &          & 4\aaa_1        &  4\aaa_3    \\
       &          &                &  8\aaa_1
\end{array}\right)
\subset 6\aaa_1\amalg 4\aaa_3
$ & $19$   & $2_{2}^{-2},8_3^{-1}$ $\ast$ \\
\hline
  &   &
$
\left(\begin{array}{cccc}
(2\aaa_1)_I & 2\aaa_3          & (6\aaa_1)_I     &  10\aaa_1  \\
            & 4\aaa_1          & (8\aaa_1)_{I}   &  12\aaa_1  \\
       &                       & 4\aaa_1         &  4\aaa_3   \\
       &                       &                 &  8\aaa_1
\end{array}\right)
\subset 6\aaa_3
$ \ $o$ & $19$   & $2_4^{-2},4_5^{-1}$ $\ast$ \\
\hline
  &   &
$
\left(\begin{array}{cccc}
(2\aaa_1)_{I}  & (6\aaa_1)_I & (6\aaa_1)_I      &  10\aaa_1 \\
               & 4\aaa_1     & 4\aaa_2          &  12\aaa_1 \\
               &             & 4\aaa_1          &  4\aaa_3  \\
       &                     &                  &  8\aaa_1
\end{array}\right)
\subset 2\aaa_1\amalg 4\ddd_4
$ & $19$   & $2_6^{+2},4_7^{+1}$ $\ast$  \\
\hline
\end{tabular}
\end{table}

%%%%%%%%%%%%%%%%%%%%%%%%%%%%%%%%%%%%%%%%%%%%%%%%%%%%%%%
%%%%%%%%%%%%%%%%%%%%%%%%%%%%%%%%%%%%%%%%%%%%%%%%%%%%%%
%%%%%%%%%%%%%%%%%%%%%%%%%%%%%%%%%%%%%%%%%%%%%%%%%%%%%%%%%
%%%%%%%%%%%%%%%%%%%%%%%%%%%%%%%%%%%%%%%%%%%%%%%%%%%%%%%%%
%%%%%%%%%%%%%%%%%%%%%%%%%%%%%%%%%%%%%%%%%%%%%%%%%%%%%%%%

\begin{table}
\label{table4}
\caption{Types and lattices $S$ of degenerations of
codimension $\ge 2$ of K\"ahlerian K3 surfaces
with symplectic automorphism group $(C_2)^3$.}

\begin{tabular}{|c||c|c|c|c|c|c|c|c|}
\hline
 {\bf n}& $|G|$ &  $G$   &$Deg$& $\rk S$ &$q_S$            \\
\hline
\hline
  $9$&$8$& $(C_2)^3$ &
  $((2\aaa_1,2\aaa_1)\subset 4\aaa_1)_I$ & $16$  &
$2_{II}^{+2},4_0^{+4}$ $\ast$\\
\hline
           &       &           &
$((2\aaa_1,2\aaa_1)\subset 4\aaa_1)_{II}$\ $o$ & $16$   &
$2_{II}^{+4},4_{II}^{+2}$ $\ast$ \\
\hline
    &       &           &
$(2\aaa_1,4\aaa_1)\subset 6\aaa_1$ & $16$   &
$2_{II}^{+4},4_7^{+1},8_1^{+1}$ $\ast$\\
\hline
    &       &           &
$(2\aaa_1,4\aaa_1)\subset 2\aaa_3$ & $16$   &
$2_{II}^{+4},4_{II}^{+2}$ $\ast$\\
\hline
    &       &           &
$(2\aaa_1,8\aaa_1)\subset 10\aaa_1$ & $16$   &
$2_{II}^{-4},4_4^{-2}$ $\ast$\\
\hline
    &       &           &
$(4\aaa_1,4\aaa_1)\subset 8\aaa_1$ & $16$   &
$2_{II}^{+4},4_{II}^{+2}$ \\
\hline
   &       &           &
$(4\aaa_1,8\aaa_1)\subset 4\aaa_3$ & $16$   &
$2_{II}^{+6}$ $\ast$\\
\hline
\end{tabular}
\end{table}

\begin{table}
\begin{tabular}{|c||c|c|c|c|c|c|c|c|}
\hline
 {\bf n}& $G$   &$Deg$& $\rk S$ &$q_S$            \\
\hline
\hline
  $9$& $(C_2)^3$ &
$
\left(\begin{array}{cccc}
 2\aaa_1 & (4\aaa_1)_I & (4\aaa_1)_I  & 10\aaa_1 \\
         & 2\aaa_1     & (4\aaa_1)_I  & 10\aaa_1 \\
         &             & 2\aaa_1      & 10\aaa_1 \\
         &             &              & 8\aaa_1
\end{array}\right)
\subset 14\aaa_1
$ & $18$   & $4_6^{+4}$ \\
\hline
    &       &
$((2\aaa_1,2\aaa_1)_I,4\aaa_1,4\aaa_1)\subset 12\aaa_1$ & $18$   &
$4_6^{+4}$ $\ast$ \\
\hline
    &       &
$((2\aaa_1,2\aaa_1)_{II},4\aaa_1,4\aaa_1)\subset 12\aaa_1$\ $o$  & $18$   &
$2_{II}^{-2},4_2^{-2}$ $\ast$ \\
\hline
     &       &
$
\left(\begin{array}{cccc}
 2\aaa_1 & 2\aaa_3     & (4\aaa_1)_I & 6\aaa_1 \\
         & 4\aaa_1     & 6\aaa_1  & 8\aaa_1 \\
         &             & 2\aaa_1      & 6\aaa_1 \\
         &             &              & 4\aaa_1
\end{array}\right)
\subset 2\aaa_3\amalg 6\aaa_1
$ & $18$   & $2_{II}^{-2},4_1^{+1},8_5^{-1}$ $\ast$ \\
\hline
     &       &
$
\left(\begin{array}{cccc}
 2\aaa_1 & 2\aaa_3     & (4\aaa_1)_I & 6\aaa_1 \\
         & 4\aaa_1     & 6\aaa_1     & 8\aaa_1 \\
         &             & 2\aaa_1     & 2\aaa_3 \\
         &             &              & 4\aaa_1
\end{array}\right)
\subset 4\aaa_3
$\ $o$ & $18$   & $2_{II}^{-2},4_2^{-2}$ $\ast$ \\
\hline
    &       &
$((2\aaa_1,2\aaa_1)_I,4\aaa_1,8\aaa_1)
\subset 4\aaa_1\amalg 4\aaa_3$ & $18$   &
$2_{II}^{-2},4_2^{-2}$ $\ast$ \\
\hline
    &       &
$(2\aaa_1,4\aaa_1,4\aaa_1,4\aaa_1)\subset 14\aaa_1$ & $18$   &
$2_{II}^{+2},4_5^{-1},8_5^{-1}$ $\ast$ \\
\hline

\end{tabular}
\end{table}

\newpage

\begin{table}
\begin{tabular}{|c||c|c|c|c|c|c|c|c|}
\hline
 {\bf n}& $G$   &$Deg$& $\rk S$ &$q_S$            \\
\hline
\hline
  $9$& $(C_2)^3$ &
$
\left(\begin{array}{ccccc}
 2\aaa_1 & (4\aaa_1)_I & (4\aaa_1)_I    & 6\aaa_1     &  10\aaa_1 \\
         & 2\aaa_1     & (4\aaa_1)_I    & 6\aaa_1     &  10\aaa_1 \\
         &             & 2\aaa_1        & 6\aaa_1     &  10\aaa_1 \\
         &             &                & 4\aaa_1     &   4\aaa_3\\
         &             &              &               &   8\aaa_1
\end{array}\right)
\subset 6\aaa_1\amalg 4\aaa_3
$ & $19$   & $4_5^{+3}$ $\ast$ \\
\hline
      &       &
$((2\aaa_1,2\aaa_1)_{II},4\aaa_1,4\aaa_1,4\aaa_1)\subset 16\aaa_1$ $o$ & $19$   &
$2_{II}^{-2},8_5^{-1}$ $\ast$ \\
\hline
     &       &
$
\left(\begin{array}{ccccc}
 2\aaa_1 & 2\aaa_3     & (4\aaa_1)_I  &  6\aaa_1     &  6\aaa_1 \\
         & 4\aaa_1     & 6\aaa_1     &  8\aaa_1     &  8\aaa_1 \\
         &             & 2\aaa_1      & 2\aaa_3     &  6\aaa_1 \\
         &             &              & 4\aaa_1     &  8\aaa_1\\
         &             &              &             &  4\aaa_1
\end{array}\right)
\subset 4\aaa_3\amalg  4\aaa_1
$\ $o$ & $19$   & $2_{II}^{-2},8_5^{-1}$ $\ast$ \\
\hline

\end{tabular}
\end{table}

V.V. Nikulin
\par Steklov Mathematical Institute,
\par ul. Gubkina 8, Moscow 117966, GSP-1, Russia;

\vskip5pt

Deptm. of Pure Mathem. The University of
Liverpool, Liverpool\par L69 3BX, UK
\par

\vskip5pt

nikulin@mi.ras.ru\, \ \
vnikulin@liv.ac.uk \, \ \  vvnikulin@list.ru
%%%%%%%%%%%%%%%%

Personal page: http://vnikulin.com

\end{document}